\documentclass[12pt,leqno]{amsart}

\usepackage{amssymb}
 \usepackage{amsthm}

\usepackage{amsmath,amsfonts,amsbsy,latexsym,multicol,fancybox}
\usepackage{graphicx}
\usepackage{epstopdf}
\usepackage{color}
\DeclareGraphicsExtensions{.pdf, .jpg, .pnp, .bmp, .fig}
\usepackage{hyperref}

\usepackage{caption}\usepackage{subcaption}

\newtheorem{theorem}{Theorem}
\newtheorem{assumption}{Assumption}
\newtheorem{definition}{Definition}
\newtheorem{lemma}{Lemma}
\newtheorem{remark}{Remark}

\newcommand{\bass}{\begin{assumption}}\newcommand{\eass}{\end{assumption}}
\newcommand{\bde}{\begin{definition}} \newcommand{\ede}{\end{definition}}
\newcommand{\ble}{\begin{lemma}} \newcommand{\ele}{\end{lemma}}
\newcommand{\bth}{\begin{theorem}} \newcommand{\ethe}{\end{theorem}}
\newcommand{\bre}{\begin{remark}} \newcommand{\ere}{\end{remark}}

\newcommand{\bpf}{\begin{proof}}\newcommand{\epf}{\end{proof}}
\newcommand{\barr}{\begin{array}}\newcommand{\earr}{\end{array}}
\newcommand{\beao}{\begin{eqnarray*}}\newcommand{\eeao}{\end{eqnarray*}\noindent}
\newcommand{\beam}{\begin{eqnarray}}\newcommand{\eeam}{\end{eqnarray}\noindent}
\newcommand{\beqq}{\begin{equation}}\newcommand{\eeqq}{\end{equation}\noindent}

 \newcommand{\un}{\underbrace}
\newcommand{\wt}{\widetilde}

\newcommand{\D}{\Delta}
\newcommand{\vep}{\varepsilon}

\newcommand{\bfE}{{\mathbb E}} 
\newcommand{\bbf}{{\mathcal F}}

 \newcommand{\bbN}{{\mathbb N}}

 \newcommand{\bbR}{{\mathbb R}}

 \usepackage{lineno}

\graphicspath{{./Images/}}

\begin{document}




\title[Lamperti Semi-Discrete method]{Lamperti Semi-Discrete method}
 \author{N. Halidias}
 \address{University of the Aegean, Department of Statistics and Actuarial-Financial Mathematics}
 \email{nick@aegean.gr}
 \author{I. S. Stamatiou}
 \address{University of West Attica, Department of Biomedical Sciences}
 \email{joniou@gmail.com, istamatiou@uniwa.gr}

\begin{abstract}
We study the numerical approximation of numerous processes, solutions of nonlinear stochastic differential equations, that appear in various applications such as financial mathematics and population dynamics. Between the investigated models are the CIR process, also known as the square root process, the constant elasticity of variance process CEV, the Heston $3/2$-model, the A\"it-Sahalia model and the Wright-Fisher model.  We propose a version of the semi-discrete method, see \cite{halidias_stamatiou:2020}, which we call Lamperti semi-discrete (LSD) method. The LSD method is domain preserving and seems to converge strongly to the solution process with order $1$ and no extra restrictions on the parameters or the step-size.
\end{abstract}

\date\today
\maketitle
\tableofcontents
\listoffigures

\section{CIR model}\label{LSD:ssec:CIR}
Let
\beqq  \label{LSD-eq:exampleSDE}
x_t =x_0  + \int_0^t (k_1 - k_2x_s) ds + \int_0^t k_3\sqrt{x_s} dW_s, \quad t\geq0.
\eeqq

SDE (\ref{LSD-eq:exampleSDE}) is known as the CIR process, or square root process with a solution remaining in the positive axis, i.e. $x_t>0$ a.s. when $(k_3)^2\leq2k_1$ and $x_0>0,$ c.f. \cite{rogers_williams_2000}. The Lamperti transformation of (\ref{LSD-eq:exampleSDE}) is $z = \frac{2}{k_3}\sqrt{x}$ and application of the  It\^o formula implies the following representation, see Appendix \ref{LSD-ap:Lamperti_tranformation},  
\beqq  \label{LSD-eq:exampleSDELamperti}
z_t =z_0 + \int_0^t \left(\frac{2k_1}{(k_3)^2}(z_s)^{-1} - \left(\frac{k_2}{2} + \frac{(k_3)^2}{8}\right)z_s\right)ds + \int_0^t  dW_s, \,\, t\geq0.
\eeqq 

To simplify notation, set $a:=\frac{2k_1}{(k_3)^2}$ and $b:=\frac{k_2}{2} + \frac{(k_3)^2}{8}.$ The new coefficients $a$ and $b$ are positive. We consider three versions of the semi-discrete method for approximating  (\ref{LSD-eq:exampleSDELamperti}). In the first two versions $(y_t)$ and $(\hat{y}_t),$ see Section \ref{LSD:subsec:CIRz1z2}, we use the semi-discrete method as originally proposed, see \cite{halidias:2012}; we discretize parts 
of the drift  coefficient producing a new differential equation in each subinterval  with known solution or a solution that is easily simulated or approximated. In the third version $(\breve{y}_t)$ we examine a new modification of the semi-discrete method, where in each subinterval $(t_n, t_{n+1}]$ we do not need to solve a new differential equation, but only an algebraic equation.

\subsection{Lamperti Semi-Discrete methods $\wt{z}^1_{n}$ and $\wt{z}^2_{n}$ for CIR}\label{LSD:subsec:CIRz1z2}

Rewrite (\ref{LSD-eq:exampleSDELamperti}) as 
\beqq  \label{LSD-eq:exampleSDELamperti_ab}
z_t =z_0 + \int_0^t \left( a(z_s)^{-1} - bz_s\right)ds + \int_0^t  dW_s, \,\, t\geq0.
\eeqq

To approximate the constant diffusion SDE (\ref{LSD-eq:exampleSDELamperti_ab}) we use the following two versions of the semi-discrete method. In the first version $(y_t)$, we discretize the linear part of the drift  coefficient and in the second $(\hat{y}_t)$ we leave it as it since it produces a differential equation with known solution.  
Let $t\in(t_n, t_{n+1}],$ where we assume the length of each subinterval to be equal to $\D$ and consider 
\beqq\label{LSD-eq:SD schemeExampleLT}
y_t = \D W_n + y_{t_n} - by_{t_n}\D +  \int_{t_n}^t a(y_s)^{-1}ds,
\eeqq
with $y_{0}=z_0$ and  
\beqq\label{LSD-eq:SD schemeExampleLT2}
\hat{y}_t = \D W_n + \hat{y}_{t_n} + \int_{t_n}^t \left(a(\hat{y}_s)^{-1} - b\hat{y}_s\right)ds,
\eeqq
with $\hat{y}_{0}=z_0$  

(\ref{LSD-eq:SD schemeExampleLT}) and (\ref{LSD-eq:SD schemeExampleLT2}) are Bernoulli type equations with solutions satisfying, see Appendix \ref{LSD-ap:Bernoulli_sol},
\beqq\label{LSD-eq:SD schemeExampleLTsol}
(y_t)^2 = \left(\D W_n +  (1-b\D )y_{t_n}\right)^2 + 2a(t-t_n)
\eeqq
and
\beqq\label{LSD-eq:SD schemeExampleLTsol2}
(\hat{y}_t)^2 = (\D W_n + \hat{y}_{t_n})^2e^{-2b(t-t_n)} + a\frac{1-e^{-2b(t-t_n)}}{b},
\eeqq
respectively.

We propose the following versions of the semi-discrete method for the approximation of (\ref{LSD-eq:exampleSDELamperti}),  
\beqq\label{LSD-eq:SD schemeExampleLT_transf}
y_{t_{n+1}} = \sqrt{\left(\D W_n +  (1-b\D)y_{t_n}\right)^2 + 2a\D}
\eeqq
and
\beqq\label{LSD-eq:SD schemeExampleLT_transf2}
\hat{y}_{t_{n+1}} = \sqrt{(\D W_n + \hat{y}_{t_n})^2e^{-2b\D} + a\frac{1-e^{-2b\D}}{b}},
\eeqq
which suggests the versions of the Lamperti semi-discrete method $(\wt{z}^1_n)_{n\in\bbN}, (\wt{z}^2_n)_{n\in\bbN}$ for the approximation of (\ref{LSD-eq:exampleSDE})
\beqq\label{LSD-eq:SD schemeExampleLToriginal}
\wt{z}^1_{t_{n+1}} = \frac{(k_3)^2}{4}\left(\left(\D W_n +  (1-b\D)y_{t_n}\right)^2 + 2a\D\right)
\eeqq
and
\beqq\label{LSD-eq:SD schemeExampleLToriginal2}
\wt{z}^2_{t_{n+1}} = \frac{(k_3)^2}{4}\left((\D W_n + \hat{y}_{t_n})^2e^{-2b\D} + a\frac{1-e^{-2b\D}}{b}\right).
\eeqq
 
\subsection{Lamperti Semi-Discrete method $\wt{z}^3_{n}$ for CIR}\label{LSD:subsec:CIRz3}

In this version of the semi-discrete method, $(\breve{y}_t)$ we examine a new modification of the semi-discrete method, where in each subinterval $(t_n, t_{n+1}]$ we do not need to solve a new differential equation, but only an algebraic equation. For $t\in(t_n, t_{n+1}]$ consider 

\beqq\label{LSD-eq:SD schemeExampleLT3}
\breve{y}_t = W_t - W_{t_n} + \breve{y}_{t_n} + a(\breve{y}_t)^{-1}\D  - b\breve{y}_t\D.
\eeqq
with $\breve{y}_{0}=z_0.$ The solution of (\ref{LSD-eq:SD schemeExampleLT3}) satisfies 

\beqq\label{LSD-eq:SD schemeExampleLTsol3}
(1 + b\D)(\breve{y}_t)^2 - (W_t - W_{t_n} + \breve{y}_{t_n})\breve{y}_t -a\D = 0.  
\eeqq

We propose the following version of the semi-discrete method for the approximation of (\ref{LSD-eq:exampleSDELamperti}), 

\beqq\label{LSD-eq:SD schemeExampleLT_transf3}
\breve{y}_{t_{n+1}} = \frac{\D W_n + \breve{y}_{t_n} +  \sqrt{(\D W_n + \breve{y}_{t_n})^2 +  4(1+b\D)a\D}}{2(1+b\D)},
\eeqq
which suggests the version of the Lamperti semi-discrete method $(\wt{z}^3_n)_{n\in\bbN}$ for the approximation of (\ref{LSD-eq:exampleSDE})

\beqq\label{LSD-eq:SD schemeExampleLToriginal3}
\wt{z}^3_{t_{n+1}} = \frac{(k_3)^2}{4}\left(\frac{\D W_n + \breve{y}_{t_n} +  \sqrt{(\D W_n + \breve{y}_{t_n})^2 +  4(1+b\D)a\D}}{2(1+b\D)}\right)^2.
\eeqq

\subsection{Numerical experiment for CIR}\label{LSD:subsec:numCIR}

For a minimal numerical experiment we present simulation paths for the numerical approximation of (\ref{LSD-eq:exampleSDE}) with  $x_0=4$ and compare with the SD method proposed in \cite{halidias:2015b}, which reads 

\beqq\label{LSD-eq:NCEV_HAL_sdesol} 
\wt{y}_{t_{n+1}}=\left(\sqrt{\wt{y}_{t_n}(1- \frac{k_2\D}{1+k_2\theta\D}) + \frac{\D}{1 + k_2\theta\D}\left(k_1  - \frac{(k_3)^2}{4(1 + k_2\theta\D)}\right)} + \frac{k_3}{2(1 + k_2\theta\D)}\D W_n\right)^{2},
\eeqq
where $\theta$ represents the level of implicitness. The case $\theta = 0 $ was studied in \cite{halidias:2012} where the idea of the semi-discrete method was originally presented. According to the results in \cite{halidias:2015b}  it is shown that SD method (\ref{LSD-eq:NCEV_HAL_sdesol}) is strongly convergent under some conditions on the coefficients $k_i,$ the level of implicitness $\theta$ and the step-size $\D.$ In particular, it strongly converges to the solution of (\ref{LSD-eq:exampleSDE}) with a logarithmic rate if also $\bfE(x_0)^p<A$ for some $p\geq2, (k_3)^2\leq 4k_1(1 + k_2\theta\D)$ and $\D(1-\theta)\leq(k_2)^{-1},$ while a polynomial rate of convergence is achieved with order at least $1/4$ for a smaller set of parameters, namely $(k_3)^2\leq 2k_1$ and $(\frac{2k_1}{(k_3)^2}-1)^2>16$ for $\theta=0,$ with $x_0\in\bbR.$ On the other hand, the LSD scheme (\ref{LSD-eq:SD schemeExampleLToriginal}) seems to work without any restriction on the step-size or on the parameters which is a very interesting result.

\bre\label{LSD-rem:Ht_CIR}
We would like to point out a mistake that escaped our attention. In the proof of the strong convergence properties of the SD scheme (\ref{LSD-eq:NCEV_HAL_sdesol}) proposed in \cite{halidias:2015b} an auxiliary process $(h_t)$ appears, see \cite[Rel. (2.3)]{halidias:2015b}
$$
h_t = x_0 + \int_0^t (\un{k_1 - k_2(1-\theta)y_{\hat{s}} - k_2\theta y_{\tilde{s}}}_{f_{\theta}(y_{\hat{s}}, y_{\tilde{s}})})ds + \int_0^t k_3\sqrt{y_s}dW_s,
$$
where $\hat{s} = t_j$ when $s\in(t_j, t_{j+1}], j =0,1,\ldots, n$ and 
$$
\wt{s}=\left\{ \barr{ll}  t_{j+1},  & \mbox{for } \, s\in[t_j,t_{j+1}], \\  t,  & \mbox{for }\,  s\in[t_n,t] \earr j=0,\ldots,n-1.\right.
$$ 
The problem is that we can not apply directly the It\^o formula on $(h_t)$ since $y_{\tilde{s}}$ is $\bbf_{t_{n+1}}$-measurable and not  $\bbf_{t_{n}}$-measurable. Nevertheless, writing the drift of $(h_t)$ as  $f_{\theta}(y_{\hat{s}}, y_{\hat{s}})$ 
we can proceed in the same way. The remainder term $f_{\theta}(y_{\hat{s}}, y_{\tilde{s}}) - f_{\theta}(y_{\hat{s}}, y_{\hat{s}}) = k_2\theta (y_{\hat{s}} - y_{\tilde{s}})$ can be easily bounded. 
\ere

We also present the implicit scheme proposed in   \cite{alfonsi:2005}, which takes the following form

\beqq\label{LSD-eq:ALF_sdesol} 
\bar{y}_{t_{n+1}}=\left( \frac{\sqrt{4(\bar{y}_{t_n}  +(k_1-\frac{(k_3)^2}{2})\D)(1 + k_2\D) + (k_3)^2(\D W_n)^2} + k_3\D W_n}{2(1 + k_2\D)}\right)^{2},
\eeqq

\begin{figure}[ht]
	\centering
	\begin{subfigure}{.47\textwidth}
 		\includegraphics[width=1\textwidth]{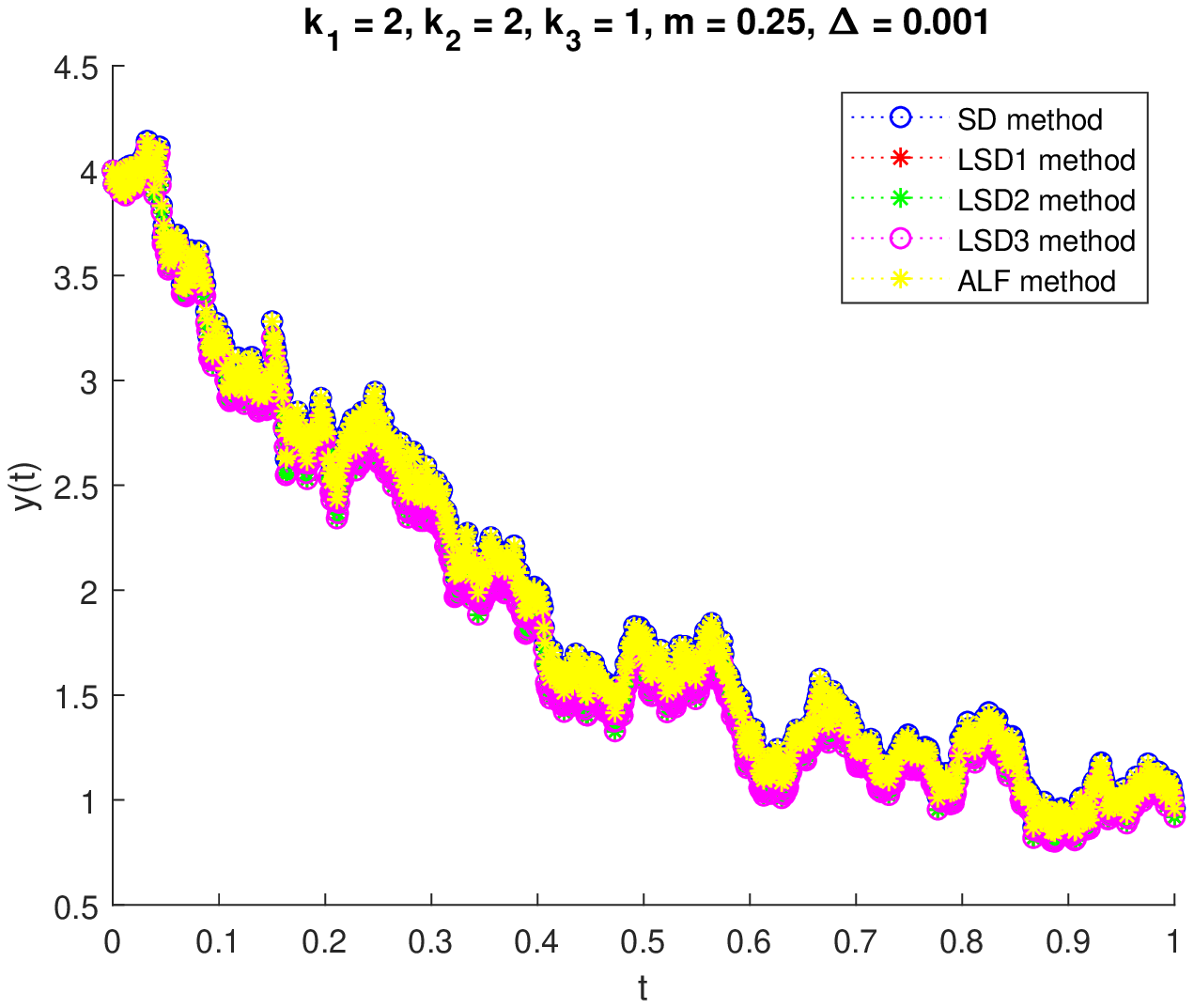}
		\caption{Trajectories  of  (\ref{LSD-eq:SD schemeExampleLToriginal})-(\ref{LSD-eq:ALF_sdesol}).}
	\end{subfigure}
	\begin{subfigure}{.47\textwidth}
		\includegraphics[width=1\textwidth]{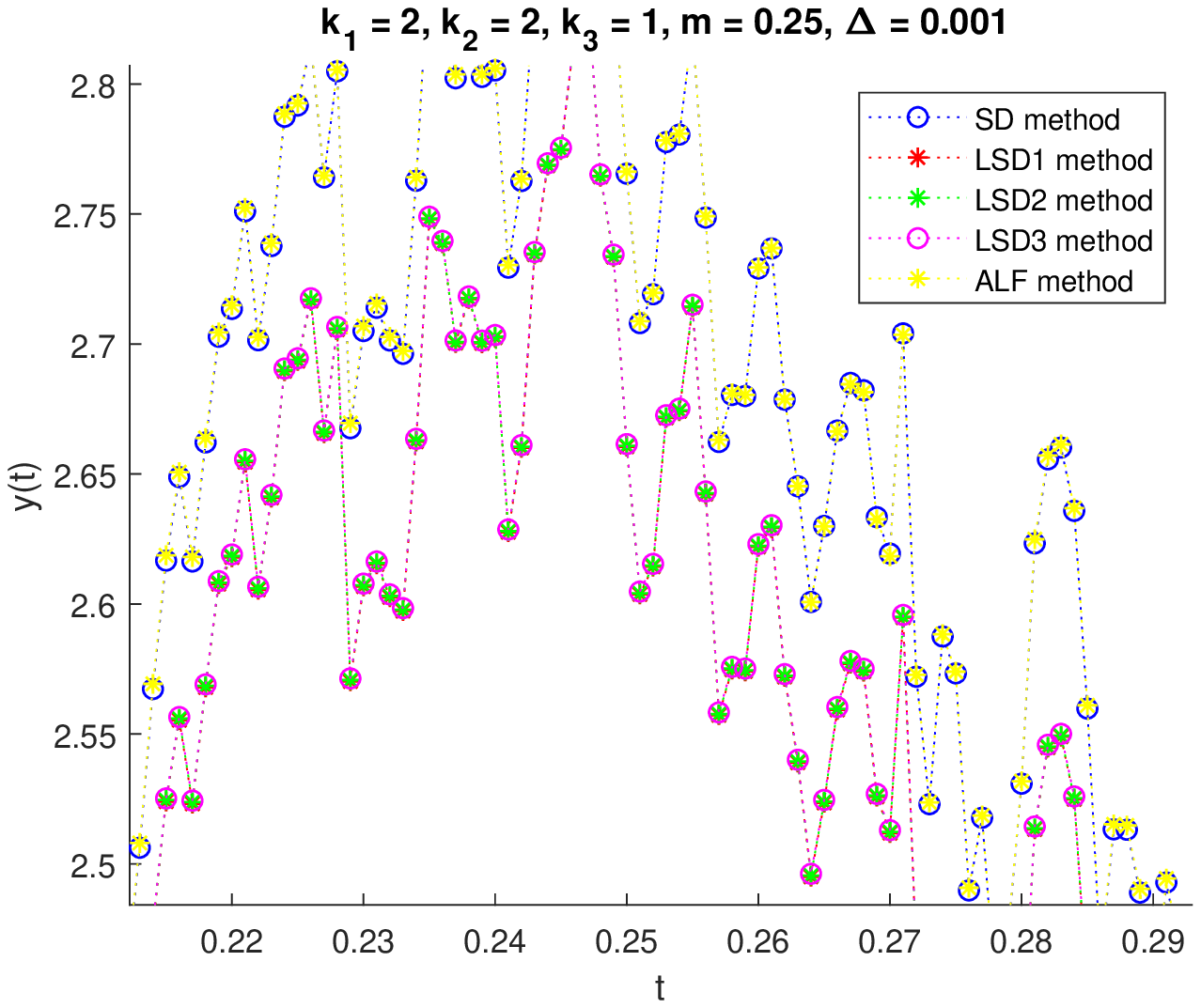}
		\caption{Zoom of Figure \ref{LSD-fig:LSDsSD}(A).}
	\end{subfigure}
	\caption{Trajectories  of  (\ref{LSD-eq:SD schemeExampleLToriginal}), (\ref{LSD-eq:SD schemeExampleLToriginal2}), (\ref{LSD-eq:SD schemeExampleLToriginal3}), (\ref{LSD-eq:NCEV_HAL_sdesol}) and (\ref{LSD-eq:ALF_sdesol}) for the approximation of (\ref{LSD-eq:exampleSDE}) with $\D=10^{-4}$.}\label{LSD-fig:LSDsSD}
\end{figure}

\begin{figure}[ht]
	\centering
	\begin{subfigure}{.47\textwidth}
		\includegraphics[width=1\textwidth]{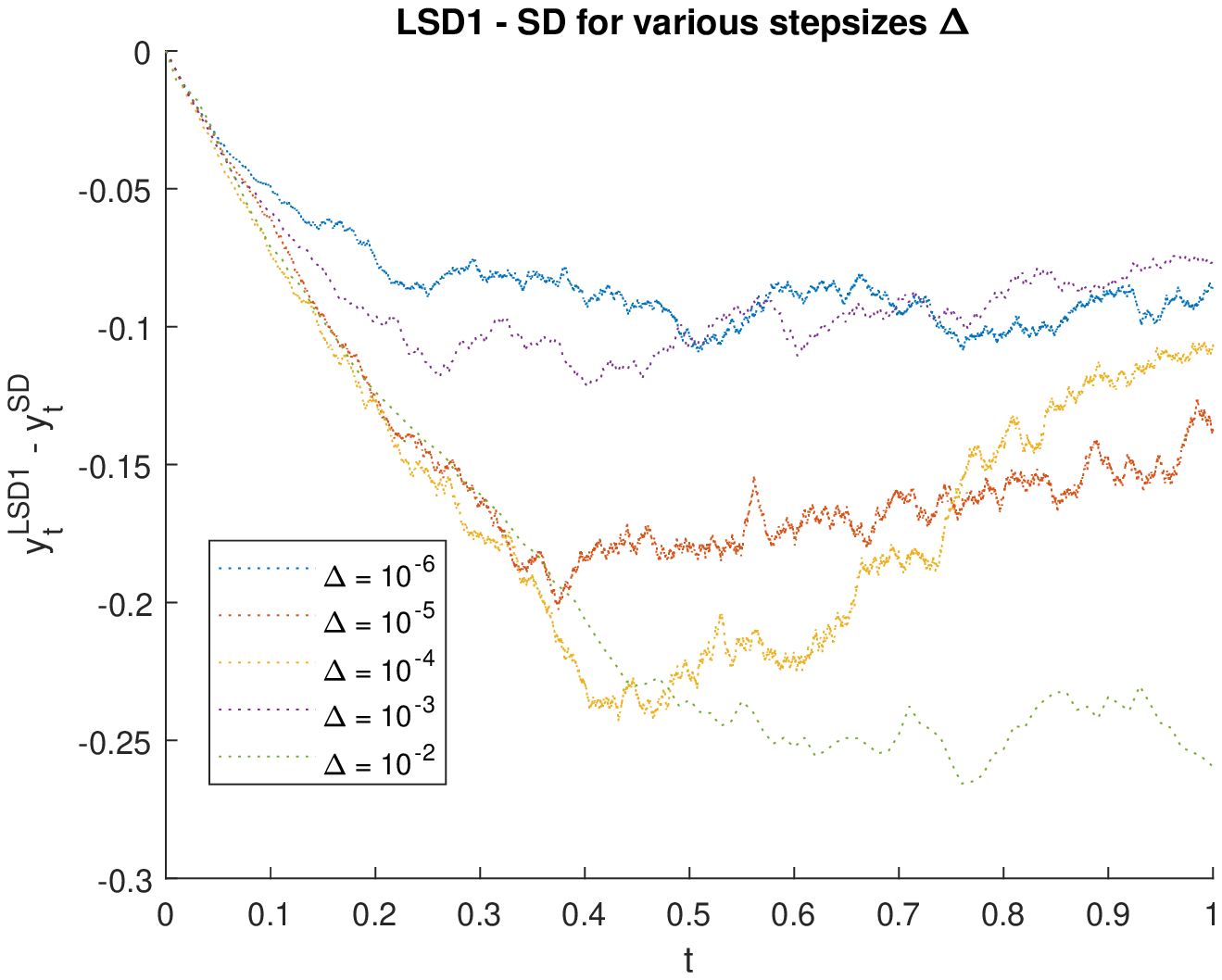}
		\caption{Difference (\ref{LSD-eq:SD schemeExampleLToriginal}) - (\ref{LSD-eq:NCEV_HAL_sdesol}) for various step sizes.}
	\end{subfigure}
	\begin{subfigure}{.47\textwidth}
		 \centering
		\includegraphics[width=1\textwidth]{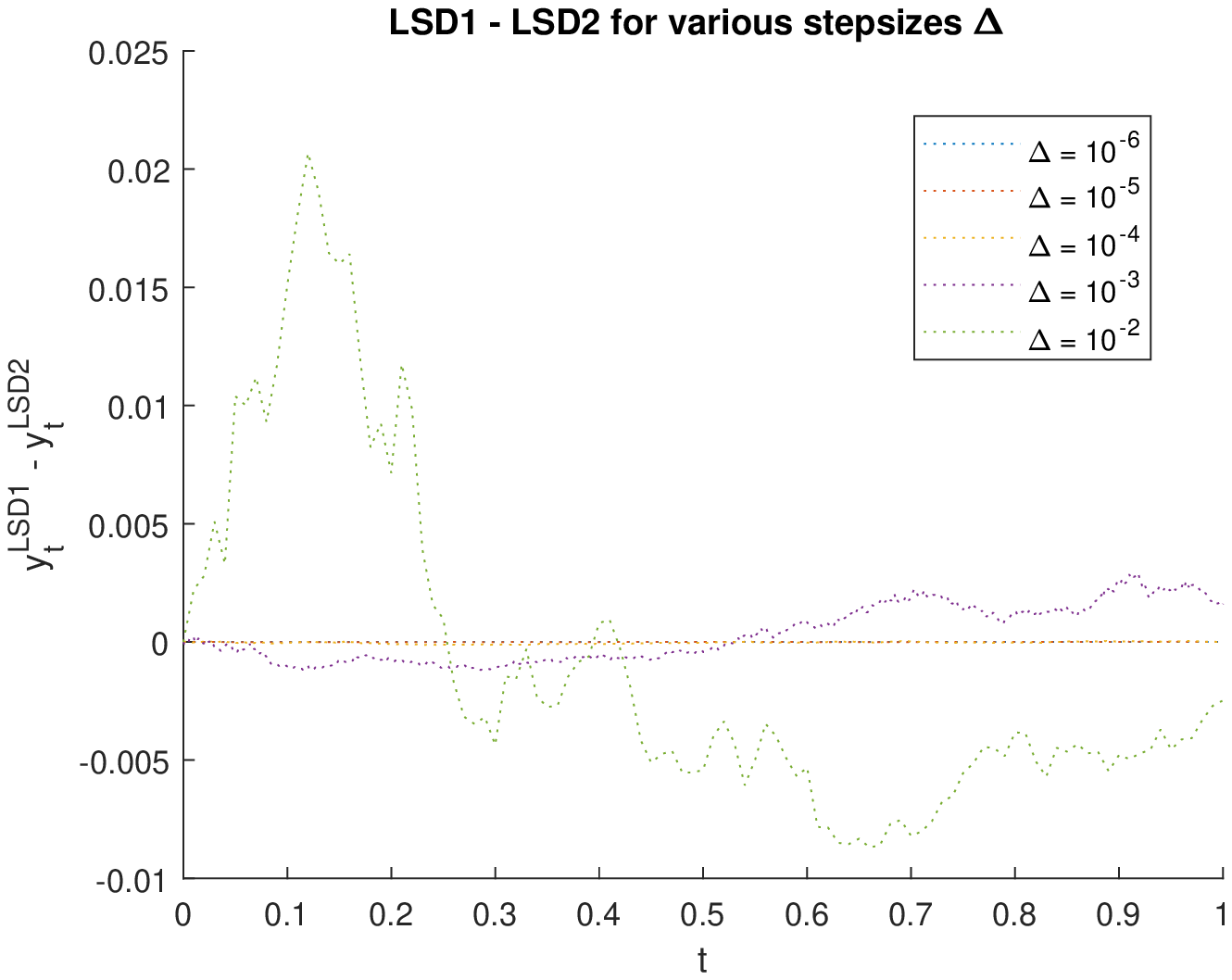}
		\caption{Difference (\ref{LSD-eq:SD schemeExampleLToriginal}) - (\ref{LSD-eq:SD schemeExampleLToriginal2}) for various step sizes.}
	\end{subfigure}
\begin{subfigure}{.55\textwidth}
	\centering
	\includegraphics[width=1\textwidth]{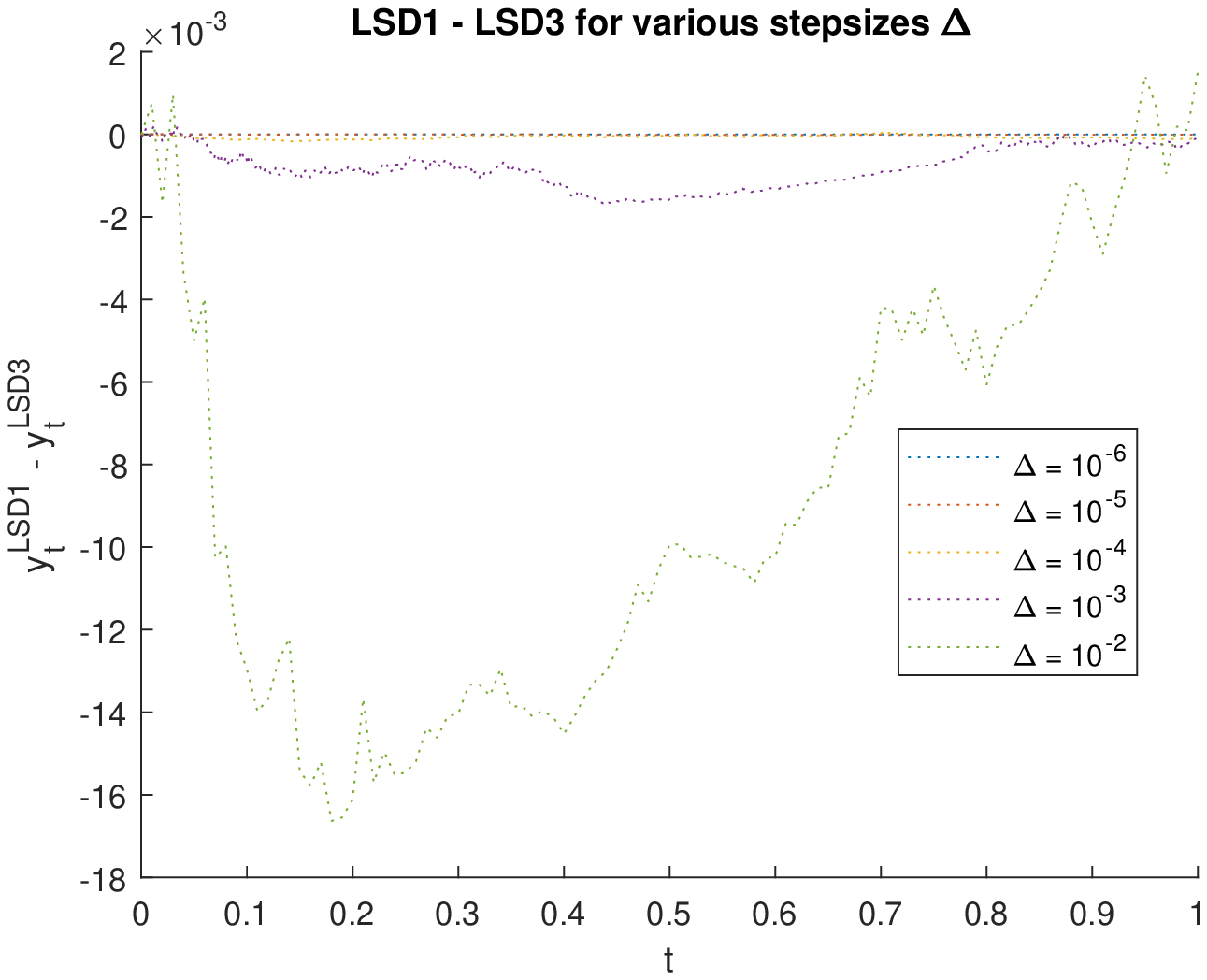}
	\caption{Difference (\ref{LSD-eq:SD schemeExampleLToriginal}) - (\ref{LSD-eq:SD schemeExampleLToriginal3}) for various step sizes.}
\end{subfigure}
	\caption{Trajectories of the differences of the semi-discrete methods (\ref{LSD-eq:SD schemeExampleLToriginal}), (\ref{LSD-eq:SD schemeExampleLToriginal2}) and  (\ref{LSD-eq:NCEV_HAL_sdesol}) for the approximation of (\ref{LSD-eq:exampleSDE}) with various step-sizes.}\label{LSD-fig:LSDsminSD}
\end{figure}

As a first graphical illustration we borrow the set of parameters from \cite[Sec.4]{halidias:2015b}; we take $k_1 = k_2 =2, k_3 = 1, T = 1$ and $\theta = 1$ with various step-sizes $\D = 10^{-4}$ and $\D = 10^{-3}, \D = 10^{-2}.$ We compare with the proposed two versions of LSD scheme (\ref{LSD-eq:SD schemeExampleLToriginal}) and (\ref{LSD-eq:SD schemeExampleLToriginal2}) and the implicit method ALF (\ref{LSD-eq:ALF_sdesol}). Figure \ref{LSD-fig:LSDsSD} shows that all the schemes perform in a similar way. We also give a presentation of the difference of the various SD approximations in Figure \ref{LSD-fig:LSDsminSD}.

\begin{figure}[ht]
	\centering
	\begin{subfigure}{.47\textwidth}
		\includegraphics[width=1\textwidth]{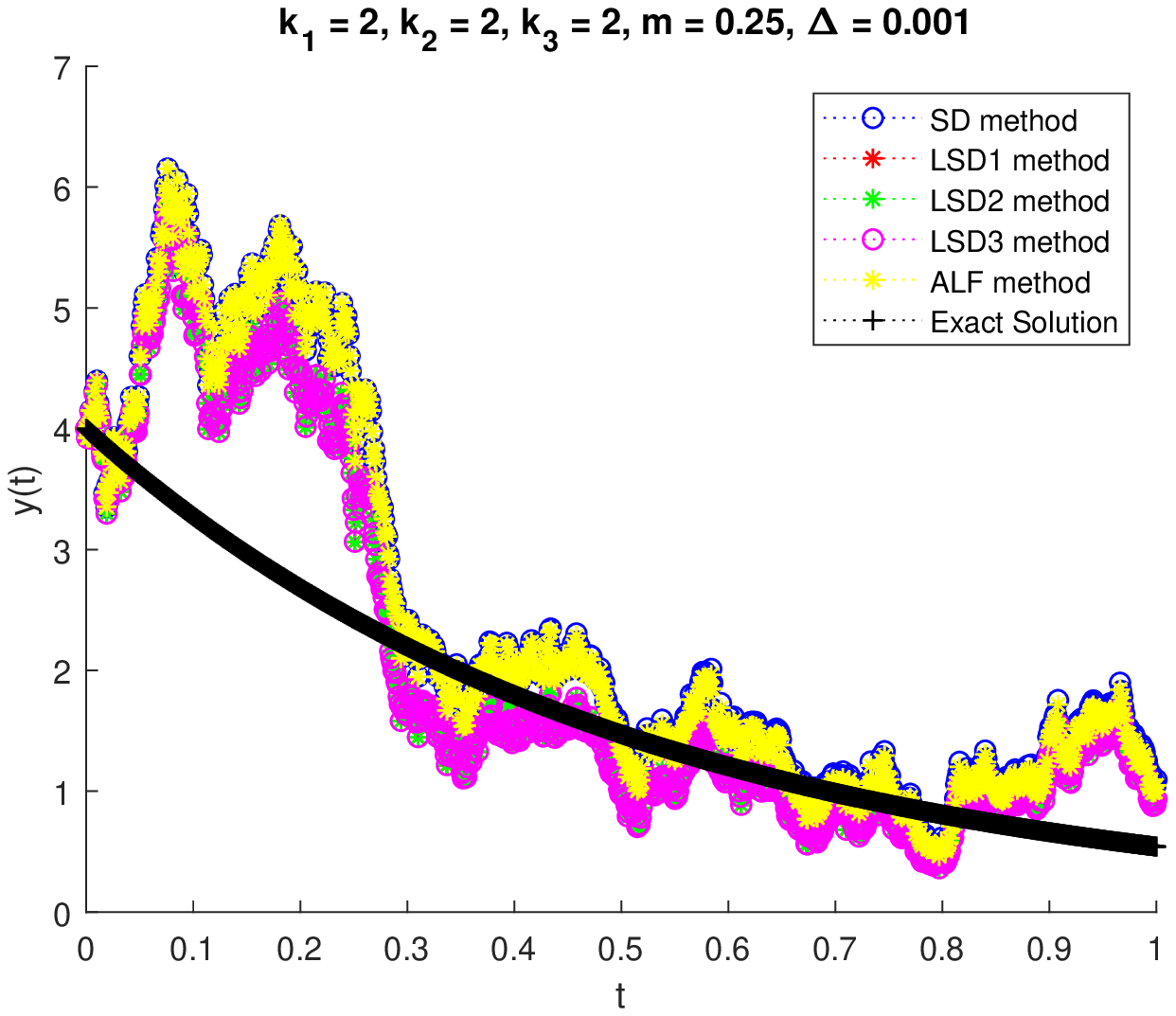}
	\end{subfigure}
	\begin{subfigure}{.47\textwidth}
		\centering
		\includegraphics[width=1\textwidth]{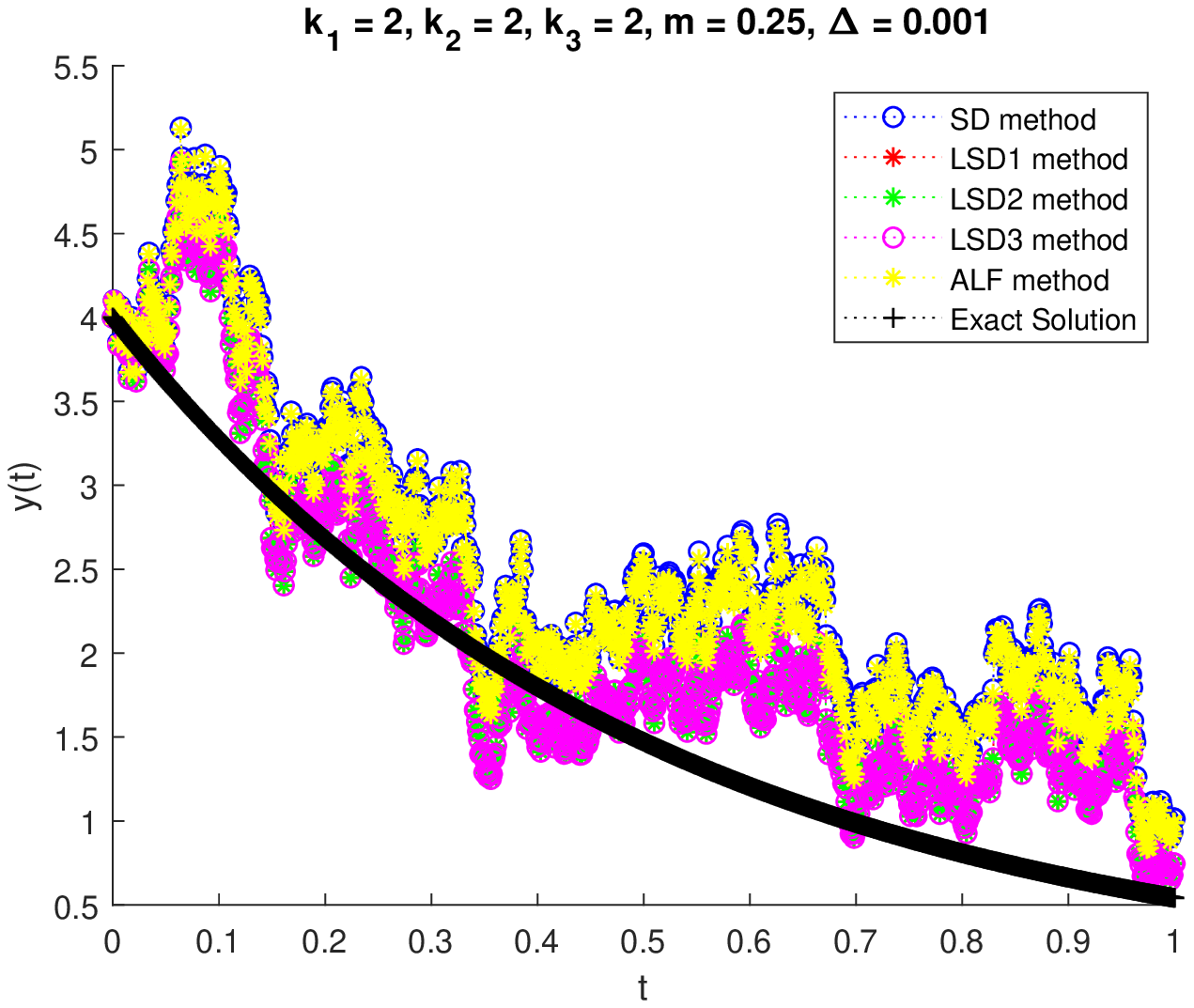}
	\end{subfigure}
	\begin{subfigure}{.47\textwidth}
	\centering
	\includegraphics[width=1\textwidth]{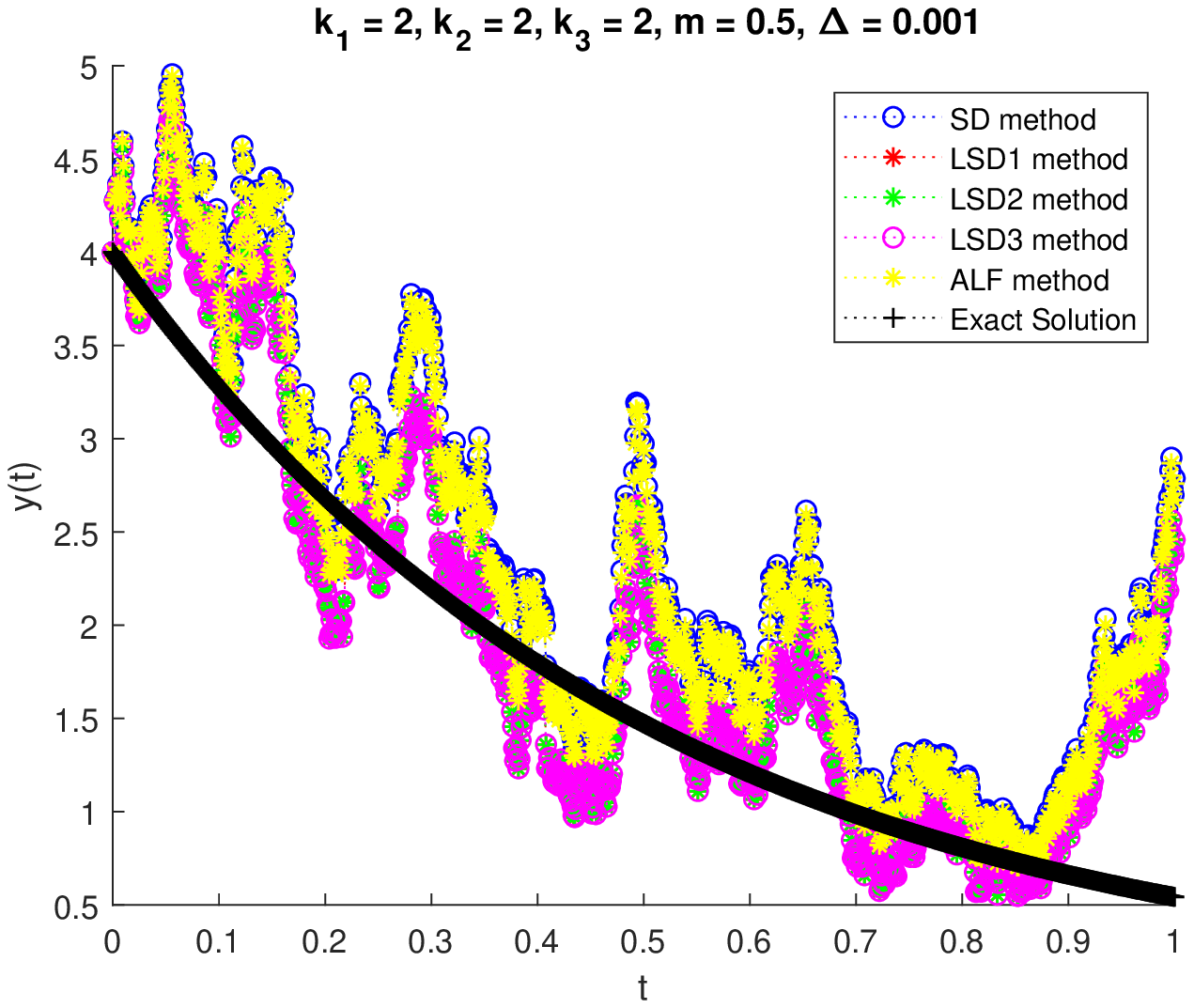}
\end{subfigure}
	\begin{subfigure}{.47\textwidth}
	\centering
	\includegraphics[width=1\textwidth]{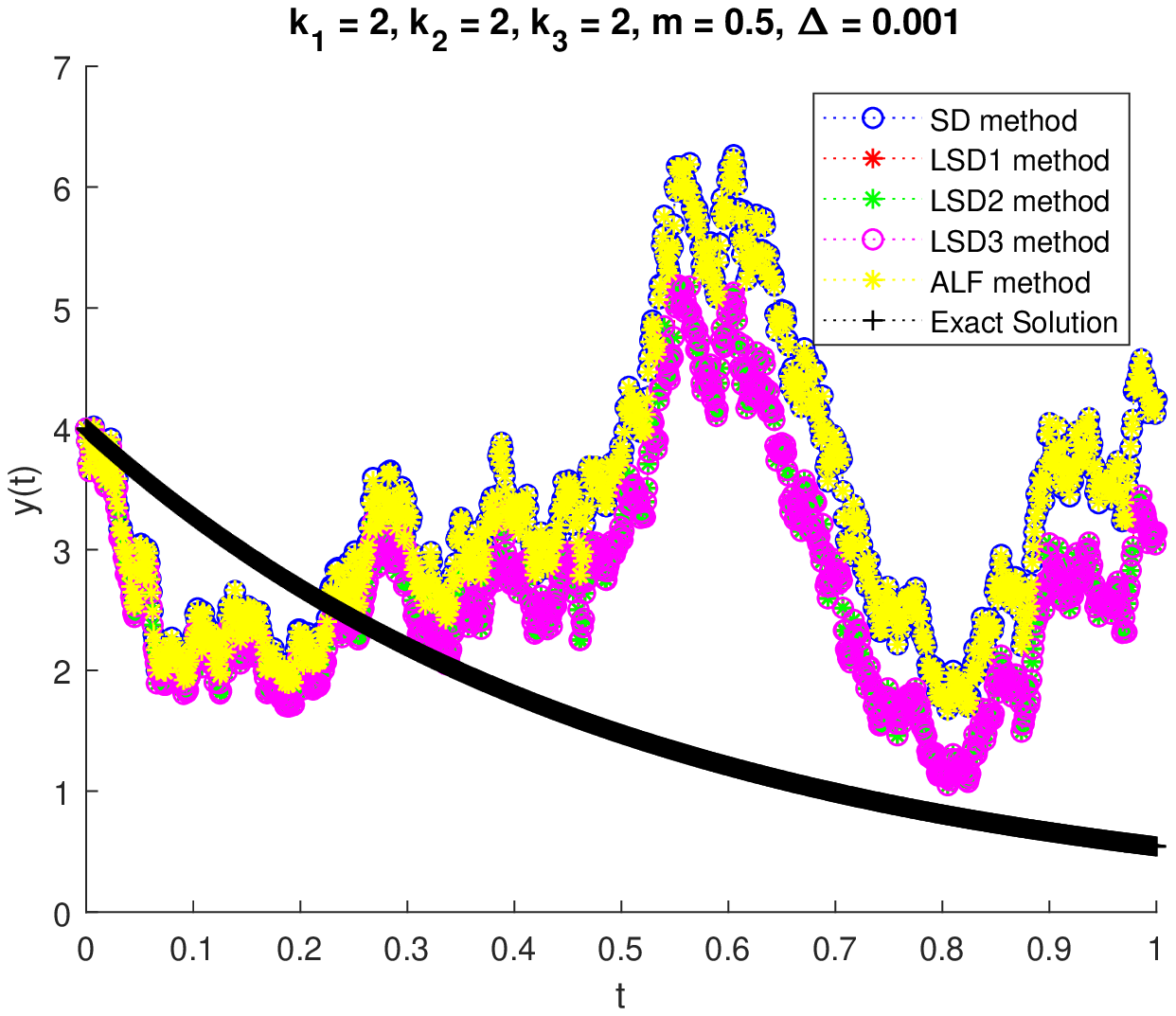}
\end{subfigure}
	\begin{subfigure}{.47\textwidth}
	\centering
	\includegraphics[width=1\textwidth]{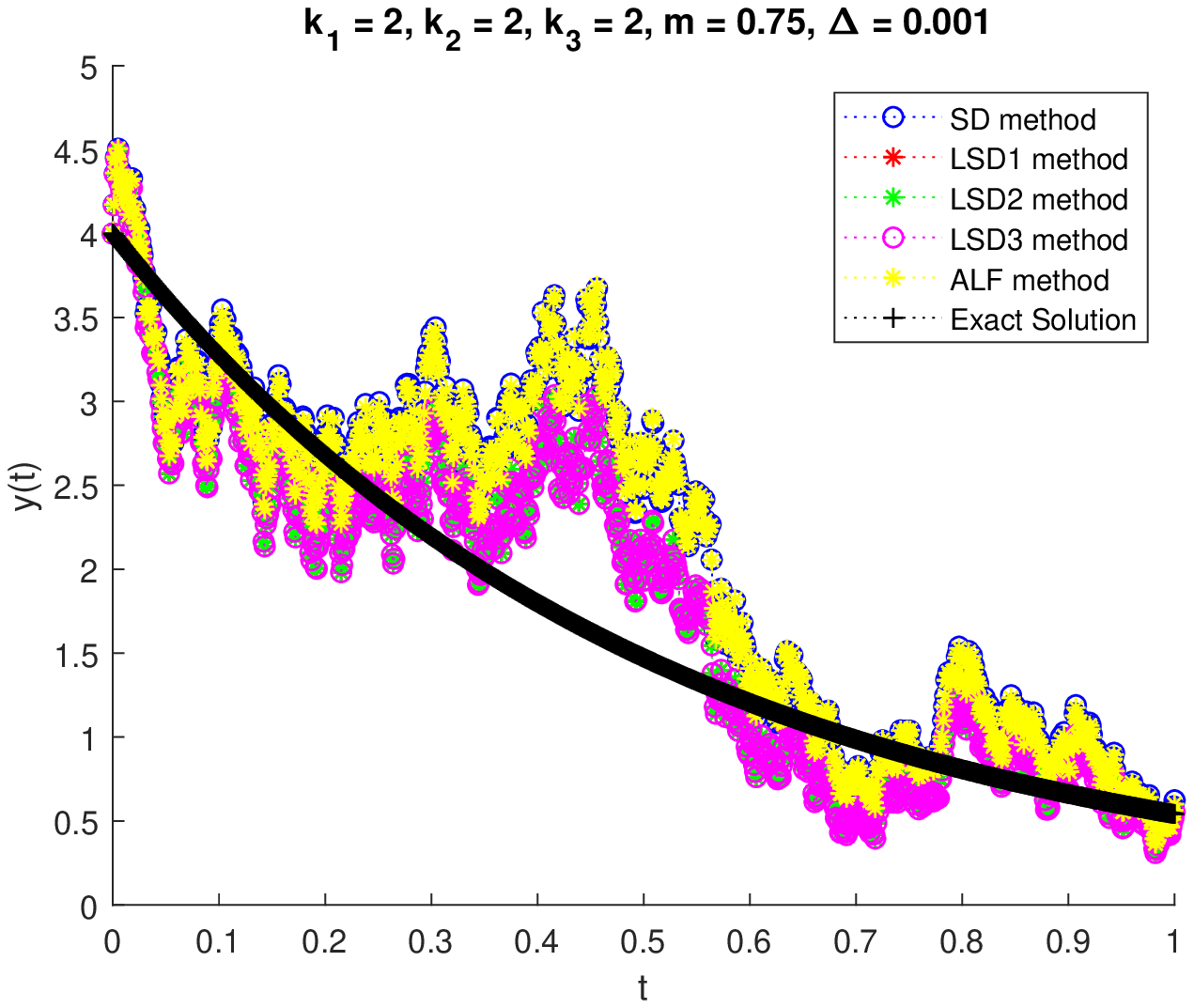}
\end{subfigure}
	\begin{subfigure}{.47\textwidth}
	\centering
	\includegraphics[width=1\textwidth]{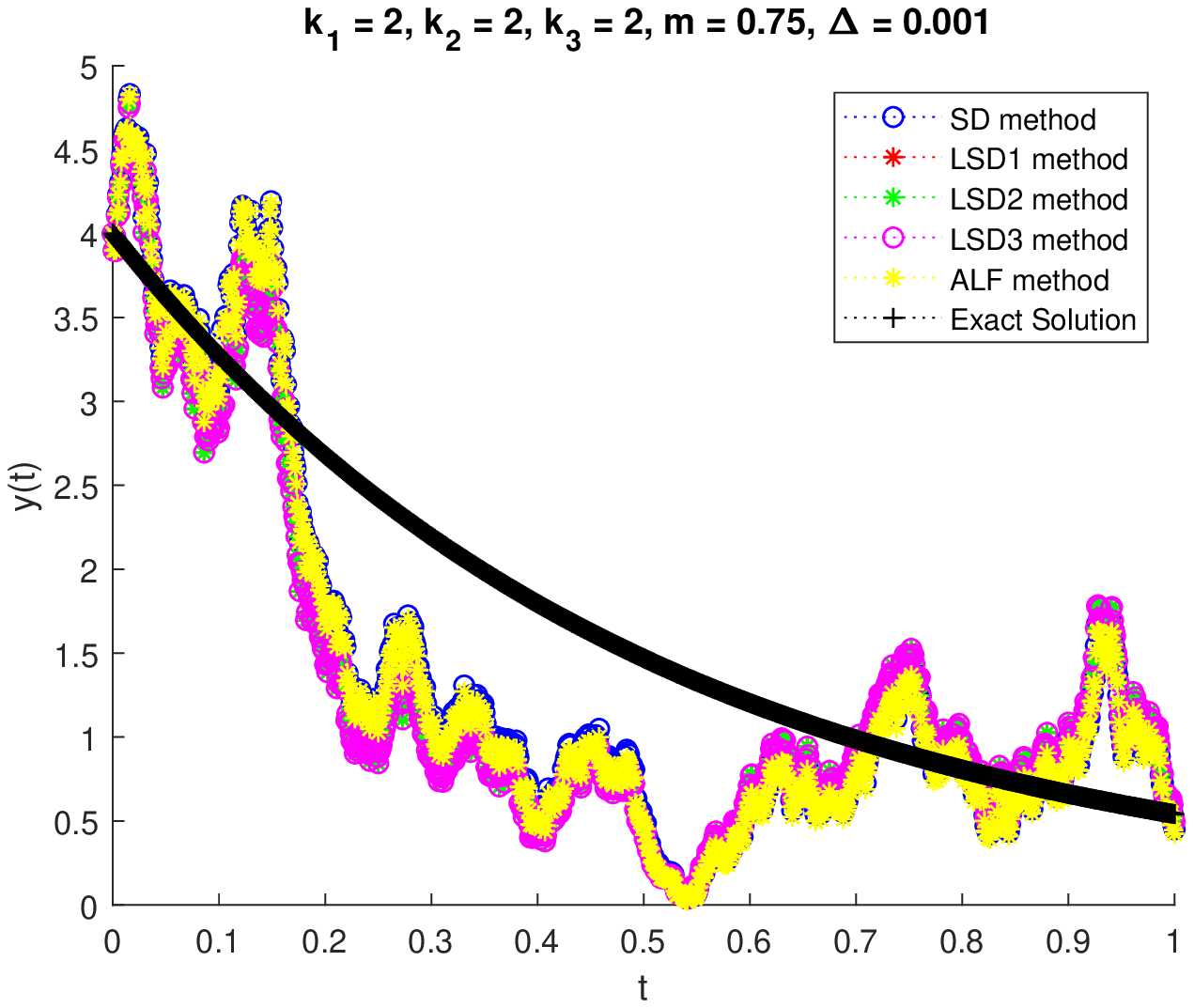}
\end{subfigure}
	\caption{Trajectories of (\ref{LSD-eq:SD schemeExampleLToriginal})-(\ref{LSD-eq:ALF_sdesol}) and the exact solution (\ref{LSD-eq:OUsolex}) with different $m$ and $\D=0.001.$}\label{LSD-fig:LSDvsExact}
\end{figure}

For a different configuration, we are able to compare with the exact solution. Note that by choosing $d = 4k_1/(k_3)^2 = 2$ then the solution of (\ref{LSD-eq:exampleSDE}) is $x_t = (x_1(t))^2 +  (x_2(t))^2,$ where 
$x_j(t)$ is the solution of the \textit{Orstein-Uhlenbelck} process
\beqq\label{LSD-eq:OU}
dx_j(t) = -\frac{1}{2}k_2x_j(t)dt + \frac{1}{2}k_3 dW^{(j)}_t
\eeqq
and the Brownian motions $W^{(j)}_t$ are independent. For $t\in[t_n,t_{n+1}]$ the solution of (\ref{LSD-eq:OU}) is 
\beqq\label{LSD-eq:OUsol}
x_j(t) = e^{-\frac{1}{2}k_2(t-t_n)}x_j(t_n) + \frac{1}{2}k_3e^{-\frac{1}{2}k_2(t-t_n)}\int_{t_n}^te^{\frac{1}{2}k_2(s-t_n)} dW^{(j)}_s.
\eeqq
Actually, we approximate in each subinterval $[t_n,t_{n+1}]$ the last stochastic integral in (\ref{LSD-eq:OUsol})  at $t_n$ producing the sequence
\beqq\label{LSD-eq:OUsolappr}
x_j(t_{n+1}) = e^{-\frac{1}{2}k_2\D}x_j(t_n) + \frac{k_3}{k_2}(1-e^{-\frac{1}{2}k_2\D})\D W^{(j)}_n
\eeqq
and therefore the solution process of (\ref{LSD-eq:exampleSDE})  at the grid points reads
\beam\nonumber
x(t_{n+1})&  = & \left(e^{-\frac{1}{2}k_2\D}x_1(t_n) + \frac{k_3}{k_2}(1-e^{-\frac{1}{2}k_2\D})\D W^{(1)}_n\right)^2\\
\label{LSD-eq:OUsolex}
&&+  \left(e^{-\frac{1}{2}k_2\D}x_2(t_n) + \frac{k_3}{k_2}(1-e^{-\frac{1}{2}k_2\D})\D W^{(2)}_n\right)^2.
\eeam
We therefore choose $k_3 = 2$ for this second experiment, with all the other parameters unchanged, so that $d = 2.$ Moreover $x(0) = (x_1(0))^2+ (x_2(0))^2.$ 
We present  in Figure \ref{LSD-fig:LSDvsExact} simulation paths of (\ref{LSD-eq:SD schemeExampleLToriginal})-(\ref{LSD-eq:ALF_sdesol}) and the exact solution (\ref{LSD-eq:OUsolappr}) choosing as initial conditions $x_1(0) = \sqrt{mx(0)}, x_2(0) = \sqrt{(1-m)x(0)}$ for different $0<m<1.$ Moreover, the driving Wiener process in this case is produced in the following way
\beqq\label{LSD-eq:wienerpr}
W(t) = \int_0^t \frac{x_1(u)dW^{(1)}_u + x_2(u)dW^{(2)}_u}{\sqrt{x(u)}}. 
\eeqq

In practice the increments of the Wiener process we use for the derivation of the paths of all the approximation methods are
\beqq\label{LSD-eq:wienerprincr}
\D W_n =  \frac{x_1(t_n)\D W^{(1)}_n + x_2(t_n)\D W^{(2)}_n}{\sqrt{x(t_n)}}. 
\eeqq

\begin{figure}[ht]
	\centering
	\begin{subfigure}{.47\textwidth}
		\centering
		\includegraphics[width=1\textwidth]{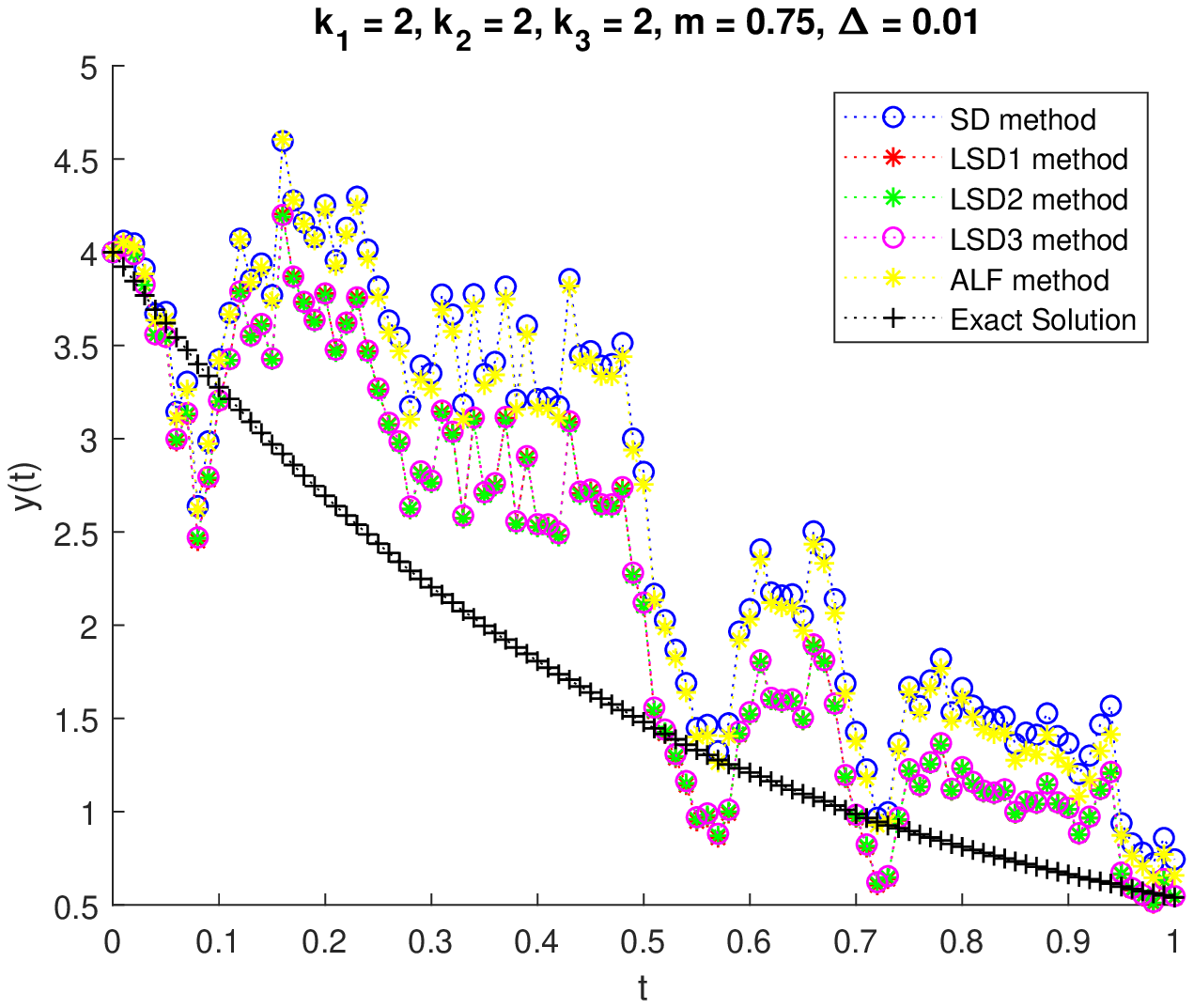}
	\end{subfigure}
	\begin{subfigure}{.47\textwidth}
		\centering
		\includegraphics[width=1\textwidth]{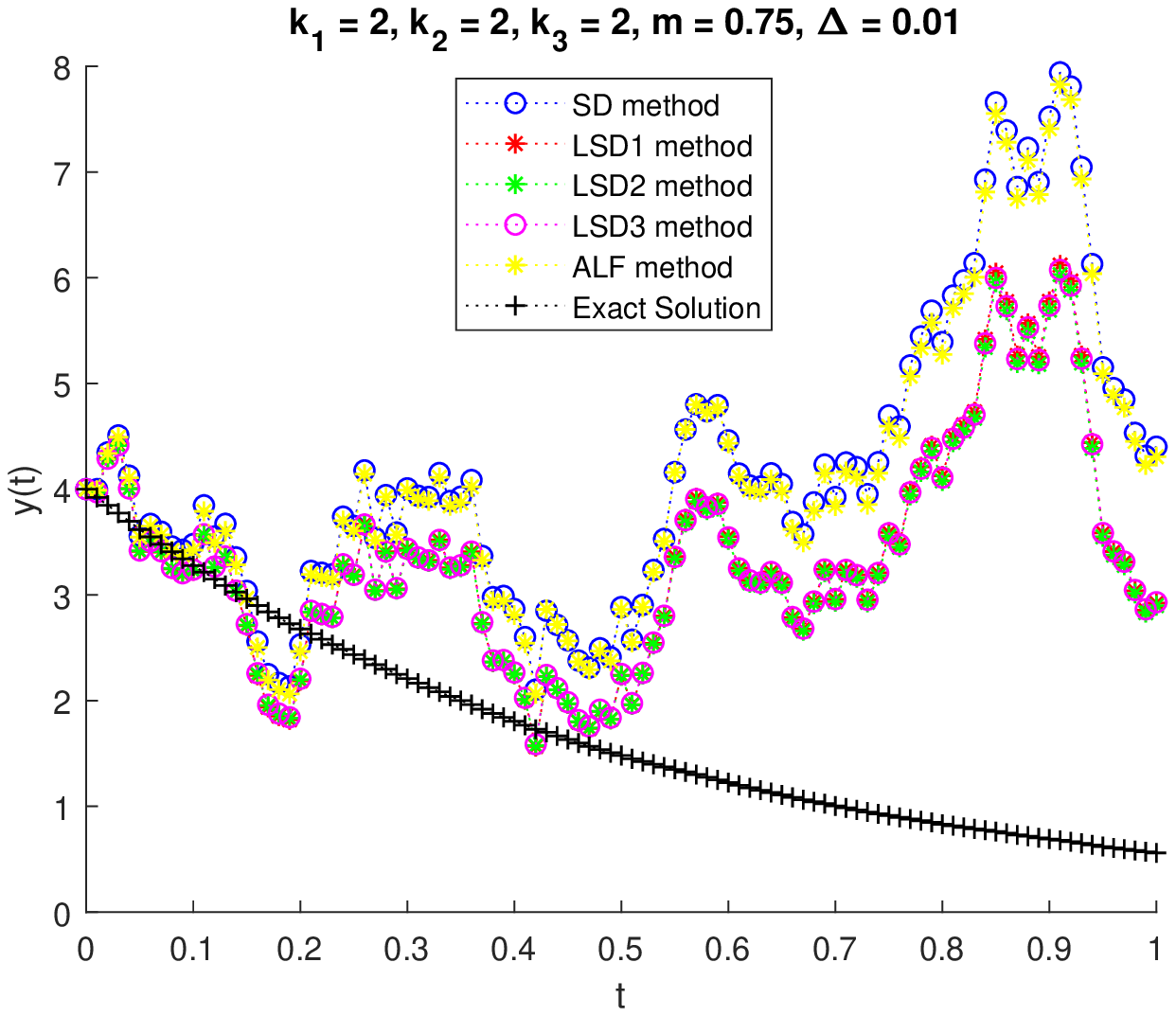}
	\end{subfigure}
	\begin{subfigure}{.47\textwidth}
	\includegraphics[width=1\textwidth]{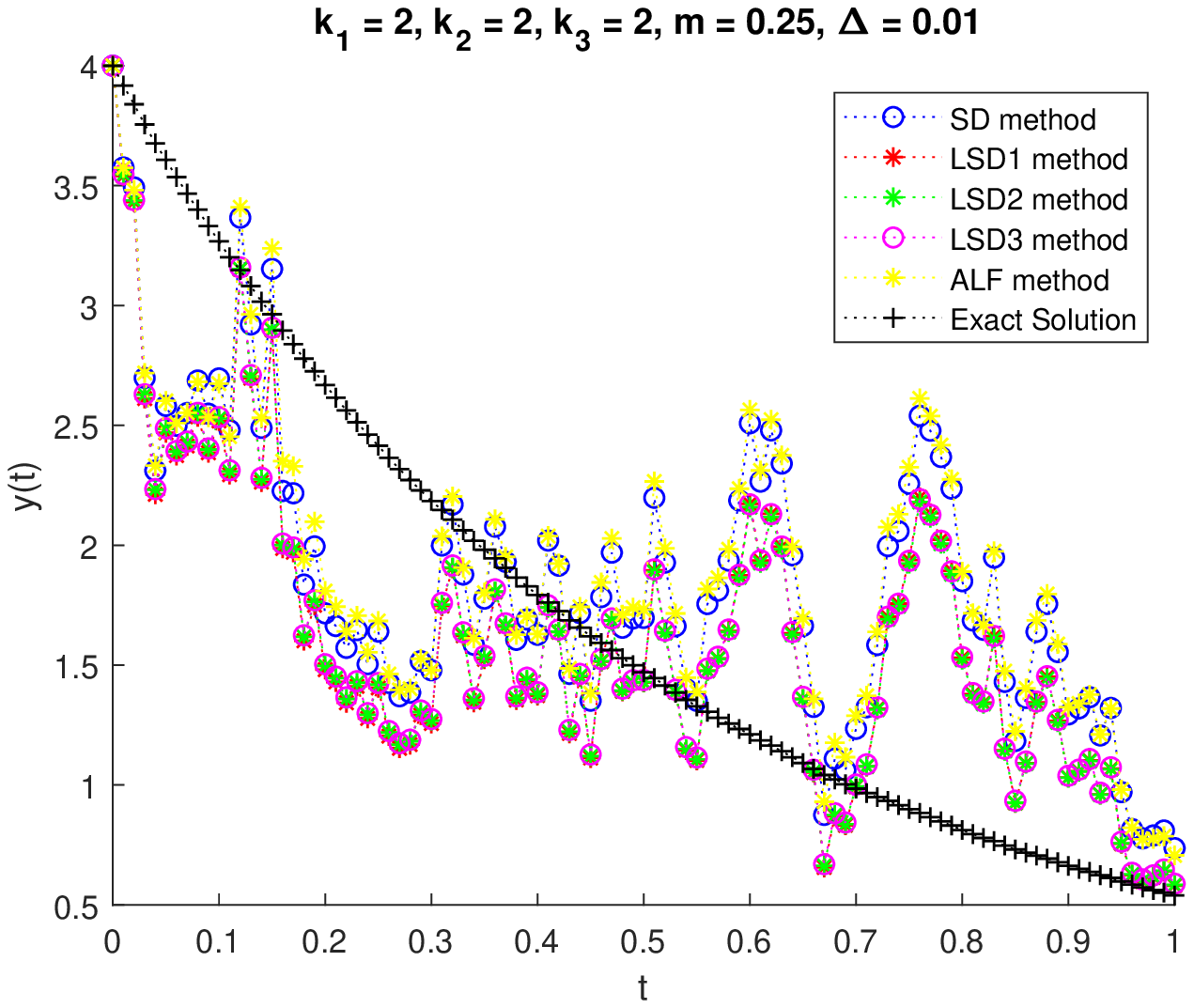}
\end{subfigure}
\begin{subfigure}{.47\textwidth}
	\centering
	\includegraphics[width=1\textwidth]{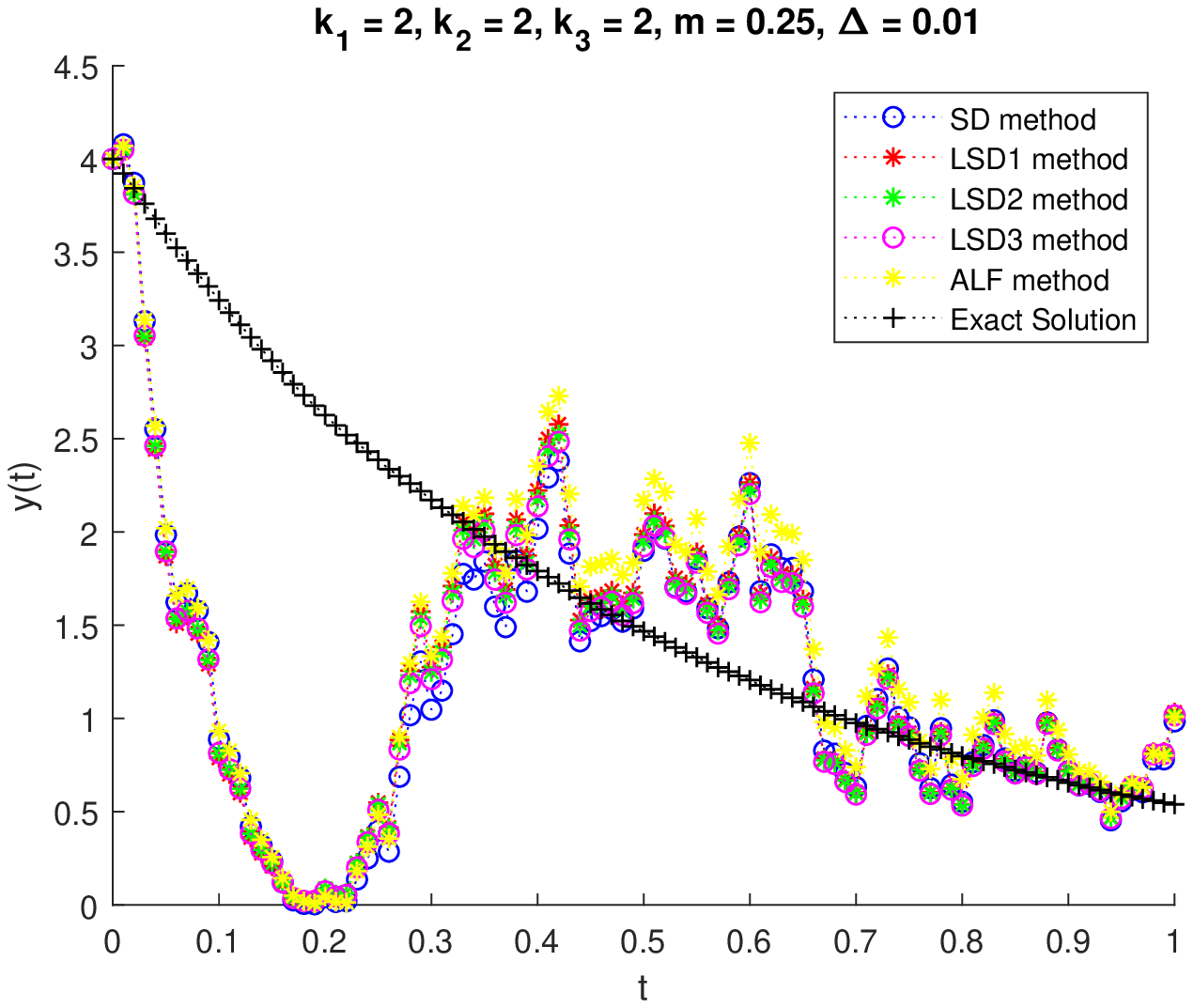}
\end{subfigure}
	\caption{Trajectories of (\ref{LSD-eq:SD schemeExampleLToriginal})-(\ref{LSD-eq:ALF_sdesol}) and the exact solution (\ref{LSD-eq:OUsolex}) with different $m$ and $\D.$}\label{LSD-fig:LSDvsExact2}
\end{figure}

In Figure \ref{LSD-fig:LSDvsExact} we cannot see the differences between the methods. By considering bigger step-sizes $\D,$ we take the picture in Figure \ref{LSD-fig:LSDvsExact2} where again we see the relation between the schemes. 

\begin{figure}[ht]
	\centering
	\begin{subfigure}{.47\textwidth}
		\centering
		\includegraphics[width=1\textwidth]{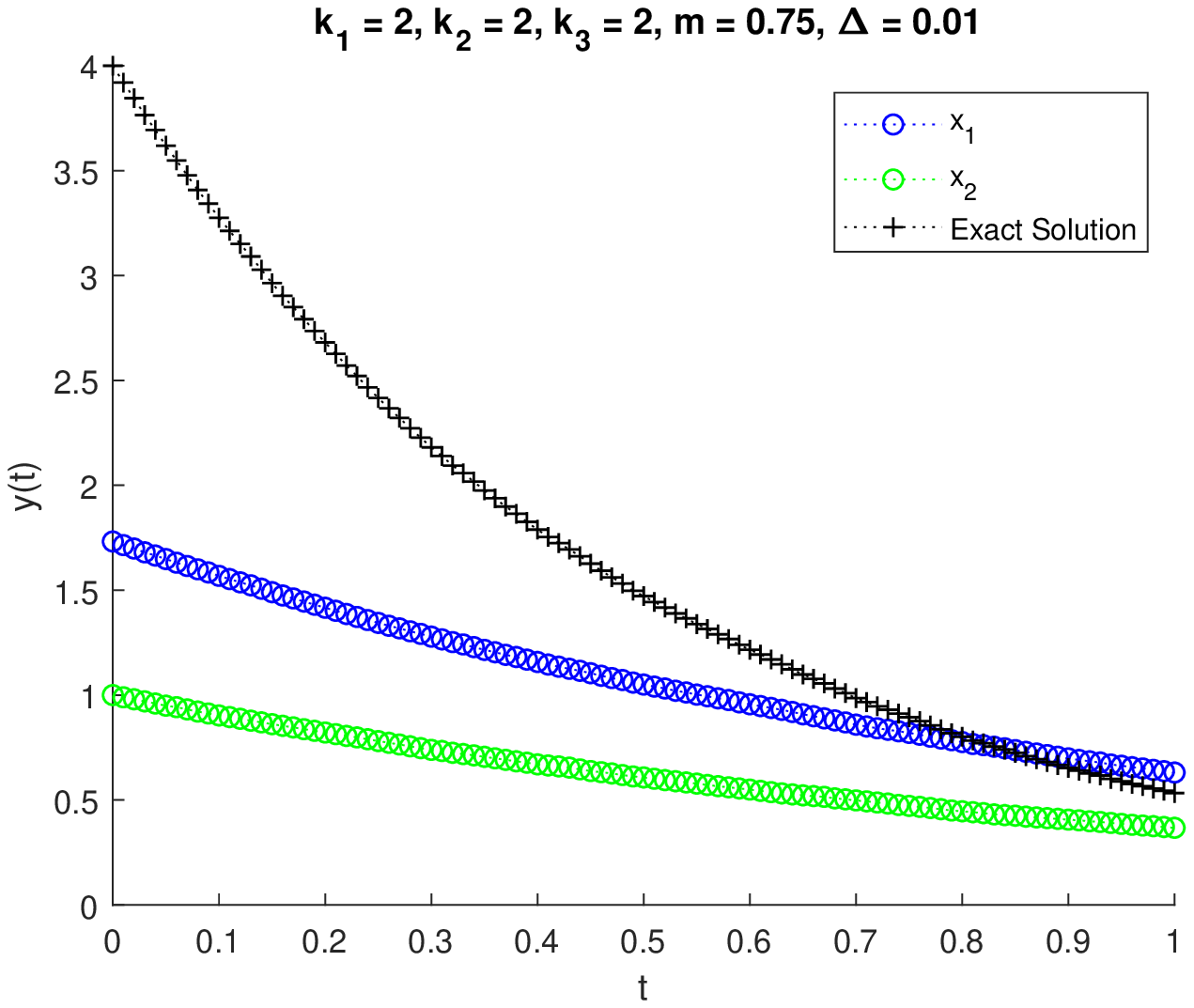}
	\end{subfigure}
	\begin{subfigure}{.47\textwidth}
		\centering
		\includegraphics[width=1\textwidth]{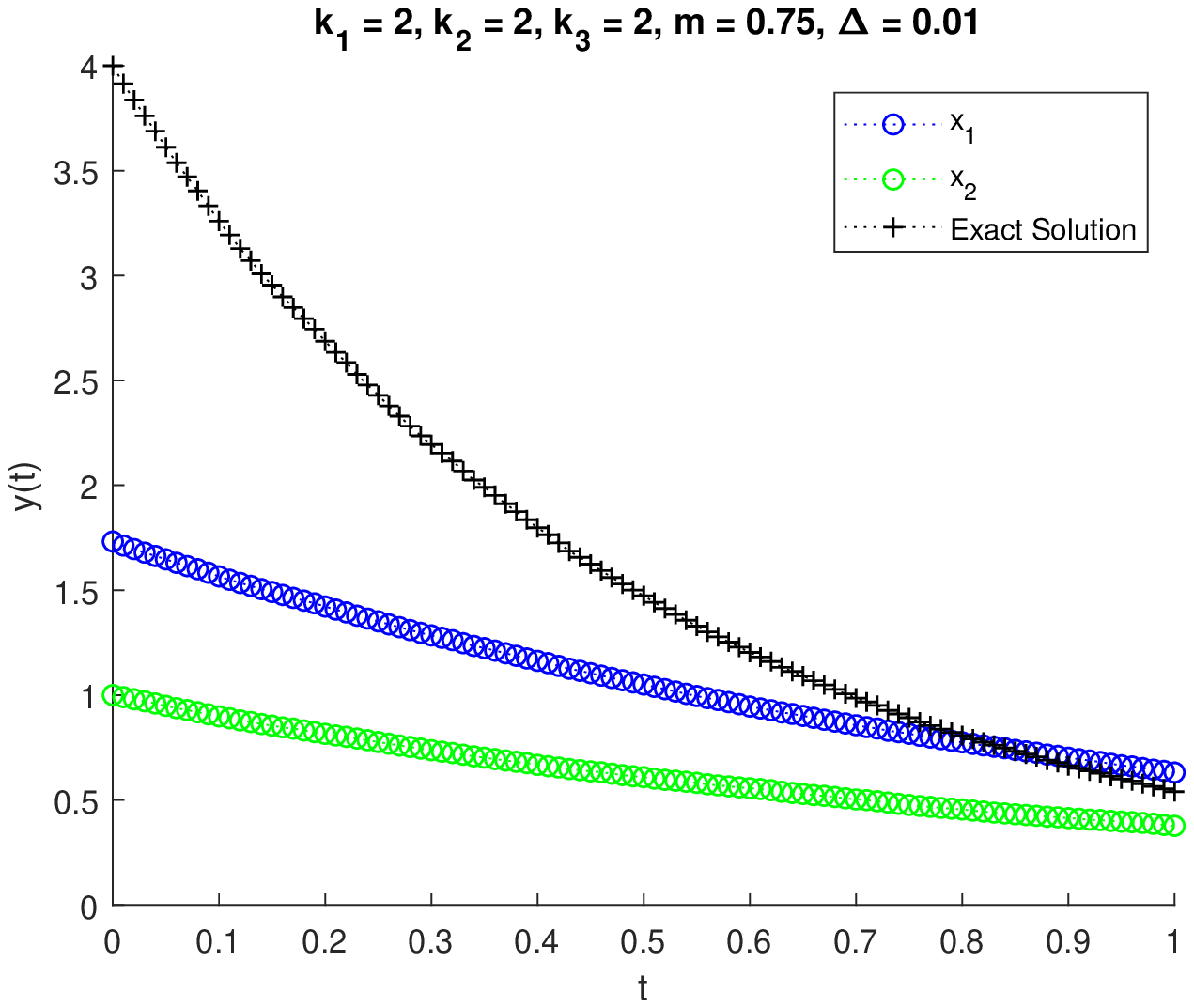}
	\end{subfigure}
	\begin{subfigure}{.47\textwidth}
		\includegraphics[width=1\textwidth]{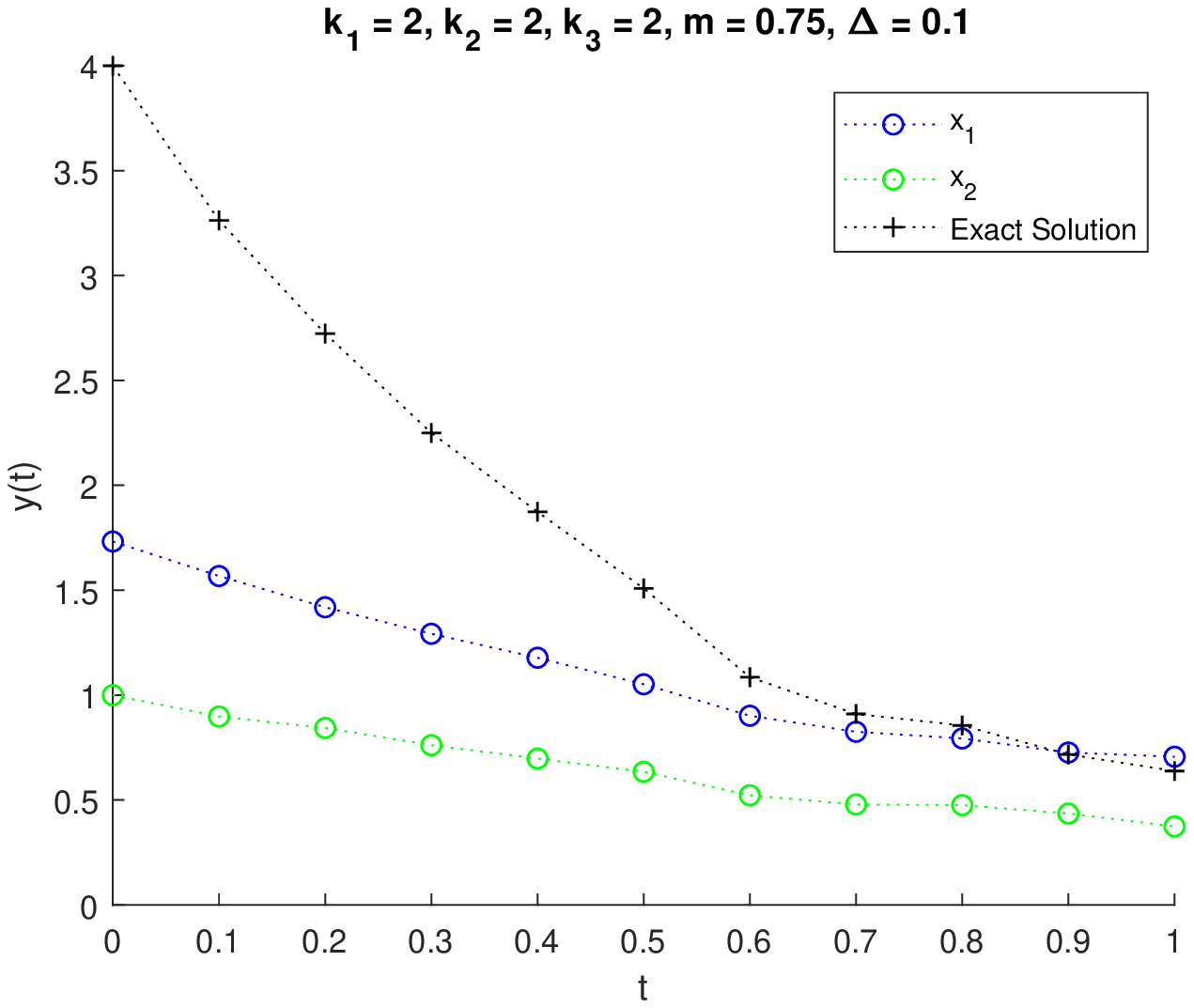}
	\end{subfigure}
	\begin{subfigure}{.47\textwidth}
		\centering
		\includegraphics[width=1\textwidth]{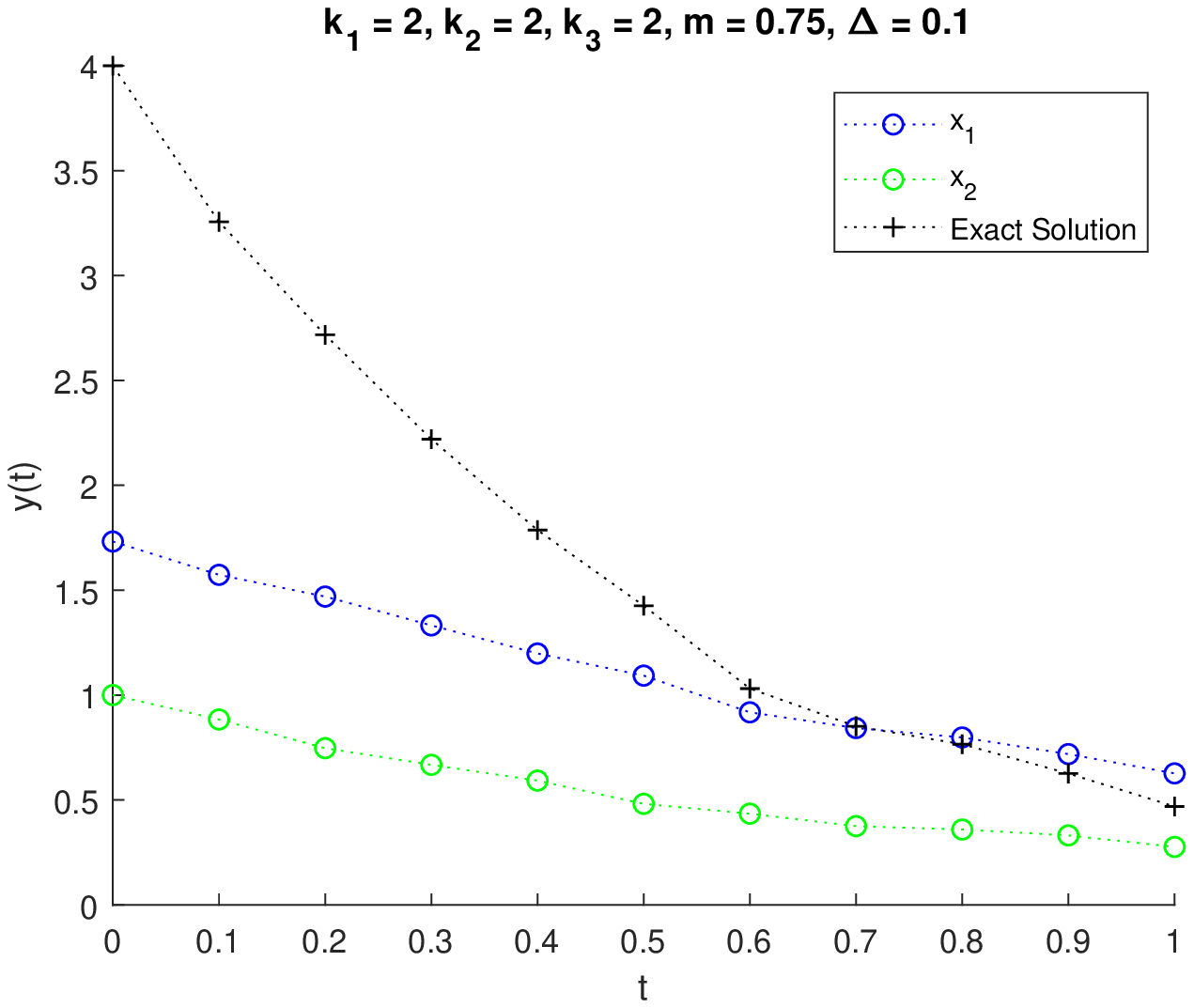}
	\end{subfigure}
	\caption{Trajectories of the processes $x_1$ and $x_2$ which give the exact solution (\ref{LSD-eq:OUsolex}) with different step-sizes $\D.$}\label{LSD-fig:Exact_constr}
\end{figure}

We note that the exact solution has a very similar behavior between different realizations of the Wiener processes. 
We need to take a bigger $\D = 0.1$ to notice a small variation of the produced solution, see Figure \ref{LSD-fig:Exact_constr}.

\begin{figure}[ht]
	\centering
	\begin{subfigure}{.47\textwidth}
		\centering
		\includegraphics[width=1\textwidth]{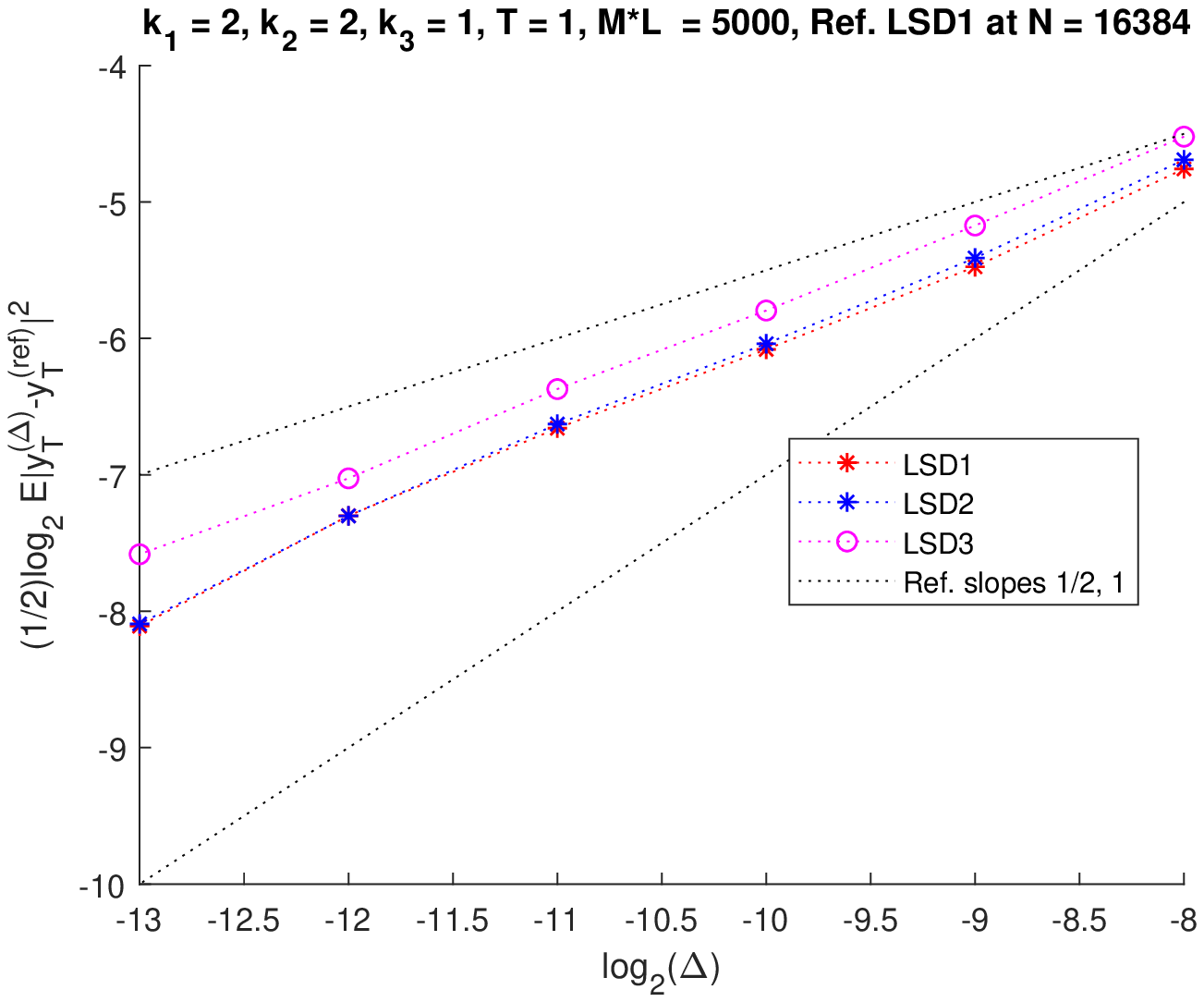}
		\caption{LSD1 as reference solution.}
	\end{subfigure}
	\begin{subfigure}{.47\textwidth}
		\centering
		\includegraphics[width=1\textwidth]{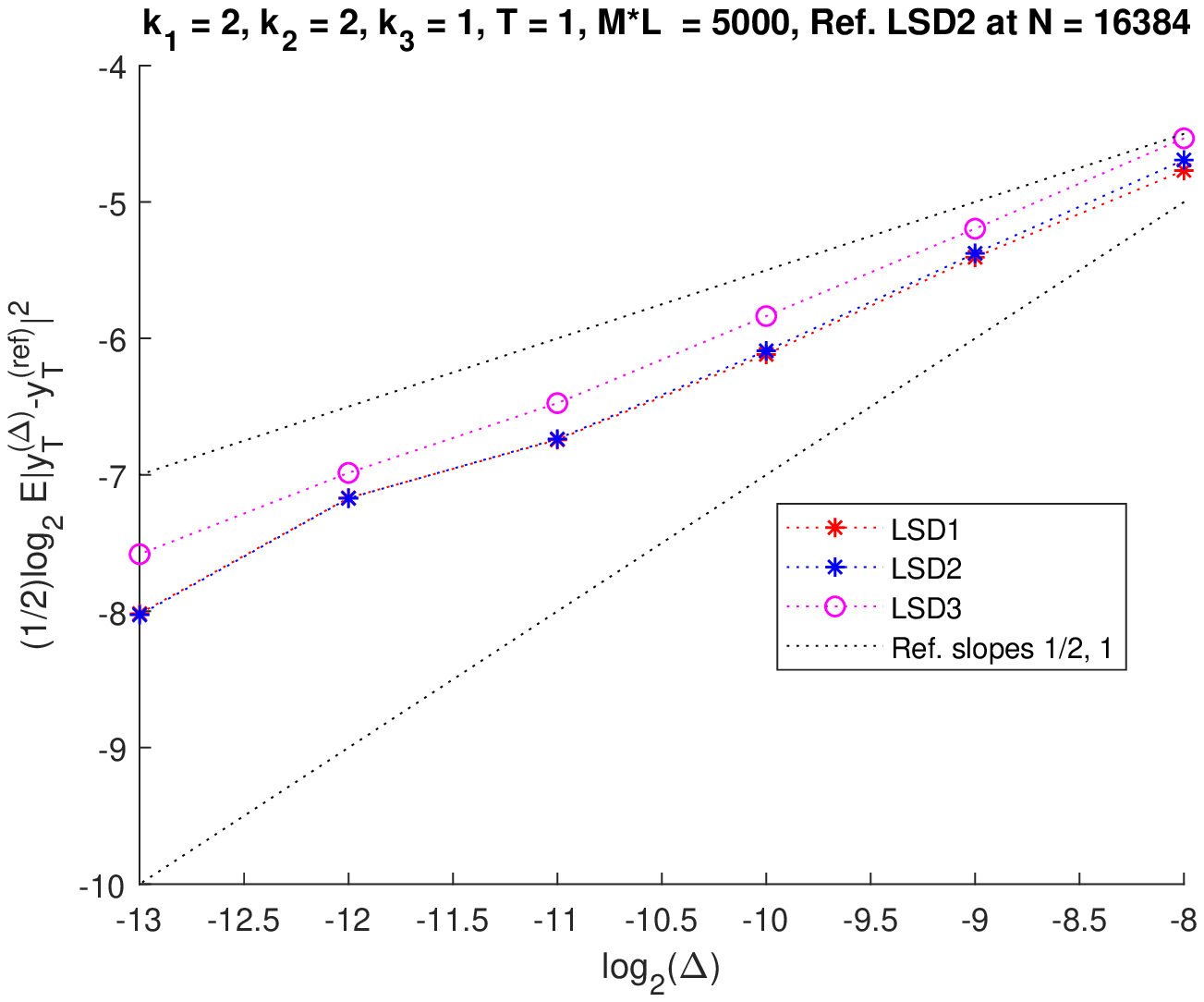}
		\caption{LSD2 as reference solution.}
	\end{subfigure}
	\begin{subfigure}{.44\textwidth}
	\centering
	\includegraphics[width=1\textwidth]{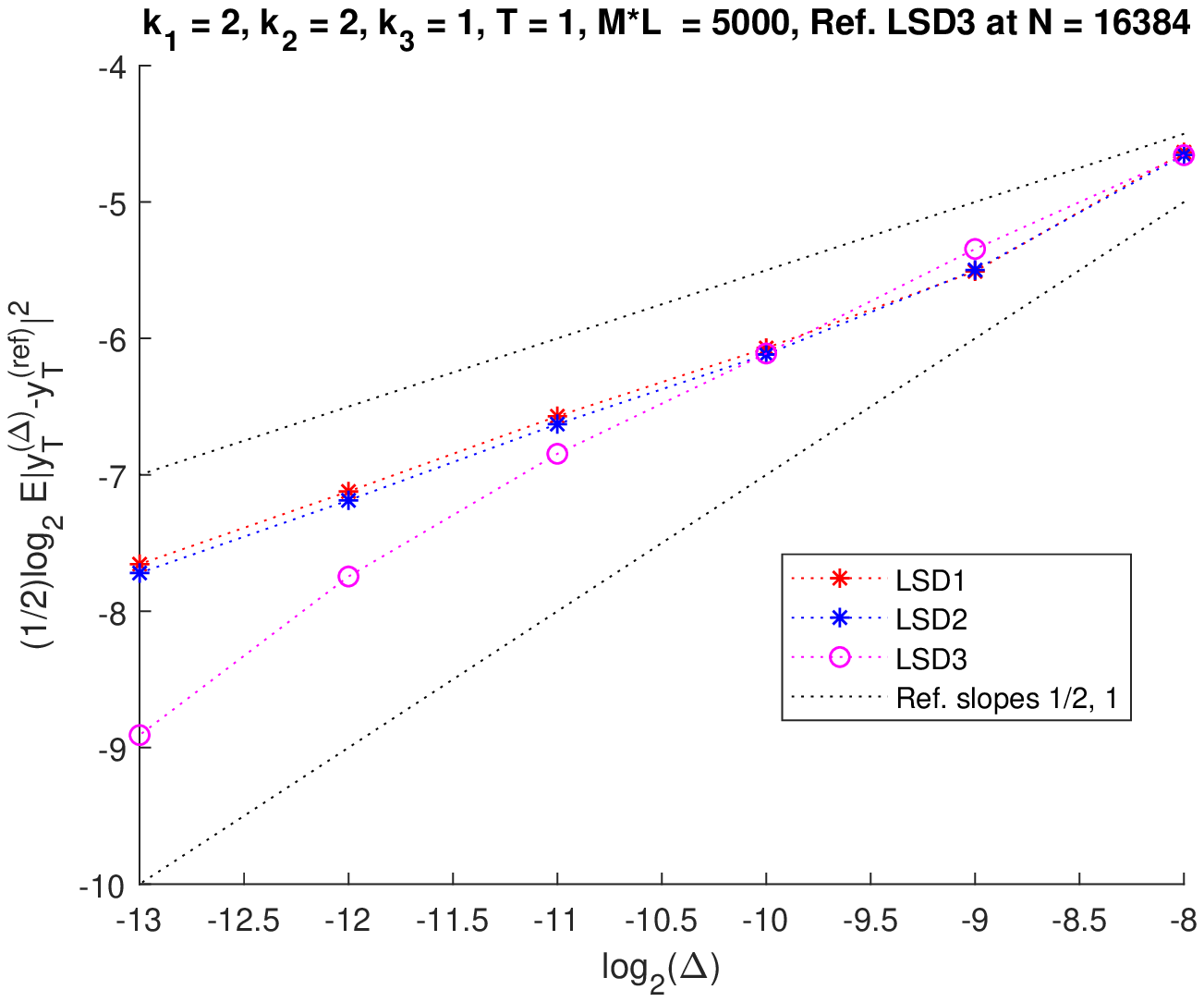}
	\caption{LSD3 as reference solution.}
\end{subfigure}
	\caption{Convergence of  Lamperti Semi-Discrete methods (\ref{LSD-eq:SD schemeExampleLToriginal}), (\ref{LSD-eq:SD schemeExampleLToriginal2}) and (\ref{LSD-eq:SD schemeExampleLToriginal3}) for the approximation of (\ref{LSD-eq:exampleSDE}) with different reference solution.}\label{LSD-fig:LSDorder}
\end{figure}

We also examine numerically the order of strong convergence of the LSD method. The numerical results suggest that the Lamperti Semi-Discrete methods converge in the mean-square sense with order close to $1,$ see Figure \ref{LSD-fig:LSDorder}. 

\begin{figure}[ht]
	\centering
	\begin{subfigure}{.47\textwidth}
		\includegraphics[width=1\textwidth]{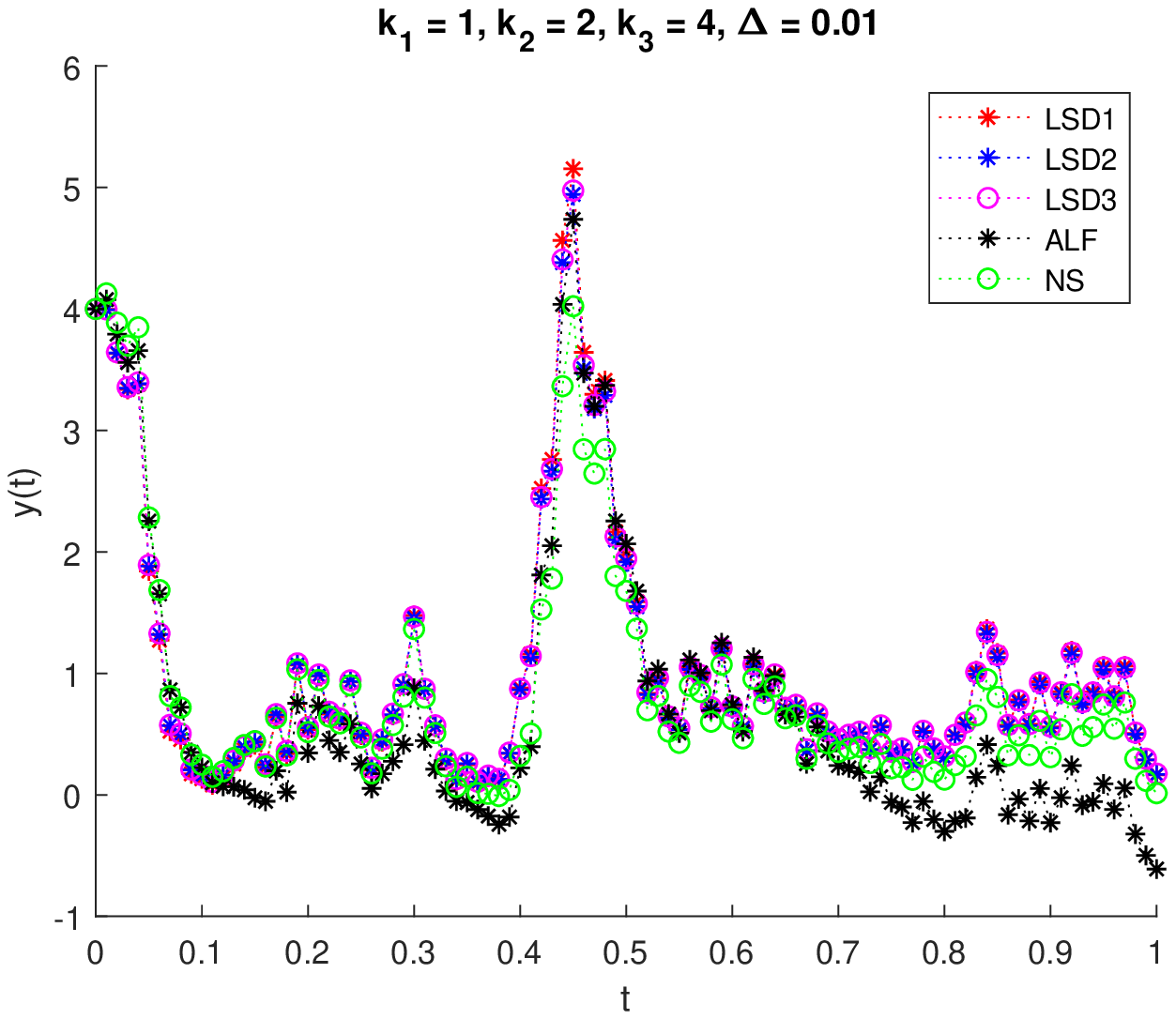}
		\caption{$k_3 = 4, \D = 10^{-2}.$}
	\end{subfigure}
	\begin{subfigure}{.47\textwidth}
		\includegraphics[width=1\textwidth]{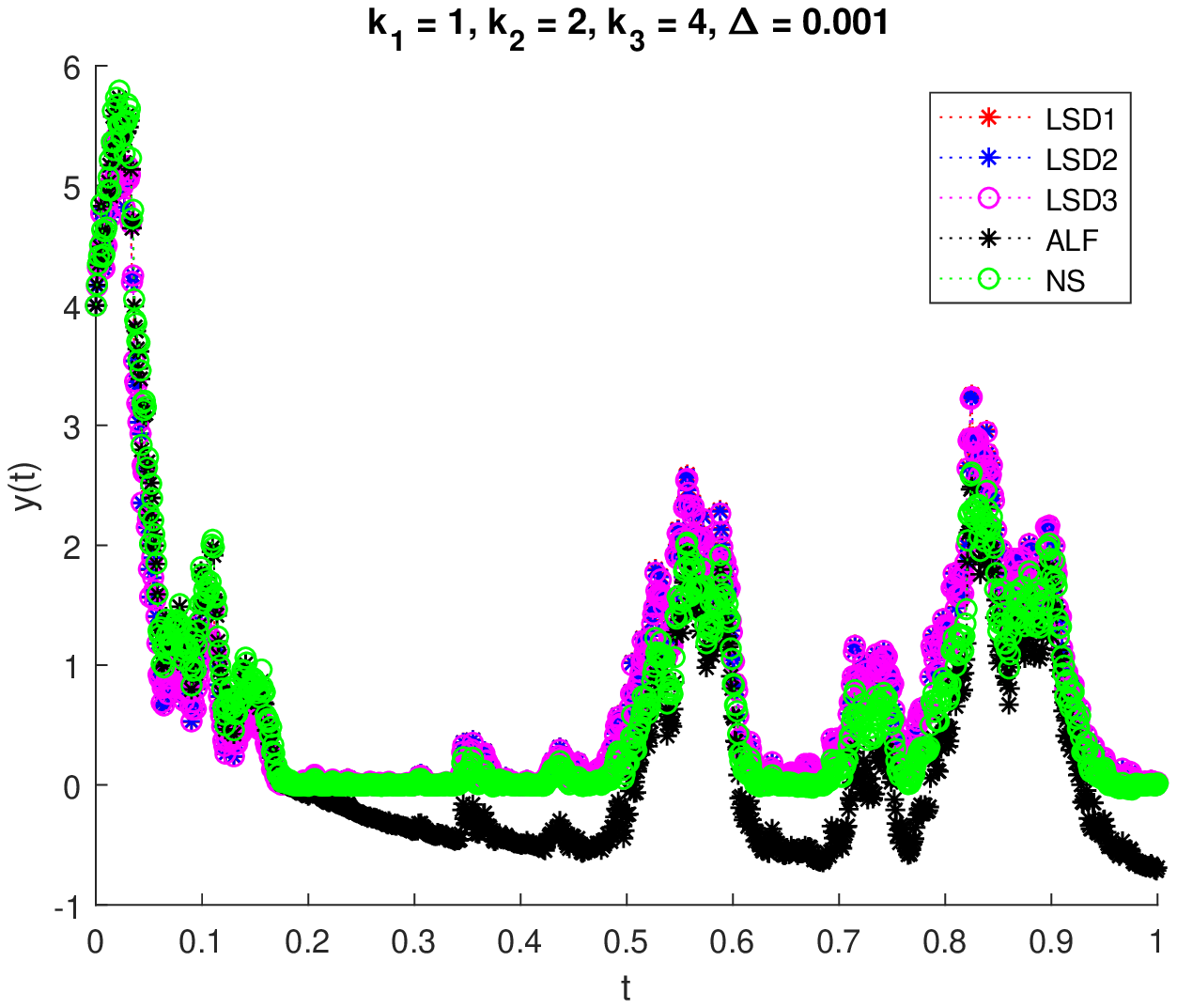}
		\caption{$k_3 = 4, \D = 10^{-3}.$}
	\end{subfigure}

	\begin{subfigure}{.47\textwidth}
	\includegraphics[width=1\textwidth]{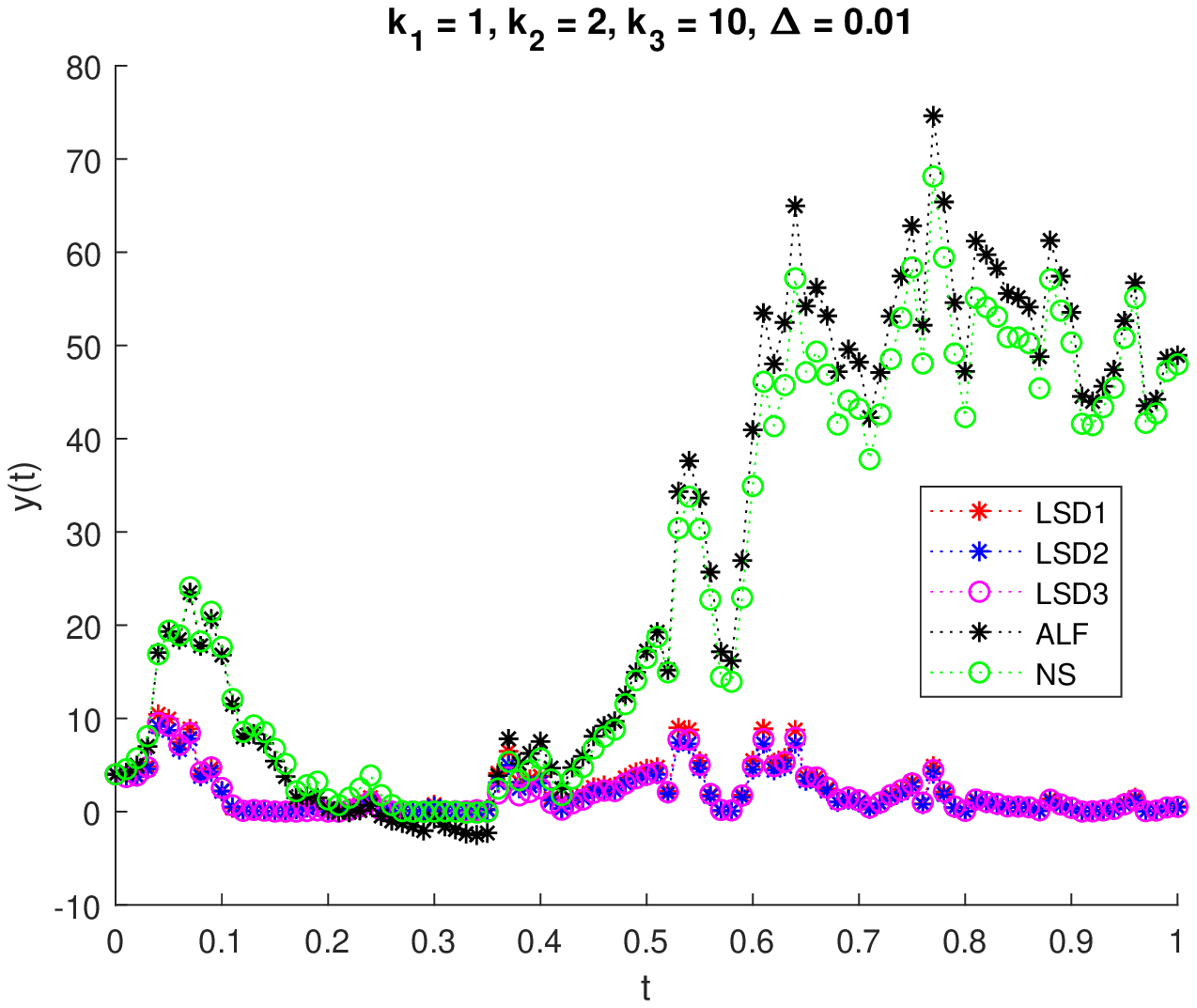}
	\caption{$k_3 = 10, \D = 10^{-2}.$}
\end{subfigure}
\begin{subfigure}{.47\textwidth}
	\includegraphics[width=1\textwidth]{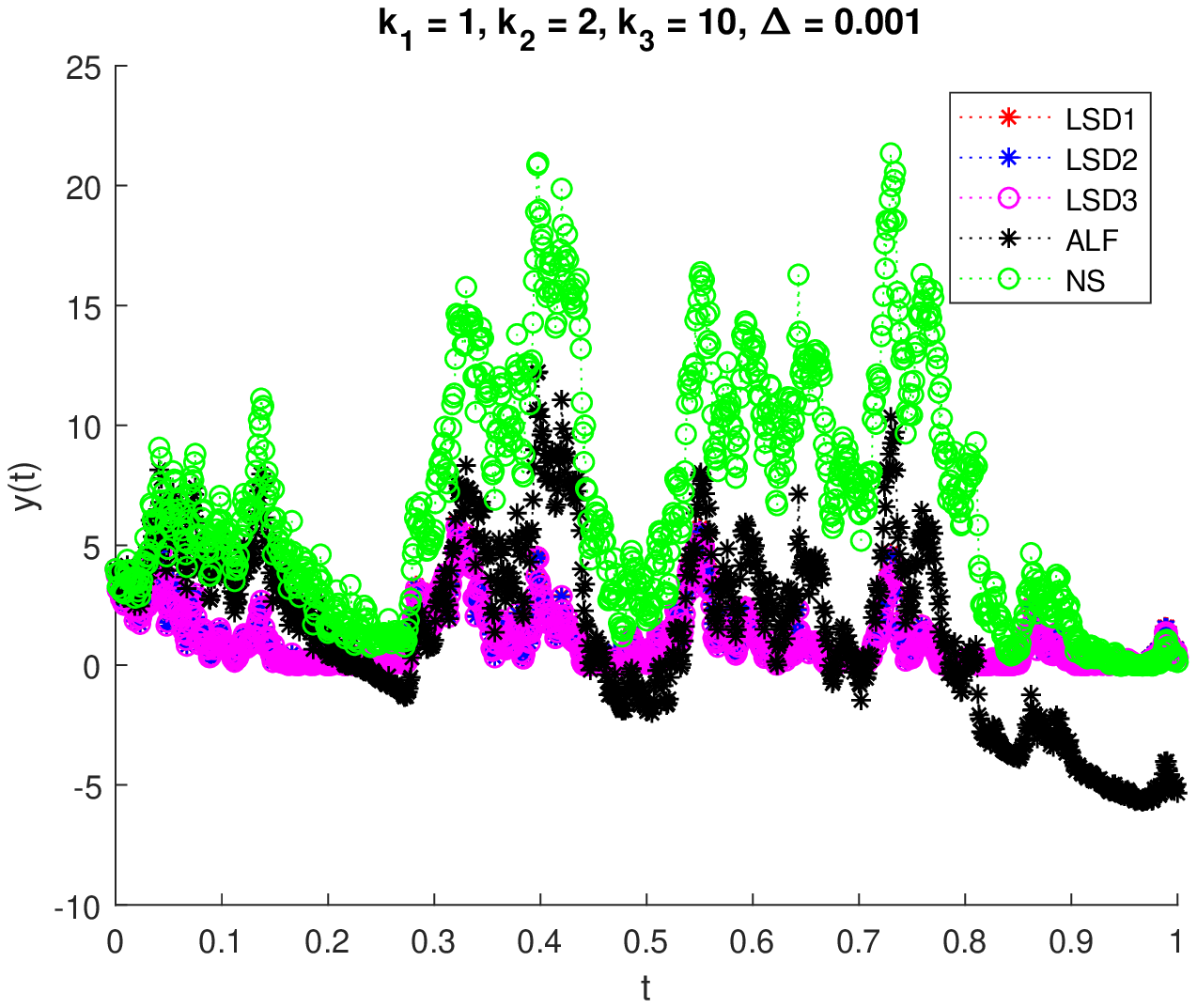}
	\caption{$k_3 = 10, \D = 10^{-3}.$}
\end{subfigure}

	\begin{subfigure}{.47\textwidth}
	\includegraphics[width=1\textwidth]{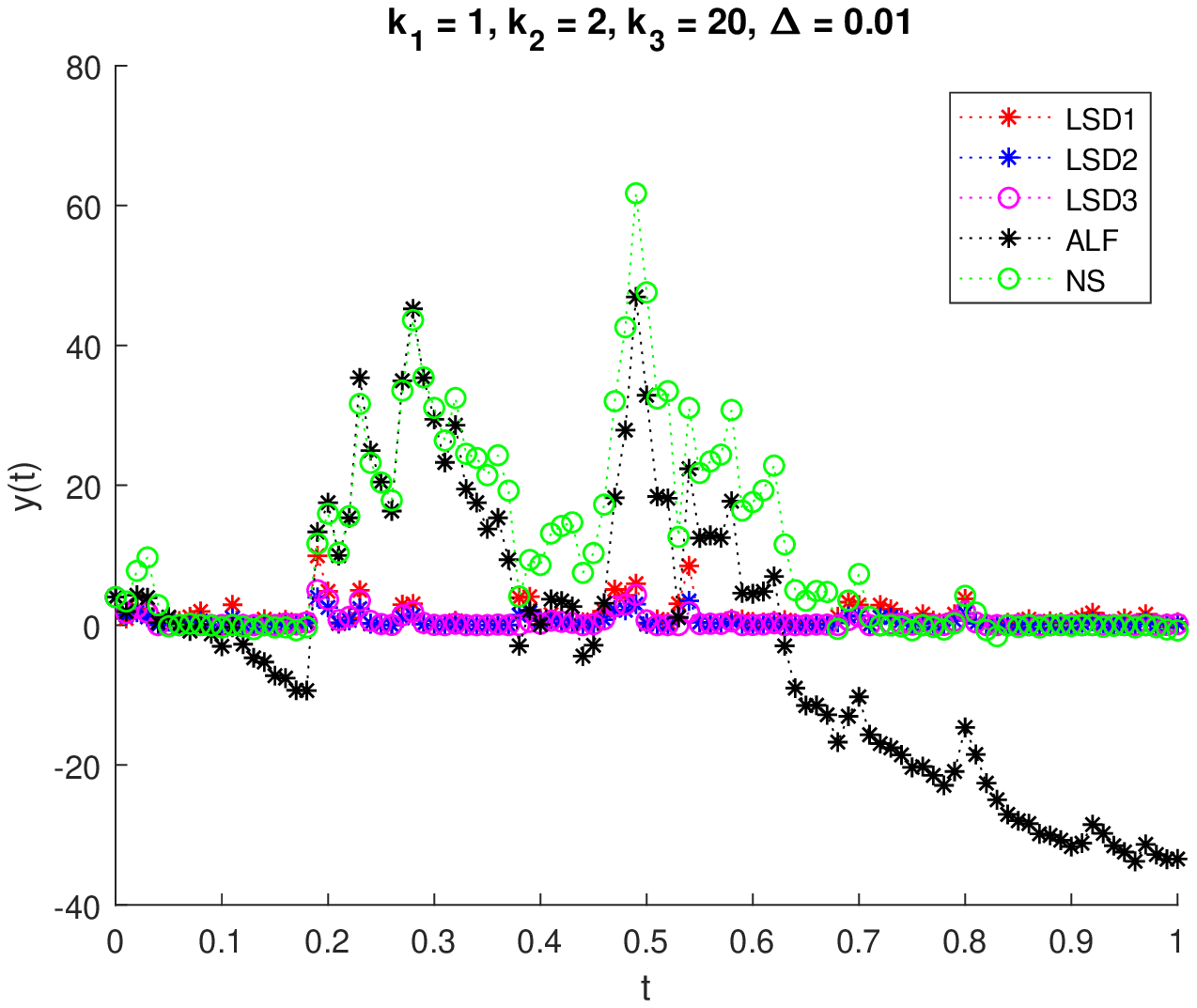}
	\caption{$k_3 = 20, \D = 10^{-2}.$}
\end{subfigure}
\begin{subfigure}{.47\textwidth}
	\includegraphics[width=1\textwidth]{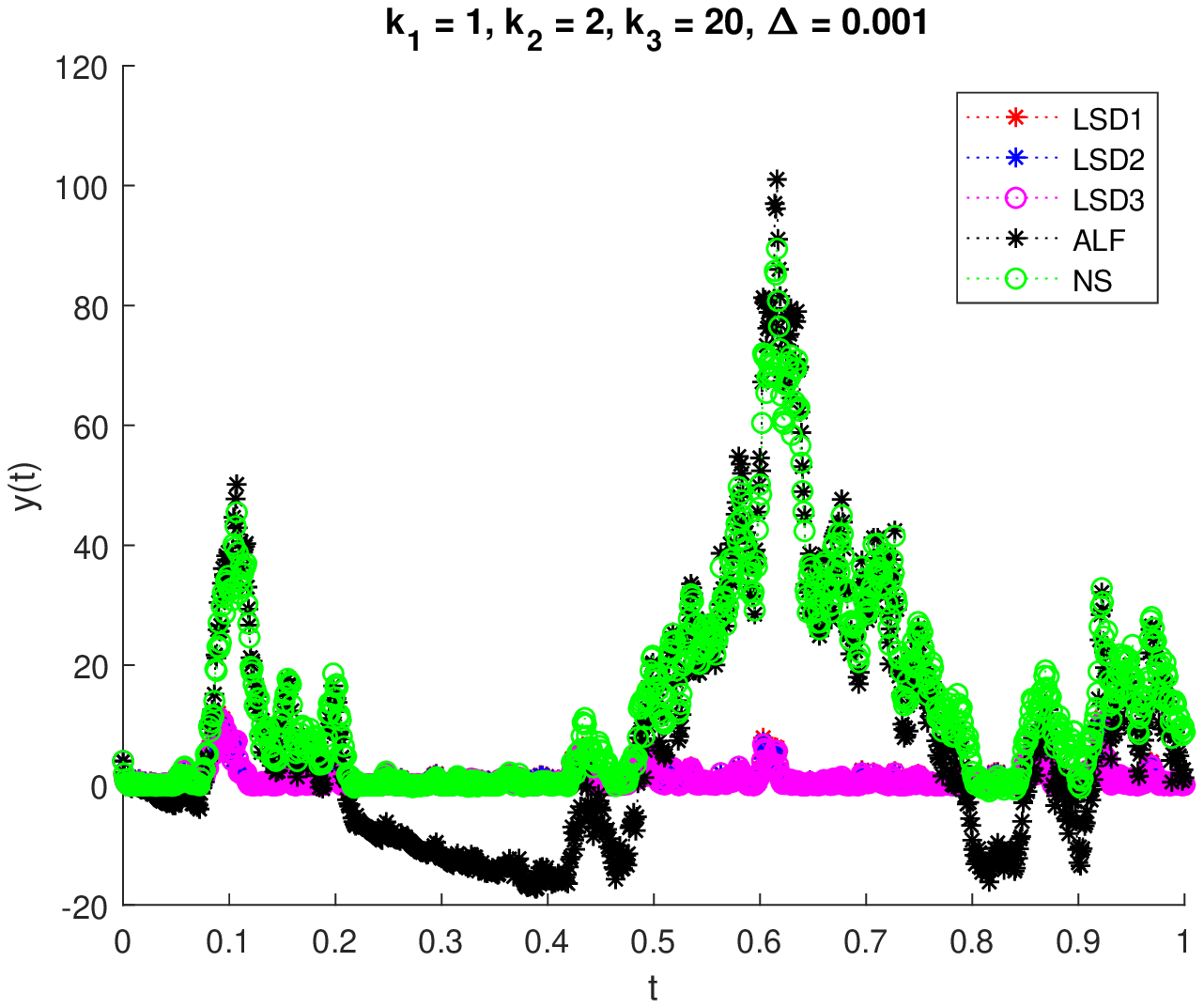}
	\caption{$k_3 = 20, \D = 10^{-3}.$}
\end{subfigure}
	\caption{Trajectories  of  (\ref{LSD-eq:SD schemeExampleLToriginal}), (\ref{LSD-eq:SD schemeExampleLToriginal2}), (\ref{LSD-eq:SD schemeExampleLToriginal3}), (\ref{LSD-eq:ALF_sdesol}) and (\ref{LSD-eq:NScir}) for the approximation of (\ref{LSD-eq:exampleSDE}) with different coefficients and various $\D.$}\label{LSD-fig:LSDCIR}
\end{figure}

Moreover, we perform one more numerical experiment to show the ability of the method to produce nonnegative solutions, outside the usual restrictions on the parameters, $(k_3)^2\leq2k_1.$ The solution of the CIR process is nonnegative, with no extra restriction on the positive parameters $k_i, i = 1, 2, 3,$ i.e. $x_t\geq0$ a.s. when $x_0>0.$ Therefore, we change a bit the parameters, taking $k_1 = 1, k_2= 2$ and different values for $k_3$ so that $(k_3)^2>2k_1.$ In this case the proposed LSD methods (\ref{LSD-eq:SD schemeExampleLToriginal}), (\ref{LSD-eq:SD schemeExampleLToriginal2}) and (\ref{LSD-eq:SD schemeExampleLToriginal3}) work in the sense that they produce nonnegative values, whereas the  implicit method (\ref{LSD-eq:ALF_sdesol}) as well as the implicit method proposed in \cite{neuenkirch_szpruch:2014}, which also takes an explicit representation in this case,

\beqq\label{LSD-eq:NScir} 
\check{y}_{t_{n+1}}=\left( \frac{\sqrt{(\check{y}_{t_n}  + \frac{1}{2}k_3\D W_n)^2 + (k_1 - \frac{(k_3)^2}{4})\D} + \check{y}_{t_n} + \frac{1}{2}k_3\D W_n}{2 + k_2\D}\right)^{2},
\eeqq
with $\check{y}_{0} = \sqrt{x_0},$ do not even produce real values, see Figure \ref{LSD-fig:LSDCIR} where for (\ref{LSD-eq:ALF_sdesol}) and (\ref{LSD-eq:NScir}) the real parts of the solution is presented. Note that as $k_3$ increases the implicit methods (\ref{LSD-eq:ALF_sdesol}) and (\ref{LSD-eq:NScir}) show an erratic behavior.

\begin{figure}[ht]
	\centering
	\begin{subfigure}{.47\textwidth}
		\centering
		\includegraphics[width=1\textwidth]{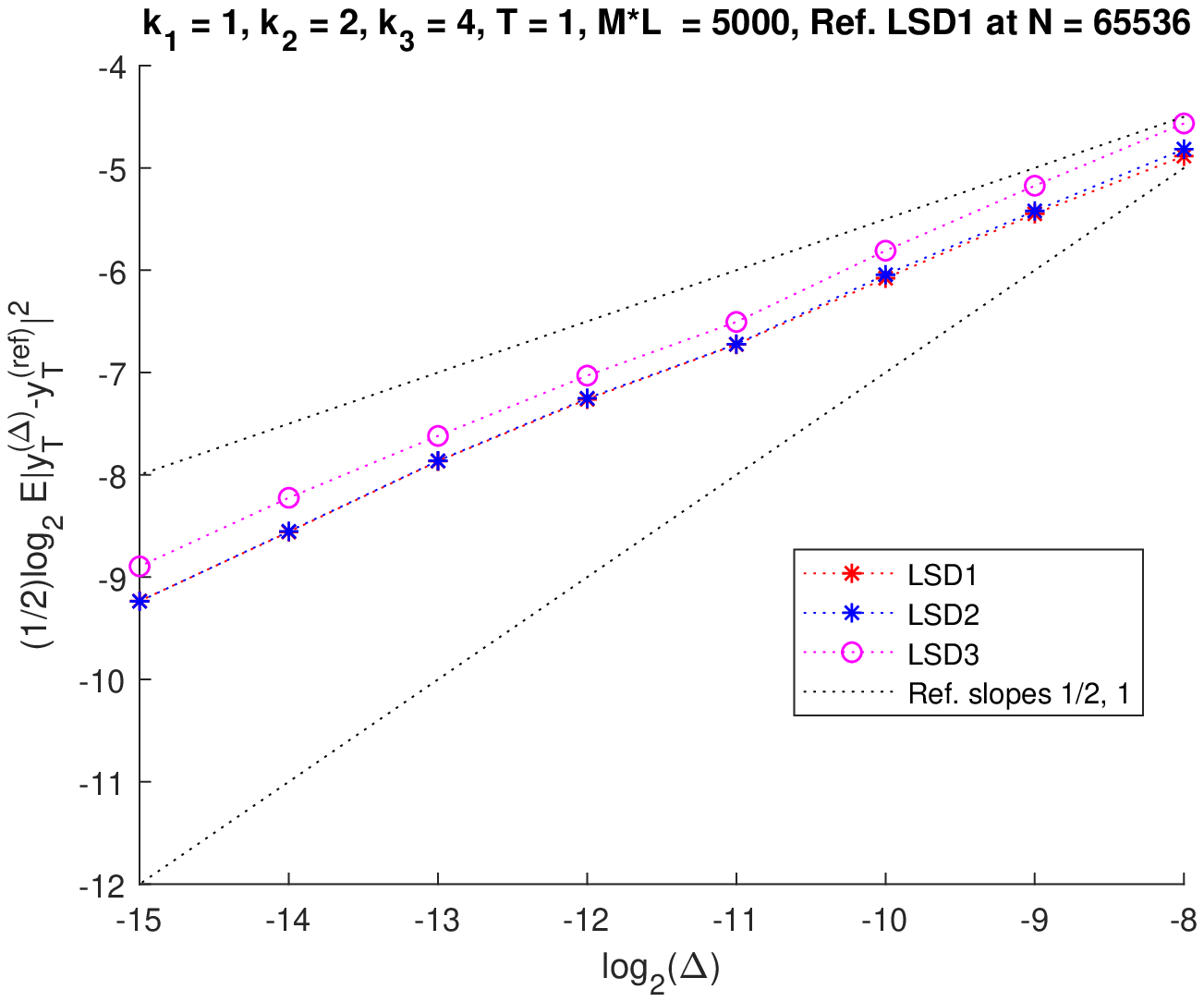}
		\caption{LSD1 as reference solution.}
	\end{subfigure}
	\begin{subfigure}{.47\textwidth}
		\centering
		\includegraphics[width=1\textwidth]{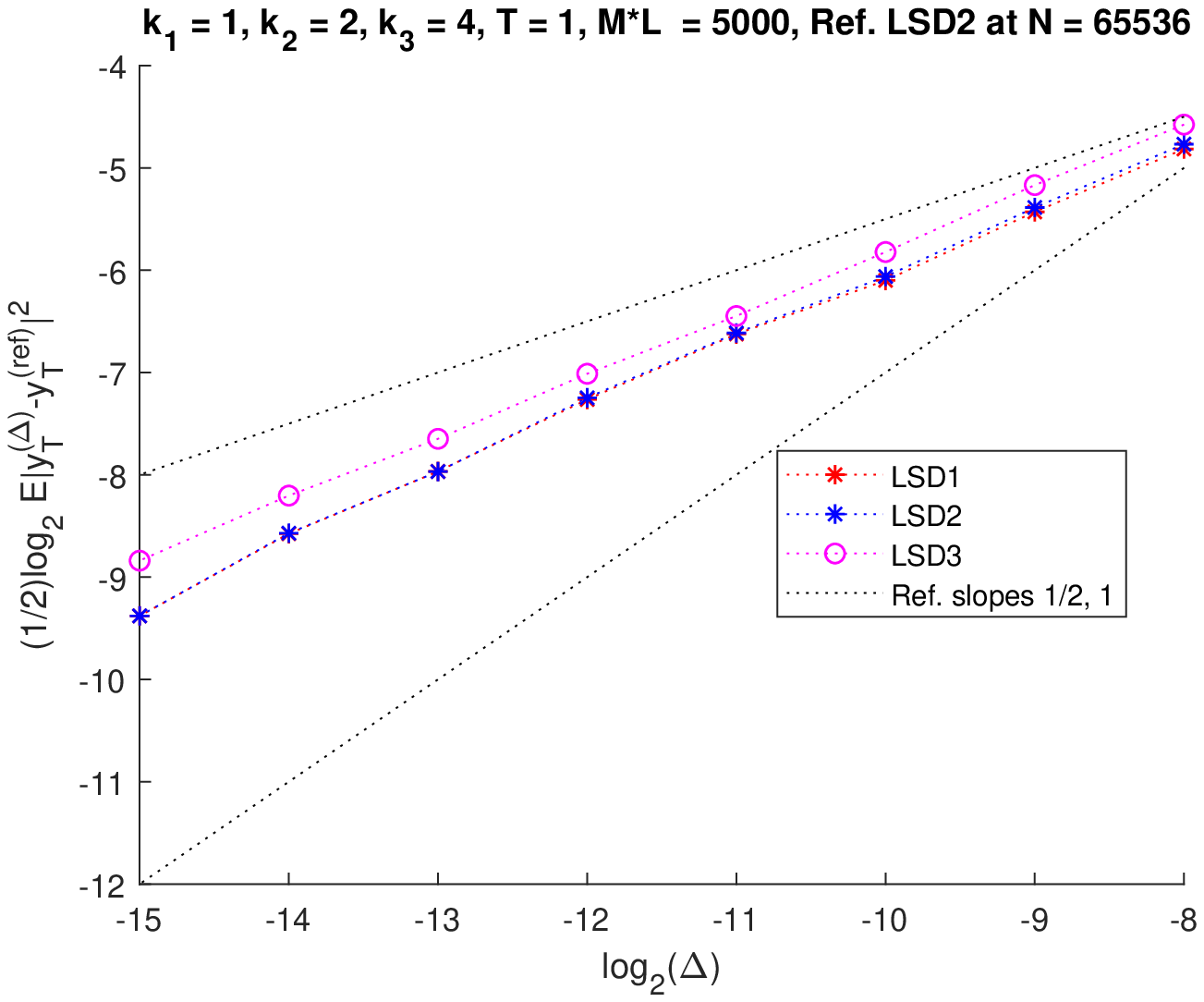}
		\caption{LSD2 as reference solution.}
	\end{subfigure}
	\begin{subfigure}{.44\textwidth}
		\centering
		\includegraphics[width=1\textwidth]{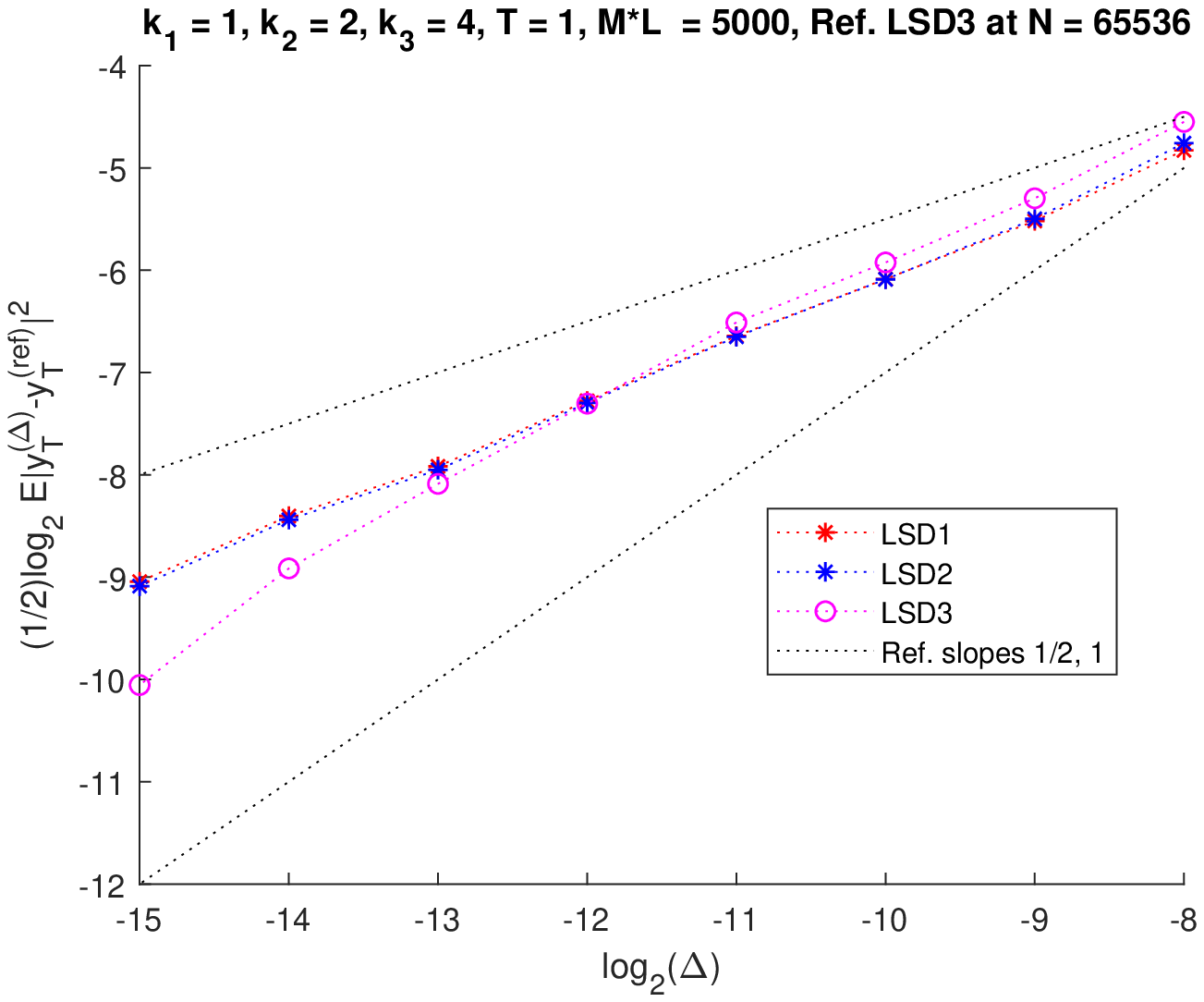}
		\caption{LSD3 as reference solution.}
	\end{subfigure}
	\caption{Convergence of  Lamperti Semi-Discrete methods (\ref{LSD-eq:SD schemeExampleLToriginal}), (\ref{LSD-eq:SD schemeExampleLToriginal2}) and (\ref{LSD-eq:SD schemeExampleLToriginal3}) for the approximation of (\ref{LSD-eq:exampleSDE}) with different reference solution.}\label{LSD-fig:LSDorder2}
\end{figure}

Finally, we present numerically the order of strong convergence of the LSD methods, see Figure \ref{LSD-fig:LSDorder2}, where we can once more see that it is close to $1.$

\section{CEV model}\label{LSD:ssec:CEV}
Let
\beqq  \label{LSD-eq:CEV}
x_t =x_0  + \int_0^t (k_1 - k_2x_s) ds + \int_0^t k_3(x_s)^{q} dW_s, \quad t\geq0.
\eeqq
where $k_1, k_2, k_3$ are positive and $1/2<q<1.$ SDE (\ref{LSD-eq:CEV}) is a mean-reverting constant elasticity of variance process (CEV) with $x_t>0$ a.s. c.f. \cite[App. A]{halidias_stamatiou:2015}. The Lamperti transformation of (\ref{LSD-eq:CEV}) is $z = \frac{1}{k_3(1-q)}x^{1-q}$ with dynamics, 
\beqq  \label{LSD-eq:CEVLamperti}
z_t =z_0 + \int_0^t \left(k_1k_3^{\frac{1-2q}{1-q}}(1-q)^{\frac{-q}{1-q}}(z_s)^{-\frac{q}{1-q}} - \frac{q}{2(1-q)}(z_s)^{-1} - k_2(1-q)z_s\right)ds + \int_0^t  dW_s.
\eeqq 
Set $a = k_1k_3^{\frac{1-2q}{1-q}}(1-q)^{\frac{-q}{1-q}}, b = q/(2-2q)$ and $c = k_2(1-q).$
We consider two versions of the semi-discrete method for approximating  (\ref{LSD-eq:CEVLamperti}). In the first version $(y_t),$ see Section \ref{LSD:subsec:CEVz1}, we use the semi-discrete method as originally proposed, and in the second version we examine the new version $(\hat{y}_t),$ see Section \ref{LSD:subsec:CEVz2}, where in each subinterval $(t_n, t_{n+1}]$ we solve an algebraic equation.

\subsection{Lamperti Semi-Discrete method $\wt{z}^1_{n}$ for CEV}\label{LSD:subsec:CEVz1}

Rewrite (\ref{LSD-eq:CEVLamperti}) as
\beqq  \label{LSD-eq:CEVLamperti_abc}
z_t =z_0 + \int_0^t \left(a(z_s)^{-\frac{q}{1-q}} - b(z_s)^{-1} - cz_s\right)ds + \int_0^t  dW_s.
\eeqq
where $a,b,c$ positive. We consider the following semi-discrete method for approximating  (\ref{LSD-eq:CEVLamperti}),  
\beam\nonumber
y_t &=& \D W_n + y_{t_n}  +  \int_{t_n}^t \left(\frac{a}{(y_{t_n})^{\frac{2q-1}{1-q}}} (y_{s})^{-1} - \frac{b}{(y_{t_n})^2}y_t - cy_s\right)ds\\
\label{LSD-eq:SD schemeCEVLT}& = & \phi_\D(y_{t_n},\D W_n) + \int_{t_n}^t \left(B_n (y_{s})^{-1}  + C_ny_s\right)ds
\eeam
with $y_{0}=z_0,$ where $\phi_\D(x,y) = \frac{x + y}{1 + \frac{b}{x^2}\D}$ and $$B_n:=\frac{a}{(y_{t_n})^{\frac{2q-1}{1-q}} + b(y_{t_n})^{\frac{4q-3}{1-q}}\D}, \,\,C_n:= - \frac{c}{1 + \frac{b\D}{(y_{t_n})^2}}.$$ 

The Bernoulli equation (\ref{LSD-eq:SD schemeCEVLT}) has a solution satisfying, see Appendix \ref{LSD-ap:Bernoulli_sol},
\beam\nonumber
(y_t)^{2} &=& \phi^2_\D(y_{t_n},\D W_n)\exp\left\{-\frac{2c}{1 + \frac{b\D}{(y_{t_n})^2}}(t-t_n)\right\} \\
\label{LSD-eq:SD schemeCEVLTsol}
&&+ \frac{a}{c(y_{t_n})^{\frac{2q-1}{1-q}}}\left(1 - \exp\left\{-\frac{2c}{1 + \frac{b\D}{(y_{t_n})^2}}(t-t_n)\right\}\right).
\eeam

We propose the following version of the semi-discrete method for the approximation of (\ref{LSD-eq:CEVLamperti}),  
\beqq\label{LSD-eq:SDschemeCVELT_transf}
y_{t_{n+1}} = \sqrt{\phi^2_\D(y_{t_n},\D W_n)e^{-\frac{2c}{1 + \frac{b\D}{(y_{t_n})^2}}\D}  + \frac{a}{c(y_{t_n})^{\frac{2q-1}{1-q}}}(1 - e^{-\frac{2c}{1 + \frac{b\D}{(y_{t_n})^2}}\D})}, 
\eeqq	
which suggests the versions of the Lamperti semi-discrete method $(\wt{z}^1_n)_{n\in\bbN},$ for the approximation of (\ref{LSD-eq:CEV}) with $\wt{z}^1_{n} = (k_3(1-q)y_{n})^{1/(1-q)}$ or  
\beam\label{LSD-eq:CEVLToriginal}
&&\wt{z}^1_{t_{n+1}} =  (k_3(1-q))^{\frac{1}{1-q}}\\
\nonumber&\times&\left| \phi^2_\D(y_{t_n},\D W_n)e^{-\frac{2c\D}{1 + \frac{b\D}{(y_{t_n})^2}}}  + \frac{a}{c(y_{t_n})^{\frac{2q-1}{1-q}}}(1 - e^{-\frac{2c\D}{1 + \frac{b\D}{(y_{t_n})^2}}})\right|^{1/(2-2q)}.
\eeam

\subsection{Lamperti Semi-Discrete methods $\wt{z}^2_{n}, \wt{z}^3_{n}$ for CEV}\label{LSD:subsec:CEVz2}

Let $t\in(t_n, t_{n+1}]$ and consider the processes $(\hat{y}_t)$ and $(\wt{y}_t)$ where
\beqq\label{LSD-eq:SD schemeCEVLT2}
\hat{y}_t =  W_t - W_{t_n} + \hat{y}_{t_n}  + a(\hat{y}_{t_n})^{\frac{1-2q}{1-q}}(\hat{y}_{t})^{-1}\D - b(\hat{y}_{t_n})^{-1}\D - c\hat{y}_t\D,
\eeqq
with $\hat{y}_{0}=z_0$ and
\beqq\label{LSD-eq:SD schemeCEVLT3}
\wt{y}_t =  W_t - W_{t_n} + \wt{y}_{t_n}  + a(\wt{y}_{t_n})^{\frac{-q}{1-q}}\D - b(\wt{y}_{t_n})^{-1}\D  - (\wt{y}_{t_n})^{-1}\D + (\wt{y}_{t})^{-1}\D- c\wt{y}_t\D,
\eeqq
with $\wt{y}_{0}=z_0.$ 

The solution of (\ref{LSD-eq:SD schemeCEVLT2}) is such that 
\beqq\label{LSD-eq:SD schemeCEVLTsol2}
(1+c\D)(\hat{y}_t)^{2} - \left(W_t - W_{t_n} + \hat{y}_{t_n}  - b(\hat{y}_{t_n})^{-1}\D\right)\hat{y}_t - a(\hat{y}_{t_n})^{\frac{1-2q}{1-q}}\D = 0,
\eeqq
while the solution of (\ref{LSD-eq:SD schemeCEVLT3}) satisfies 
\beqq\label{LSD-eq:SD schemeCEVLTsol3}
(1+c\D)(\wt{y}_t)^{2} - \left(W_t - W_{t_n} + \wt{y}_{t_n} + (a(\wt{y}_{t_n})^{\frac{-q}{1-q}} - b(\wt{y}_{t_n})^{-1})\D\right)\wt{y}_t - \D = 0
\eeqq

We propose the following versions of the semi-discrete method for the approximation of (\ref{LSD-eq:CEVLamperti}),  
\beqq\label{LSD-eq:SDschemeCVELT_transf2}
\hat{y}_{t_{n+1}} = \frac{\hat{\phi}_\D(\hat{y}_{t_n},\D W_n) +  
	\sqrt{\hat{\phi}^2_\D(\hat{y}_{t_n},\D W_n) + 4 a\D(1+c\D)(\hat{y}_{t_n})^{\frac{1-2q}{1-q}}}}{2(1+c\D)},
\eeqq
with $\hat{\phi}_\D(x,y) = (y + x  - bx^{-1}\D),$ and 
\beqq\label{LSD-eq:SDschemeCVELT_transf3}
\wt{y}_{t_{n+1}} = \frac{\wt{\phi}_\D(\wt{y}_{t_n},\D W_n) +  
	\sqrt{\wt{\phi}^2_\D(\wt{y}_{t_n},\D W_n) + 4 \D(1+c\D)}}{2(1+c\D)},
\eeqq
with $\wt{\phi}_\D(x,y) = (y + x  + ax^{\frac{-q}{1-q}}\D - bx^{-1}\D),$ which suggest the versions of the Lamperti semi-discrete method $(\wt{z}^2_n)_{n\in\bbN}, (\wt{z}^3_n)_{n\in\bbN}$ for the approximation of (\ref{LSD-eq:CEV}) with   
\beam\label{LSD-eq:CEVLToriginal2}
&&\wt{z}^2_{t_{n+1}} = (k_3(1-q))^{\frac{1}{1-q}}\\
\nonumber&\times&  \left| \frac{\hat{\phi}_\D(\hat{y}_{t_n},\D W_n) +  
	\sqrt{\hat{\phi}^2_\D(\hat{y}_{t_n},\D W_n) + 4 a\D(1+c\D)(\hat{y}_{t_n})^{\frac{1-2q}{1-q}}}}{2(1+c\D)} \right|^{1/(1-q)},
\eeam
and
\beam\label{LSD-eq:CEVLToriginal3}
&&\wt{z}^3_{t_{n+1}} = (k_3(1-q))^{\frac{1}{1-q}}\\
\nonumber&\times&  \left| \frac{\wt{\phi}_\D(\wt{y}_{t_n},\D W_n) +  
	\sqrt{\wt{\phi}^2_\D(\wt{y}_{t_n},\D W_n) + 4 \D(1+c\D)}}{2(1+c\D)} \right|^{1/(1-q)}.
\eeam
\subsection{Numerical experiment for CEV}\label{LSD:subsec:CEVnum}

For a minimal numerical experiment we present simulation paths for the numerical approximation of (\ref{LSD-eq:CEV}) with  $x_0=1/16$ and compare with the SD method proposed in \cite{halidias_stamatiou:2015}, which reads 

\beam\nonumber
\wt{y}_{t_{n+1}}&=&\Big( \sqrt{\wt{y}_{t_n}\left(1- \frac{k_2\D}{1+k_2\theta\D}\right) + \frac{k_1\D}{1+k_2\theta\D} - \frac{(k_3)^2\D}{4(1+k_2\theta\D)^2}(\wt{y}_{t_n})^{2q-1}}\\
\label{LSD-eq:NCEV_SD}&& +\frac{k_3}{2(1+k_2\theta\D)}(\wt{y}_{t_n})^{q -\frac{1}{2}}\D W_n\Big)^2,
\eeam
where $\theta$ represents the level of implicitness. In particular we choose the coefficients as in \cite[Sec.6]{halidias_stamatiou:2015}; we take $k_1 = \frac{1}{16}, k_2 = 1$ and $k_3 = 0.4$ and $q =3/4$ for the fully implicit SD scheme (\ref{LSD-eq:NCEV_SD}) with $\theta = 1$ and compare with the proposed versions of LSD scheme (\ref{LSD-eq:CEVLToriginal}) and (\ref{LSD-eq:CEVLToriginal2}). 

\bre\label{LSD-rem:Ht_CEV}
As in Remark \ref{LSD-rem:Ht_CIR}, we note that in the proof of the strong convergence properties of the SD scheme (\ref{LSD-eq:NCEV_SD}) proposed in \cite{halidias_stamatiou:2015} an auxiliary process $(h_t)$ appears, see \cite[Rel. (33)]{halidias_stamatiou:2015}. Following the same lines the results in \cite{halidias_stamatiou:2015} as well as in \cite{STAMATIOU:2019} where a more general CIR/CEV-type model is examined with delay, are true.
\ere

Moreover, we compare with the implicit method, see \cite{neuenkirch_szpruch:2014}.
Set $$G(x) = x - (1-q)\left(k_1x^{-q/(1-q)} - k_2x - q(k_3)^2x^{-1}/2   \right)\D$$ and compute
$$
y_{n+1} = G^{-1}(y_{n} + k_3(1-q)\D W_n)
$$
and then transform back to get the following scheme
\beqq\label{LSD-eq:Implicitscheme}
y_{n+1}^{Impl}= (y_{n+1})^{1/(1-q)}.
\eeqq

\begin{figure}[ht]
	\centering
	\begin{subfigure}{.47\textwidth}
		\includegraphics[width=1\textwidth]{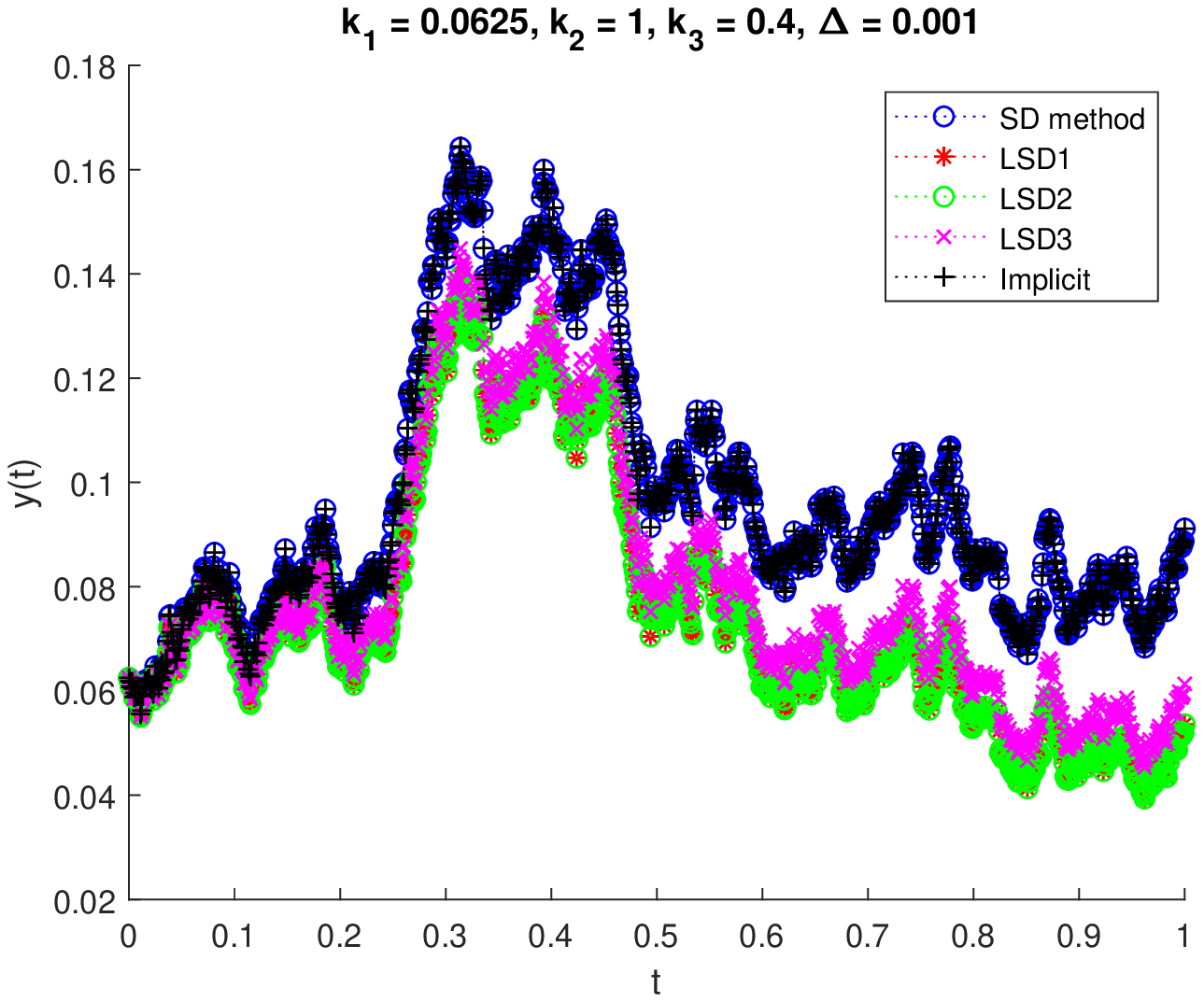}
		\caption{Trajectories  with $\D=10^{-3}$.}
	\end{subfigure}
	\begin{subfigure}{.47\textwidth}
		\includegraphics[width=1\textwidth]{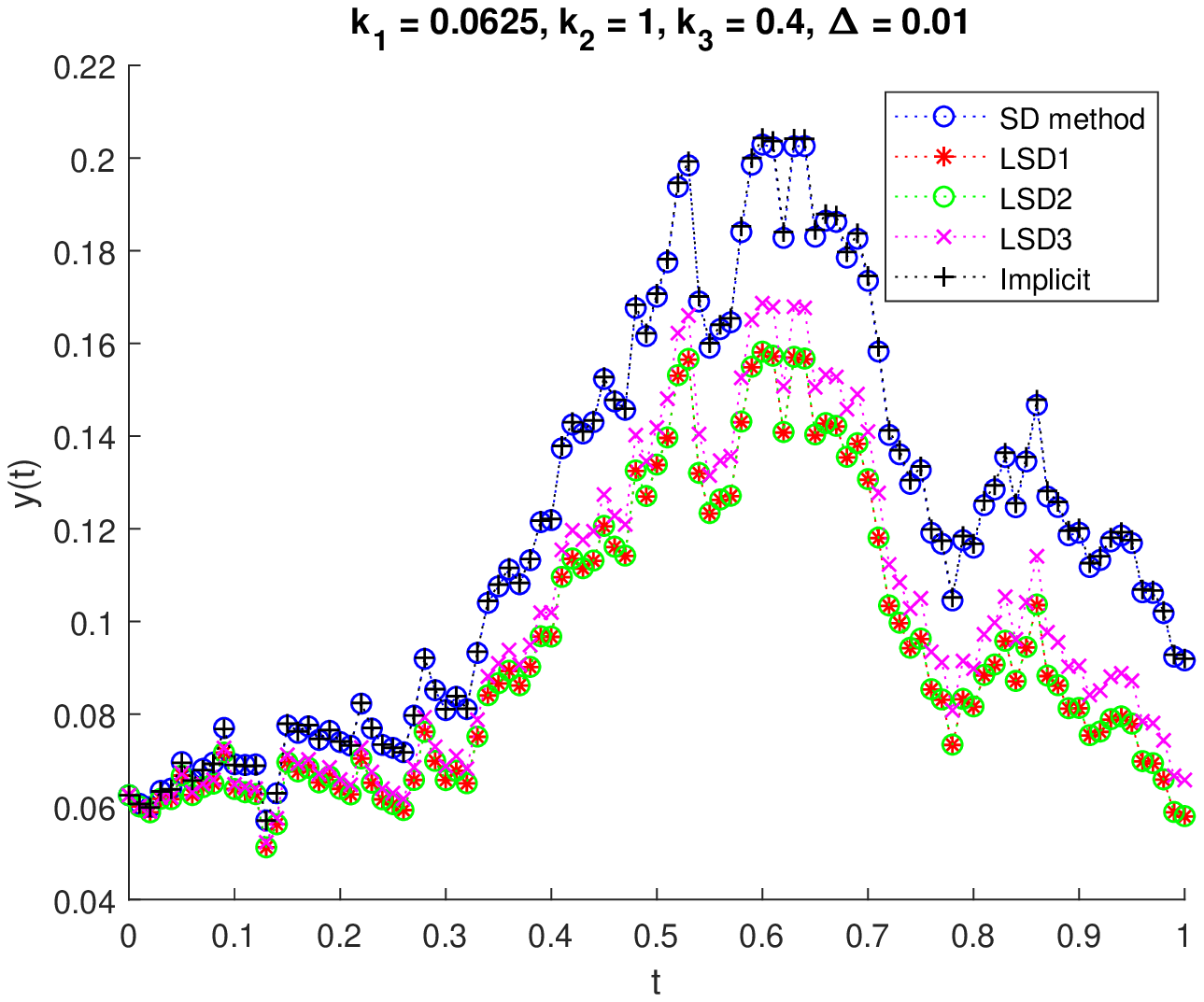}
		\caption{Trajectories  with $\D=10^{-2}$.}
	\end{subfigure}
	\caption{Trajectories  of  (\ref{LSD-eq:CEVLToriginal}), (\ref{LSD-eq:CEVLToriginal2}), (\ref{LSD-eq:CEVLToriginal3}), (\ref{LSD-eq:NCEV_SD}) and (\ref{LSD-eq:Implicitscheme}) for the approximation of (\ref{LSD-eq:CEV}) for different step-sizes.}\label{LSD-fig:LSDCEVSD}
\end{figure}

\begin{figure}[ht]
	\centering
	\begin{subfigure}{.47\textwidth}
		\includegraphics[width=1\textwidth]{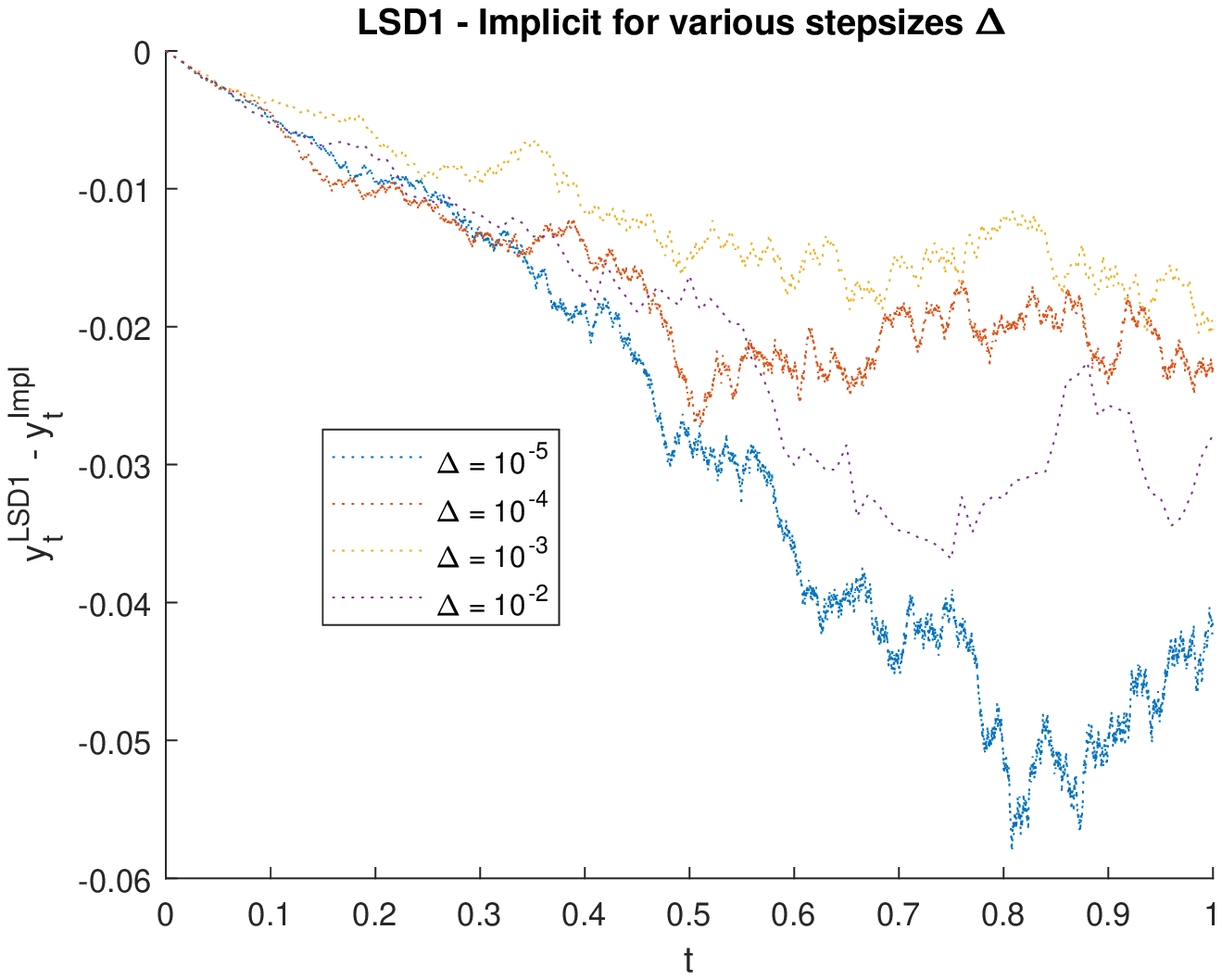}
		\caption{LSD1 - Implicit}
	\end{subfigure}
	\begin{subfigure}{.47\textwidth}
		\includegraphics[width=1\textwidth]{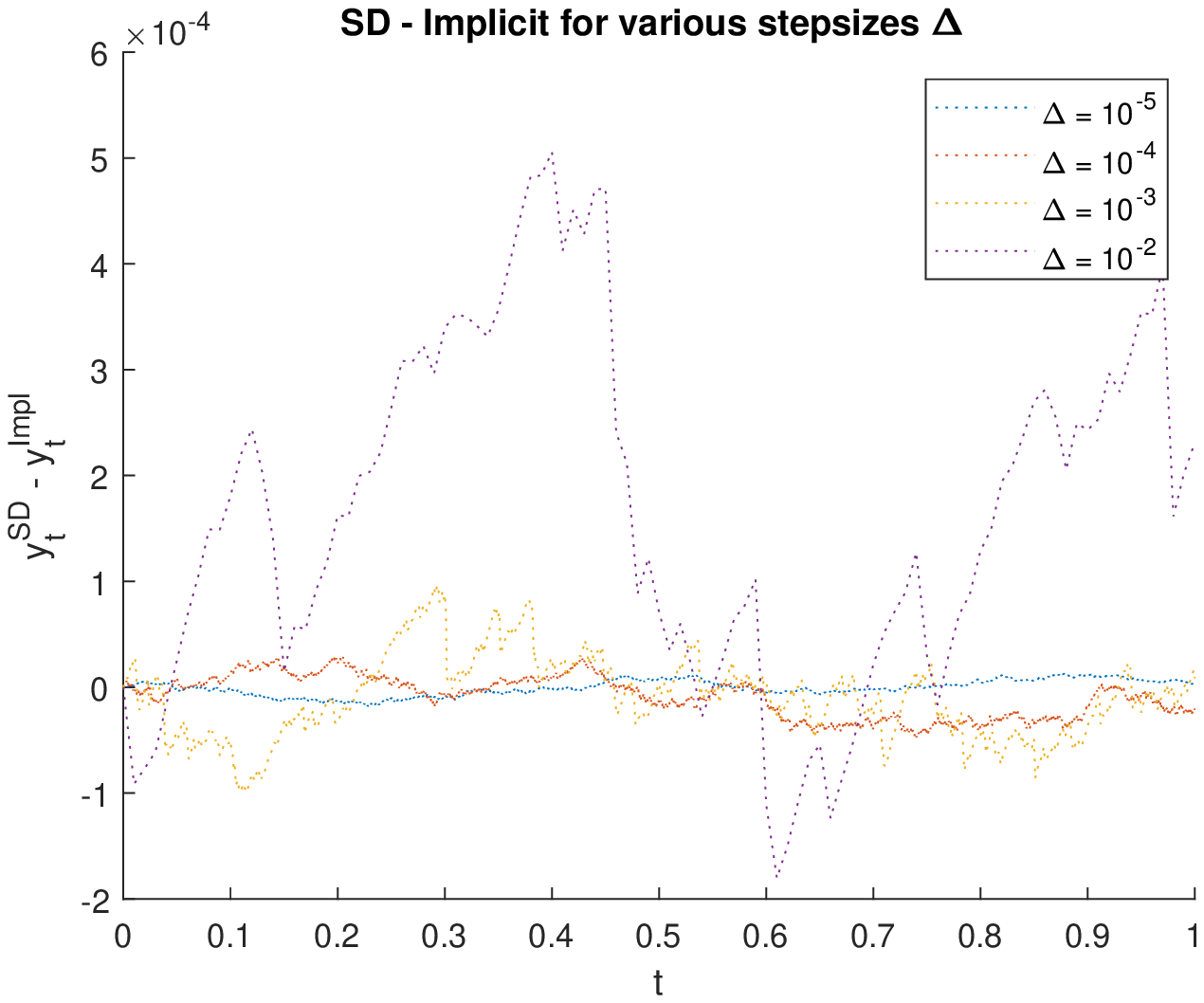}
		\caption{SD - Implicit}
	\end{subfigure}	
\begin{subfigure}{.47\textwidth}
	\includegraphics[width=1\textwidth]{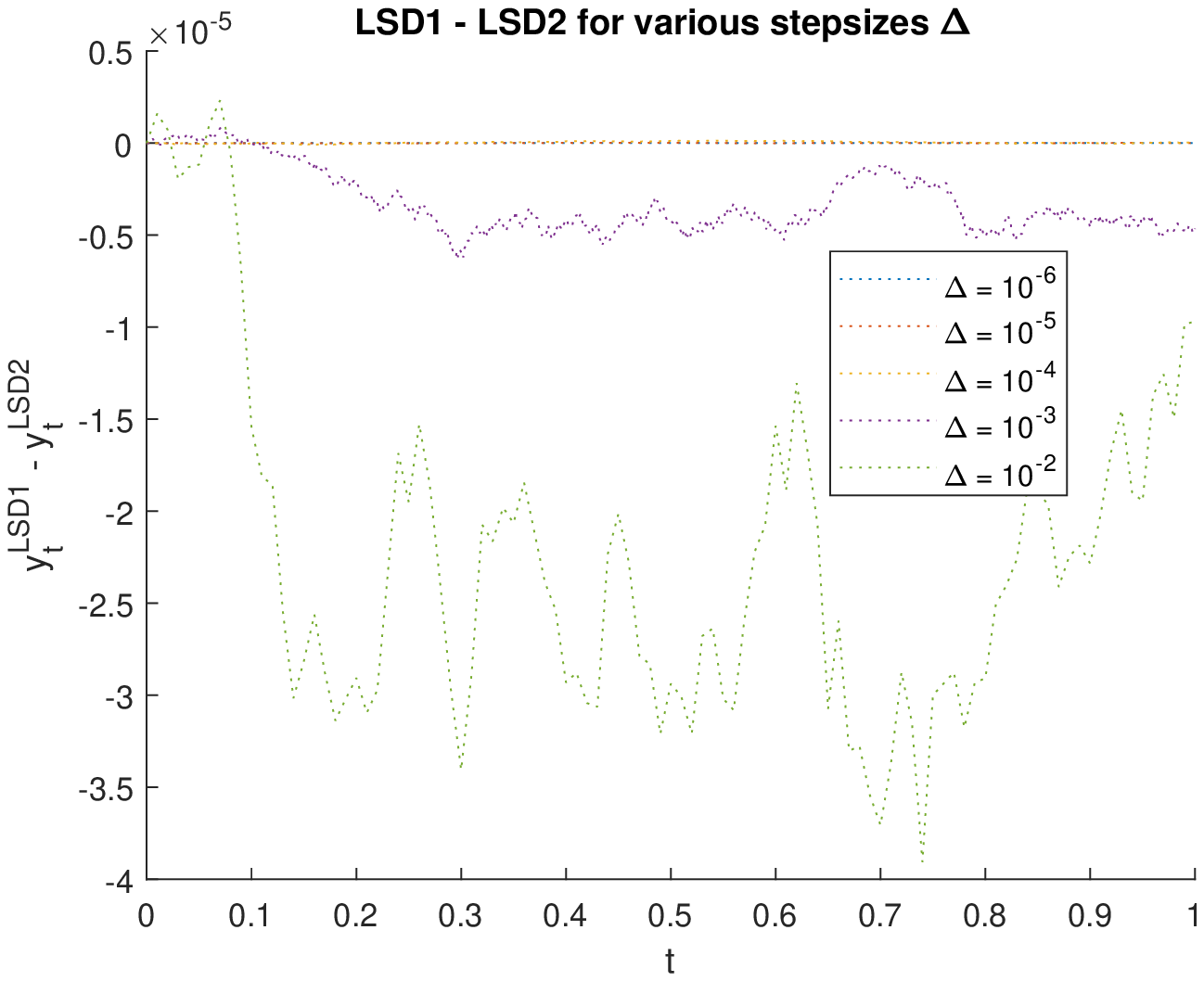}
	\caption{LSD1 - LSD2}
\end{subfigure}
	\begin{subfigure}{.47\textwidth}
	\includegraphics[width=1\textwidth]{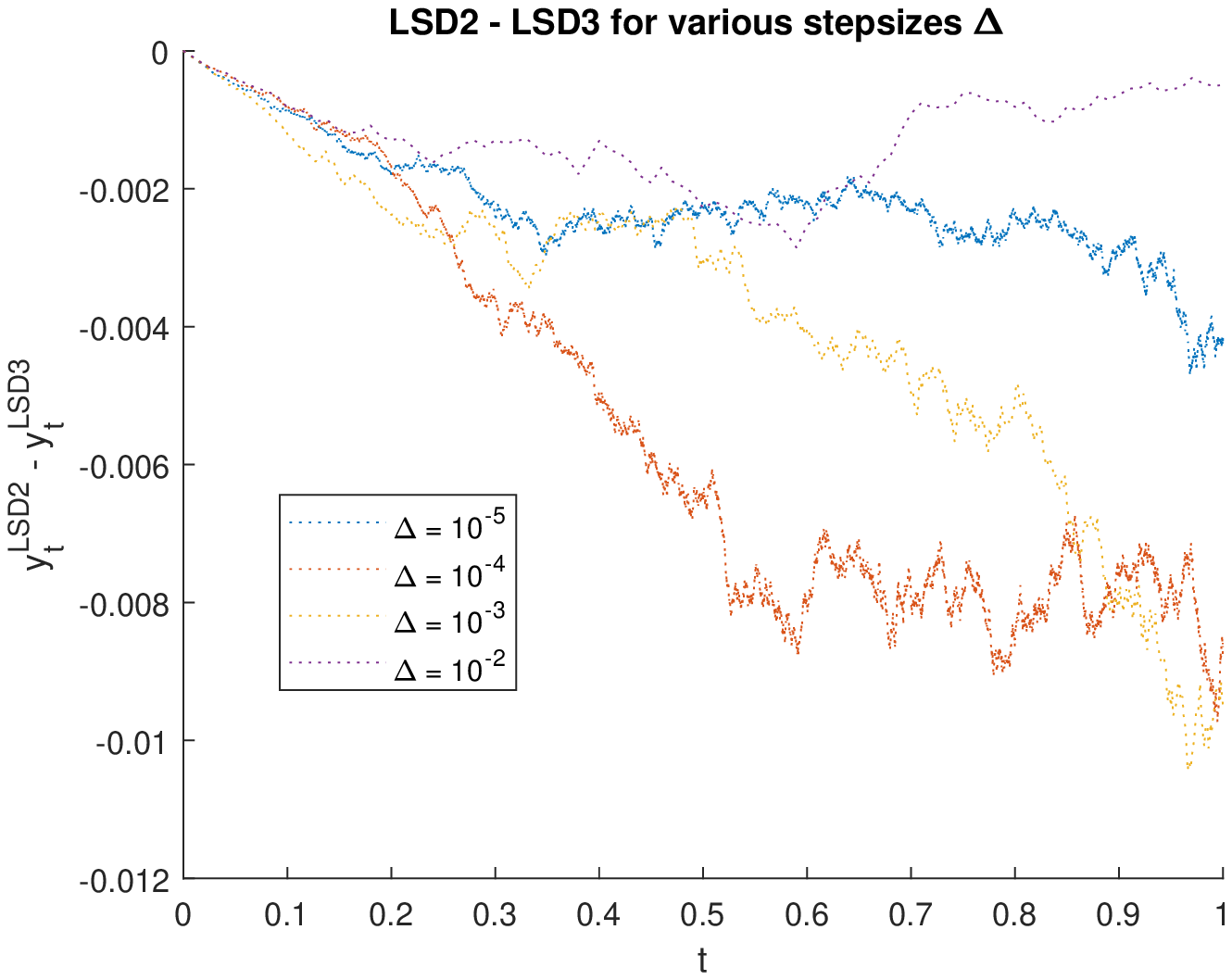}
	\caption{LSD2 - LSD3}
\end{subfigure}
	\caption{Trajectories of the differences of the numerical methods for the approximation of (\ref{LSD-eq:CEV}) with various step-sizes.}\label{LSD-fig:LSDCEVminSD}
\end{figure}

The simulation paths are presented in Figures \ref{LSD-fig:LSDCEVSD} and \ref{LSD-fig:LSDCEVminSD}.

\begin{figure}[ht]
	\centering
	\begin{subfigure}{.47\textwidth}
		\includegraphics[width=1\textwidth]{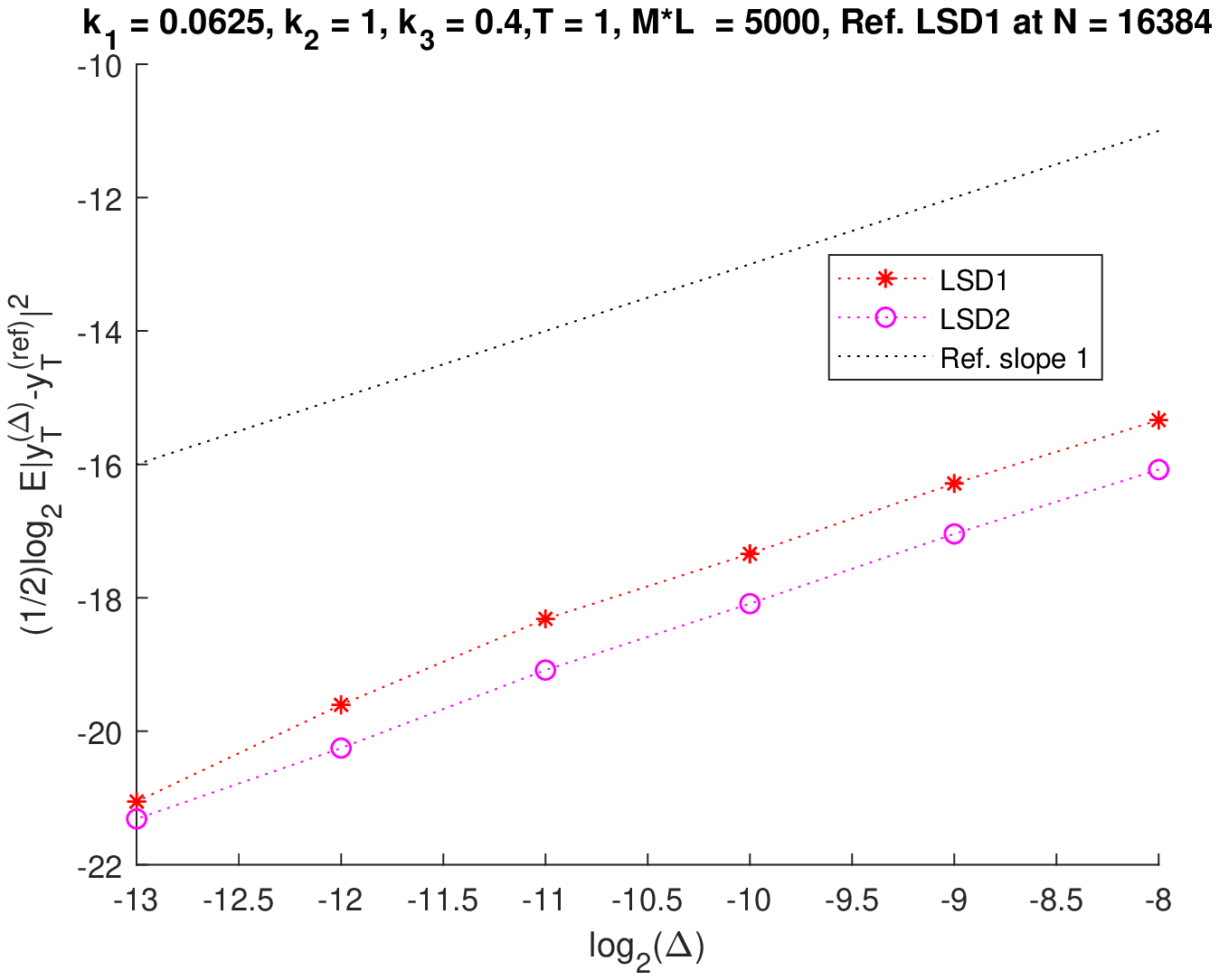}
		\caption{LSD1 as reference solution}
	\end{subfigure}
	\begin{subfigure}{.47\textwidth}
		\includegraphics[width=1\textwidth]{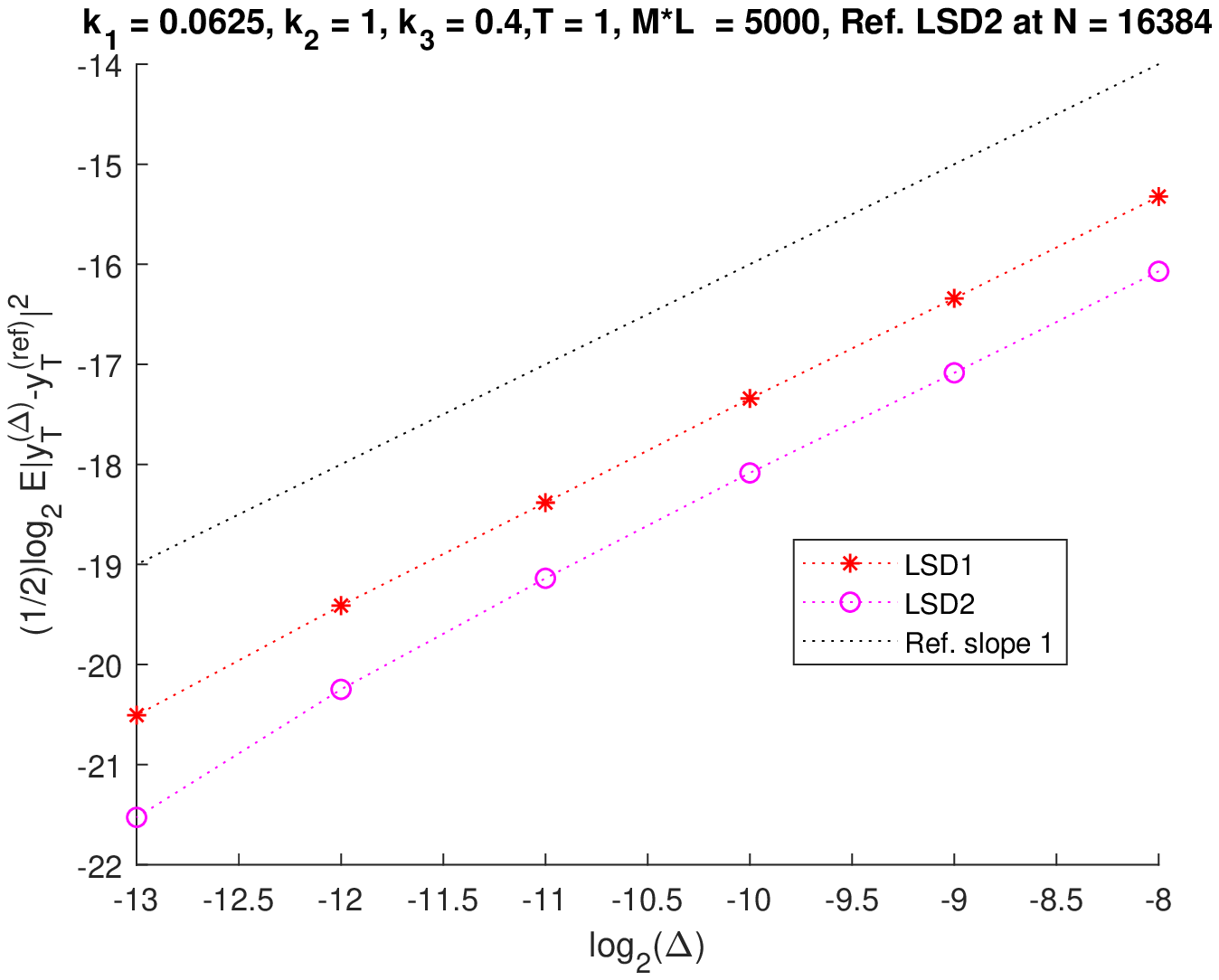}
		\caption{LSD2 as reference solution}
	\end{subfigure}
	\begin{subfigure}{.55\textwidth}
	\includegraphics[width=1\textwidth]{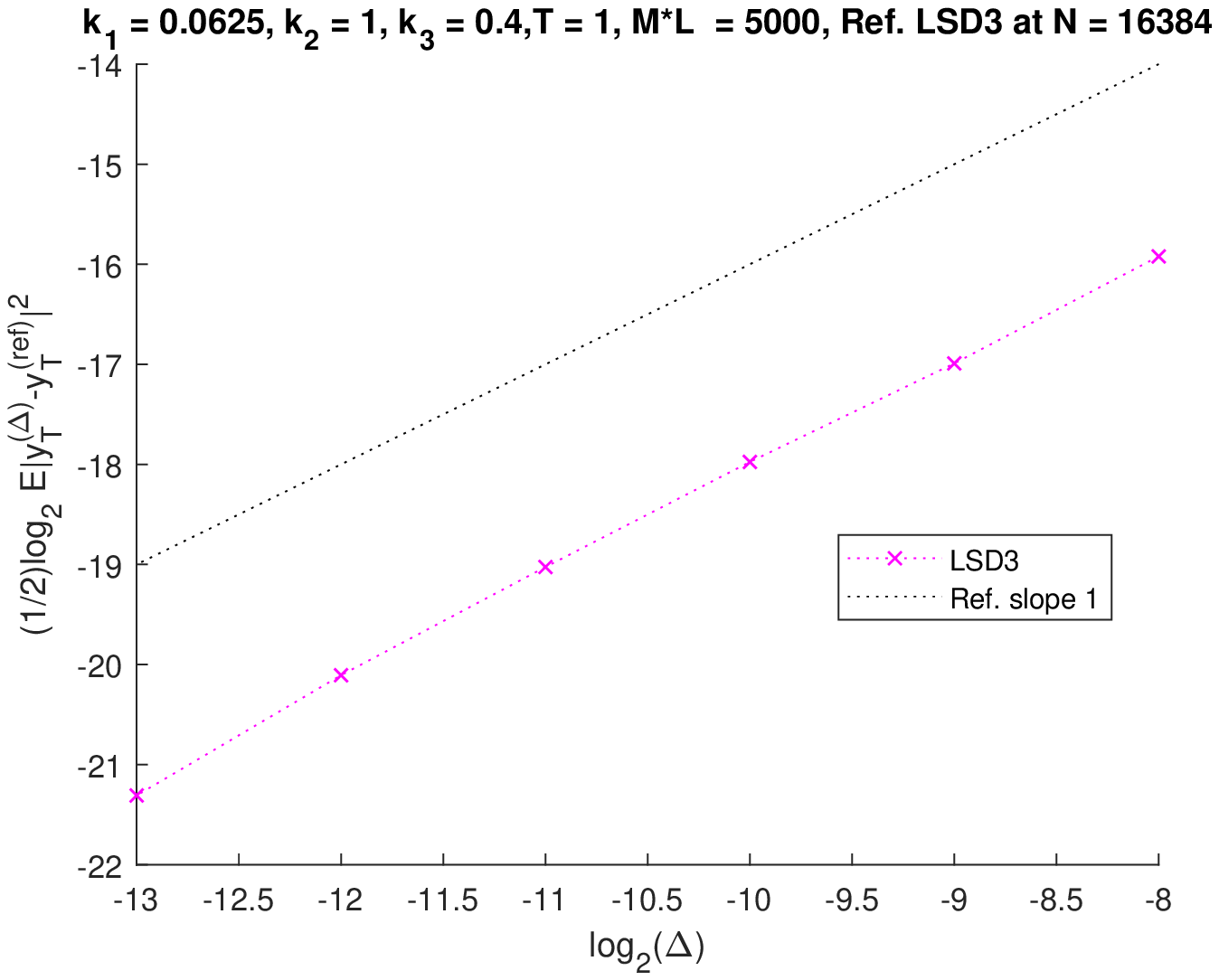}
	\caption{LSD3 as reference solution}
\end{subfigure}
	\caption{Convergence of (\ref{LSD-eq:CEVLToriginal}), (\ref{LSD-eq:CEVLToriginal2}) and (\ref{LSD-eq:CEVLToriginal3}) for the approximation of (\ref{LSD-eq:CEV}) with different reference solutions.}\label{LSD-fig:LSDCEVorder}
\end{figure}

We also examine numerically the order of strong convergence of the LSD methods. The numerical results suggest that the LSD schemes converge in the mean-square sense with order close to $1,$ see Figure \ref{LSD-fig:LSDCEVorder}. 

\section{Wright - Fisher model}\label{LSD:ssec:WF}

Let
\beqq  \label{LSD-eq:WFmodel}
x_t =x_0 + \int_0^t (k_1 - k_2x_s)ds + k_3\int_0^t \sqrt{x_s(1-x_s)}dW_s,
\eeqq
where $k_i>0, i=1,2,3.$ If $x_0\in(0,1)$ and $2k_1\geq (k_3)^2, 2(k_2-k_1)\geq (k_3)^2,$ then $0<x_t<1$ a.s., see \cite{karlin:1981}. SDE (\ref{LSD-eq:WFmodel}) appears in population dynamics to describe fluctuations in gene frequency of reproducing individuals among finite populations \cite{ewens:2012} and ion channel dynamics within cardiac and neuronal cells, (c.f. \cite{dangerfield:2012}, \cite{goldwyn:2011}, \cite{dangerfield:2012b} and references therein).

The transformed process $(z_t)$ of (\ref{LSD-eq:WFmodel}) with $z = 2\arcsin(\sqrt{x})$ has dynamics, 
\beam\nonumber  
z_t &=& z_0 + \int_0^t \left( \left(k_1 - \frac{(k_3)^2}{4}\right)\cot(\frac{z_s}{2}) - \left(k_2 - k_1 - \frac{(k_3)^2}{4}\right)\tan(\frac{z_s}{2})\right)ds\\
\label{LSD-eq:WFmodelLamperti}&& + k_3\int_0^t  dW_s.
\eeam
Set $a := k_1 - \frac{(k_3)^2}{4}$ and $b:=k_2 - k_1 - \frac{(k_3)^2}{4}.$ The conditions on the parameters imply $a>0$ and $b>0.$    We consider three versions of the semi-discrete method for approximating  (\ref{LSD-eq:CEVLamperti}). In the first version $(y_t),$ see Section \ref{LSD:subsec:WFz1}, we use the standard semi-discrete method and in the other two versions $(\hat{y}_t)$ and $(\wt{y}_t)$ we study the new versions  see Section \ref{LSD:subsec:WFz2z3}, where in each subinterval $(t_n, t_{n+1}]$ we solve an algebraic equation.

\subsection{Lamperti Semi-Discrete method $\wt{z}^1_{n}$ for Wright-Fisher}\label{LSD:subsec:WFz1}
	
Rewrite (\ref{LSD-eq:WFmodelLamperti}) as	

\beqq\label{LSD-eq:WFmodelLamperti_ab}
z_t = z_0 + \int_0^t \left( a\cot(\frac{z_s}{2}) - b\tan(\frac{z_s}{2})\right)ds + k_3\int_0^t  dW_s.
\eeqq
	
For $t\in(t_n, t_{n+1}]$ consider the process $(y_t)$ where  
\beam\nonumber
y_t &=& k_3\D W_n + y_{t_n} -b\tan(\frac{y_{t_n}}{2})\frac{y_{t}}{y_{t_n}}\D + \int_{t_n}^t a\cot(\frac{y_s}{2})ds\\
\label{LSD-eq:WFmodelLT}&=& \phi_\D(y_{t_n},\D W_n) + \int_{t_n}^t\frac{a}{1 + \frac{b}{y_{t_n}}\tan(\frac{y_{t_n}}{2})\D}\cot(\frac{y_s}{2}) ds,
\eeam
with $\hat{y}_{0}=z_0 $ and $\phi_\D(x, y) := \frac{k_3y + x}{1 + \frac{b}{x}\tan(x/2)\D}.$ 
Equation (\ref{LSD-eq:WFmodelLT}) has solution with the property, see Appendix \ref{LSD-ap:WFmodelLT_sol},
\beqq\label{LSD-eq:SD WFmodelLTsol2}
|\cos (\frac{y_t}{2}) | = |\cos (\frac{\phi_\D(y_{t_n},\D W_n)}{2}) | \exp\left\{-\frac{a/2}{1 + \frac{b}{y_{t_n}}\tan(\frac{y_{t_n}}{2})\D}    (t-t_n)\right\}.
\eeqq
Note that $0<z_t<\pi$ a.s therefore since $(y_t)$ converges to $(z_t)$  the process $(y_t)/2$ and belong to $(0,\pi/2).$
The proposed semi-discrete method $(y_t)$ for the approximation of (\ref{LSD-eq:WFmodelLamperti}) satisfies  
\beqq\label{LSD-eq:WFmodelLT_transf}
\cos (\frac{y_{t_{n+1}}}{2})  = |\cos (\frac{\phi_\D(y_{t_n},\D W_n)}{2})| \exp\left\{-\frac{a\D/2}{1 + \frac{b}{y_{t_n}}\tan(\frac{y_{t_n}}{2})\D}\right\}, 
\eeqq
which suggests the Lamperti semi-discrete method $(\wt{z}^1_n)_{n\in\bbN}$  for the approximation of (\ref{LSD-eq:WFmodel})
\beqq\label{LSD-eq:WFmodeloriginal}
\wt{z}^1_{t_{n+1}} = 1 - \cos^2(\frac{\phi_\D(y_{t_n},\D W_n)}{2})\exp\left\{-\frac{a\D}{1 + \frac{b}{y_{t_n}}\tan(\frac{y_{t_n}}{2})\D}\right\}.
\eeqq
Note that $(\wt{z}^1_{t_{n}})_{n\in\bbN}\in(0,1)$ when $x_0\in(0,1).$

\subsection{Lamperti Semi-Discrete methods $\wt{z}^2_{n}, \wt{z}^3_{n}$ and $\wt{z}^4_{n}$  for Wright-Fisher}\label{LSD:subsec:WFz2z3}

For $t\in(t_n, t_{n+1}]$ consider the processes $(\hat{y}_t), (\wt{y}_t)$ and $(\bar{y}_t)$ where  
\beam\nonumber
\hat{y}_t &=& k_3(W_t - W_{t_n}) + \hat{y}_{t_n} + \hat{y}_{t}\left(\frac{a}{\hat{y}_{t_n}}\cot(\frac{\hat{y}_{t_n}}{2}) -\frac{b}{\hat{y}_{t_n}}\tan(\frac{\hat{y}_{t_n}}{2})\right)\D \\
\label{LSD-eq:WFmodelLT2}&=& \frac{k_3(W_t - W_{t_n}) + \hat{y}_{t_n}}{1 -\frac{a}{\hat{y}_{t_n}}\cot(\frac{\hat{y}_{t_n}}{2})\D +\frac{b}{\hat{y}_{t_n}}\tan(\frac{\hat{y}_{t_n}}{2})\D},
\eeam
with $\hat{y}_{0}=z_0,$

\beam\nonumber
\wt{y}_t &=& k_3(W_t - W_{t_n}) + \wt{y}_{t_n} + a\cot(\frac{\wt{y}_{t_n}}{2})\D -\frac{b}{\wt{y}_{t_n}}\tan(\frac{\wt{y}_{t_n}}{2})\D\wt{y}_t \\
\label{LSD-eq:WFmodelLT3}&=& \frac{k_3(W_t - W_{t_n}) + \wt{y}_{t_n} + a\cot(\frac{\wt{y}_{t_n}}{2})\D}{1  + \frac{b}{\wt{y}_{t_n}}\tan(\frac{\wt{y}_{t_n}}{2})\D},
\eeam
with $\wt{y}_{0}=z_0$ and
$$
\bar{y}_t = k_3(W_t - W_{t_n}) + \bar{y}_{t_n} + \left(a\cot(\frac{\bar{y}_{t_n}}{2}) -b\tan(\frac{\bar{y}_{t_n}}{2}) - \frac{1}{\bar{y}_{t_n}}\right)\D + \frac{\D}{\bar{y}_{t}}\\
$$
with $\bar{y}_{0}=z_0.$ The solution of $(\bar{y}_t)$ of the above equation satisfies
\beqq\label{LSD-eq:WFmodelLT4} 
(\bar{y}_t)^2 - \left(k_3(W_t - W_{t_n}) + \bar{y}_{t_n} - \frac{\D}{\bar{y}_{t_n}} + \left(a\cot(\frac{\bar{y}_{t_n}}{2}) -b\tan(\frac{\bar{y}_{t_n}}{2})\right)\D\right)(\bar{y}_t) -\D =0.
\eeqq

The proposed semi-discrete methods $(\hat{y}_t)$ and $(\wt{y}_t)$ for the approximation of (\ref{LSD-eq:WFmodelLamperti}) read

\beqq\label{LSD-eq:WFmodelLT2discr}
\hat{y}_{t_{n+1}} = \frac{k_3\D W_{n} + \hat{y}_{t_n}}{1 -\frac{a}{\hat{y}_{t_n}}\cot(\frac{\hat{y}_{t_n}}{2})\D +\frac{b}{\hat{y}_{t_n}}\tan(\frac{\hat{y}_{t_n}}{2})\D}
\eeqq
and
\beqq\label{LSD-eq:WFmodelLT3discr}
\wt{y}_{t_{n+1}} =  \frac{k_3\D W_n + \wt{y}_{t_n} + a\cot(\frac{\wt{y}_{t_n}}{2})\D}{1  + \frac{b}{\wt{y}_{t_n}}\tan(\frac{\wt{y}_{t_n}}{2})\D},
\eeqq
respectively, while for $(\bar{y}_t)$ we have that
\beqq\label{LSD-eq:WFmodelLT4discr}
\bar{y}_{t_{n+1}} = \frac{\bar{\phi}_\D(\bar{y}_{t_n},\D W_n) + \sqrt{\bar{\phi}^2_\D(\bar{y}_{t_n},\D W_n) + 4\D}}{2},
\eeqq
with $\bar{\phi}_\D(x,y) = (k_3y + x - \D x^{-1} + a\left(\cot(x/2) -b\tan(x/2)\right)\D.$ Therefore, the versions of the Lamperti semi-discrete method $(\wt{z}^2_n)_{n\in\bbN}, (\wt{z}^3_n)_{n\in\bbN}$ and $(\wt{z}^4_n)_{n\in\bbN}$ for the approximation of (\ref{LSD-eq:WFmodel}) are 
\beqq\label{LSD-eq:WFmodeloriginal2}
\wt{z}^2_{t_{n+1}} = \sin^2\left(\frac{k_3\D W_{n} + \hat{y}_{t_n}}{2 -\frac{2a}{\hat{y}_{t_n}}\cot(\frac{\hat{y}_{t_n}}{2})\D +\frac{2b}{\hat{y}_{t_n}}\tan(\frac{\hat{y}_{t_n}}{2})\D}\right),
\eeqq
\beqq\label{LSD-eq:WFmodeloriginal3}
\wt{z}^3_{t_{n+1}} = \sin^2\left(\frac{k_3\D W_{n} + \wt{y}_{t_n}  + a\cot(\frac{\wt{y}_{t_n}}{2})\D}{2 + \frac{2b}{\wt{y}_{t_n}}\tan(\frac{\wt{y}_{t_n}}{2})\D}\right)
\eeqq
and 
\beqq\label{LSD-eq:WFmodeloriginal4}
\wt{z}^4_{t_{n+1}} = \sin^2\left(\frac{\bar{\phi}_\D(\bar{y}_{t_n},\D W_n) + \sqrt{\bar{\phi}^2_\D(\bar{y}_{t_n},\D W_n) + 4\D}}{4}\right),
\eeqq
respectively. Note that $(\wt{z}^j_{t_{n}})_{n\in\bbN}\in(0,1), j=2,3,4$ when $x_0\in(0,1).$

\subsection{Numerical experiment for Wright-Fisher}\label{LSD:subsec:WFnum}

The semi-discrete method we proposed in \cite{stamatiou:2018} reads
\beqq\label{LSD-eq:WFmodelSD}
y_{t_{n+1}}=\sin^2 \left(\frac{k_3}{2}\D W_n + \arcsin(\sqrt{\un{y_{t_n} + \left(k_1 - \frac{(k_3)^2}{4} + y_{t_n}\left(\frac{(k_3)^2}{2}-k_2\right)\right)\cdot\D}_{y_n}})\right),
\eeqq
which also possesses the qualitative property of domain preservation. Method (\ref{LSD-eq:WFmodelSD}) is well defined for all sufficiently small $\D$ such that $0<y_n <1.$ Let $\beta:=\frac{(k_3)^2}{2}-k_2$ with $\beta<0.$ We require $\D$ small enough so that $0<y_{t_n}(1+\beta\D) + a\D.$ To simplify the conditions on $a, \beta, \D$, when necessary, we may adopt the strategy presented in \cite{stamatiou:2018} and consider the SD method
\beqq\label{LSD-eq:SD_processNWalt}
\bar{y}_{t_{n+1}} = \sin^2 \left(\frac{k_3}{2}\D W_n + \arcsin(\sqrt{\bar{y}_n})\right),
\eeqq
with 
$$
\bar{y}_n:=\frac{\bar{y}_{t_n}(1 + \beta \D) + a\D}{1  + (a + \beta)\D}.
$$
The numerical scheme (\ref{LSD-eq:SD_processNWalt}) is mean square convergent when $(k_3)^2<2k_2$ for $\D<-1/\beta,$ see \cite[Prop. 2.5]{stamatiou:2018}  where the order of strong convergence was not theoretically proved. 

The Balance Implicit Split Step (BISS) method suggested in \cite[(4.8)]{dangerfield:2012} reads 
\beqq\label{LSD-eq:BISSscheme}
y_{n+1}^{BISS}=y_n + (k_1-k_2y_n)\D + \frac{C\sqrt{y_n(1-y_n)}\D W_n }{1+ d^1(y_n)|\D W_n |}(1-k_2\D), 
\eeqq
where $\D$ is the step-size of the equidistant discretization of the interval $[0,1]$, the control function $d^1$ is given by
$$
d^1(y) = \begin{cases}
k_3\sqrt{(1-\vep)/\vep} \quad \,\text{if}\,\, y< \vep,\\        
k_3\sqrt{(1-y)/y} \quad \text{if}\,\, \vep\leq y<1/2,\\
k_3\sqrt{y/(1-y)} \quad \text{if}\,\, 1/2\leq y\leq 1-\vep,\\  
k_3\sqrt{(1-\vep)/\vep} \quad \,\text{if}\,\, y>1 - \vep,
\end{cases}
$$
and 
$$
\vep= \min \{ k_1\D, (k_2-k_1)\D, 1-k_1\D, 1-(k_2-k_1)\D\}.
$$

The hybrid (HYB)  scheme as proposed in \cite[(11)]{dangerfield:2012c} is the result of a splitting method and reads
\beqq\label{LSD-eq:Hybscheme}
y_{n+1}^{HYB}=\frac{a}{\beta}(e^{\beta\D}-1) + e^{\beta\D}\sin^2 \left(\frac{k_3}{2}\D W_n + \arcsin(\sqrt{y_n})\right).
\eeqq
with the restriction that 
$$\frac{k_1}{k_2}\in\left(\frac{(k_3)^2}{4k_2},1-\frac{(k_3)^2}{4k_2}\right).$$ 
Moreover, we compare with the implicit method, see \cite{neuenkirch_szpruch:2014}.
Set $$G(x) = x - a\cot(x/2)\D - b(\tan(x/2))\D$$ and compute
$$
y_{n+1} = G^{-1}(y_{n} + k_3\D W_n)
$$
and then transform back to get the following scheme
\beqq\label{LSD-eq:ImplicitschemeWF}
y_{n+1}^{Impl}= \sin^2(y_{n+1}/2),
\eeqq

We use the set of parameters from \cite[Sec.4]{stamatiou:2018} where all the methods work well, i.e. we take $(k_1, k_2, k_3)=(1, 2, 0.20101),$ with $T = 1$ and various step-sizes and compare the proposed versions of LSD schemes (\ref{LSD-eq:WFmodeloriginal}), (\ref{LSD-eq:WFmodeloriginal2}), (\ref{LSD-eq:WFmodeloriginal3}) and (\ref{LSD-eq:WFmodeloriginal4}) with the BISS, the HYB, the SD method (\ref{LSD-eq:WFmodelSD}) and the implicit method (\ref{LSD-eq:ImplicitschemeWF}). The initial condition is chosen to be the steady state of the deterministic part, i.e. $x_0 = k_1/k_2.$ 

\begin{figure}[ht]
	\centering
	\begin{subfigure}{.47\textwidth}
		\includegraphics[width=1\textwidth]{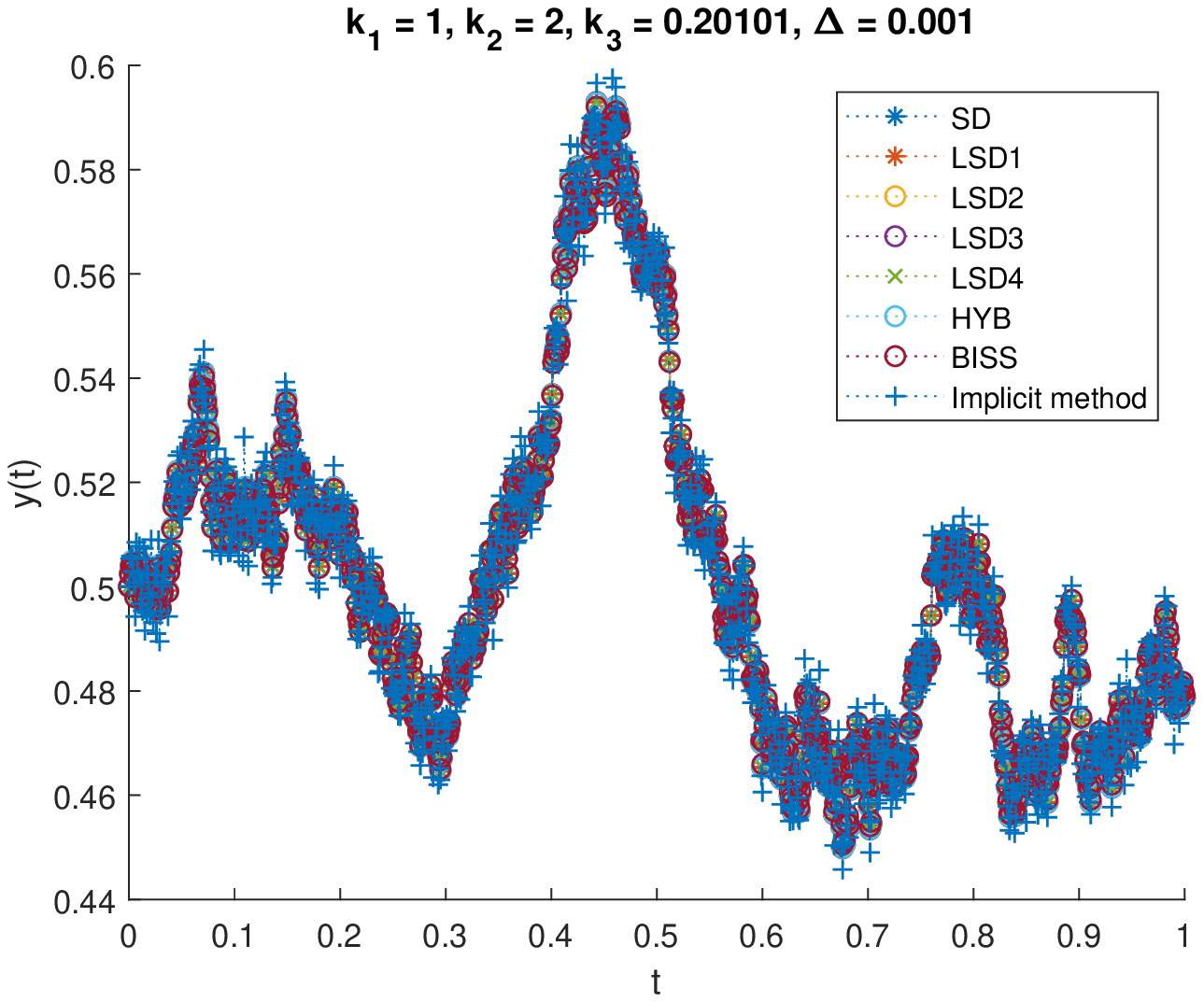}
		\caption{With $\D = 10^{-3}$.}
	\end{subfigure}
	\begin{subfigure}{.47\textwidth}
		\includegraphics[width=1\textwidth]{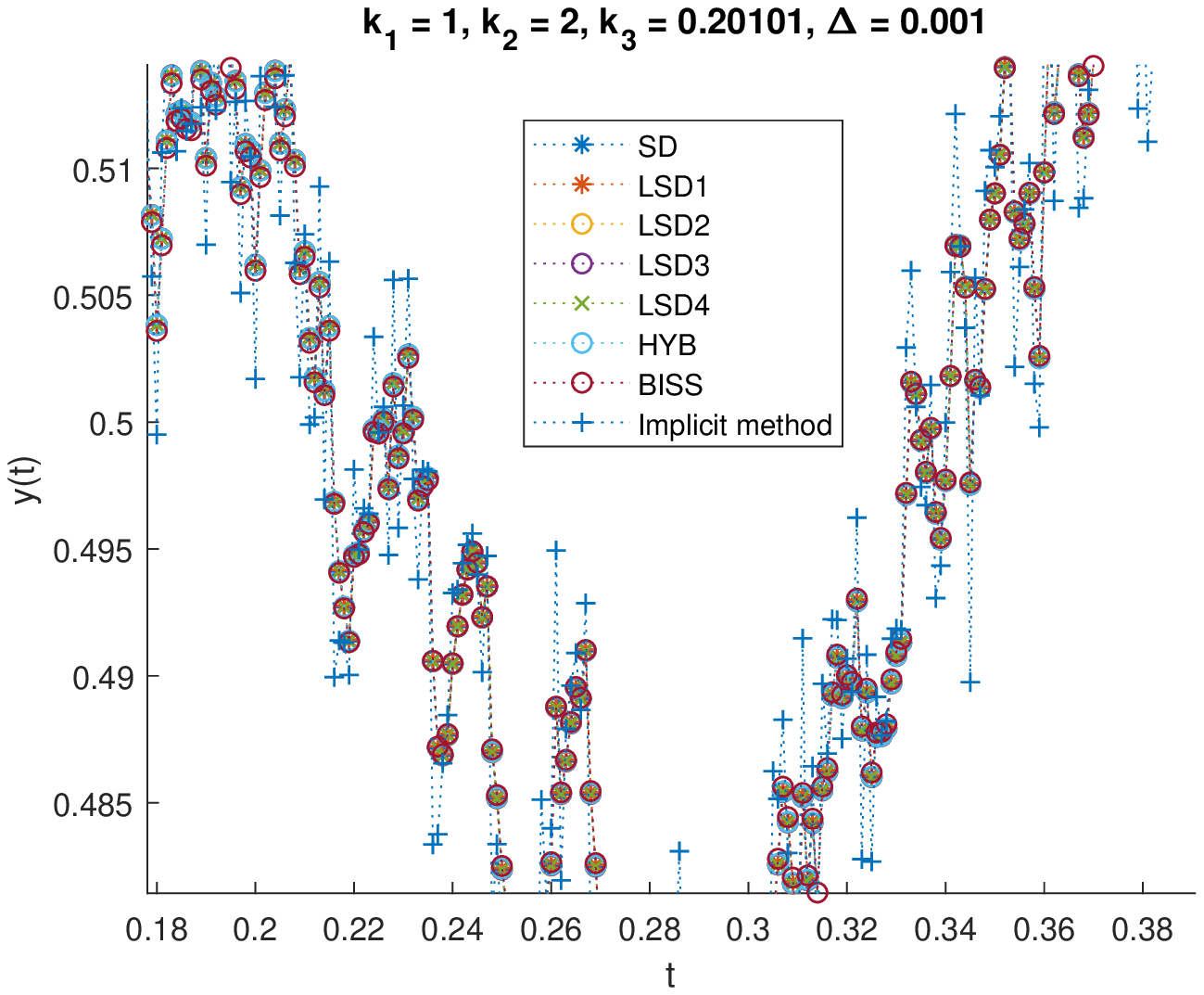}
		\caption{Zoom of Figure \ref{LSD-fig:LSDWFSDI}(A).}
	\end{subfigure}
	\caption{Trajectories  of (\ref{LSD-eq:WFmodeloriginal}), (\ref{LSD-eq:WFmodeloriginal2}), (\ref{LSD-eq:WFmodeloriginal3}), (\ref{LSD-eq:WFmodeloriginal4}), (\ref{LSD-eq:WFmodelSD}), (\ref{LSD-eq:BISSscheme}), (\ref{LSD-eq:Hybscheme}) and (\ref{LSD-eq:ImplicitschemeWF}) for the approximation of (\ref{LSD-eq:WFmodel}).}\label{LSD-fig:LSDWFSDI}
\end{figure}

\begin{figure}[h]
	\centering
	\begin{subfigure}{.3\textwidth}
		\includegraphics[width=1\textwidth]{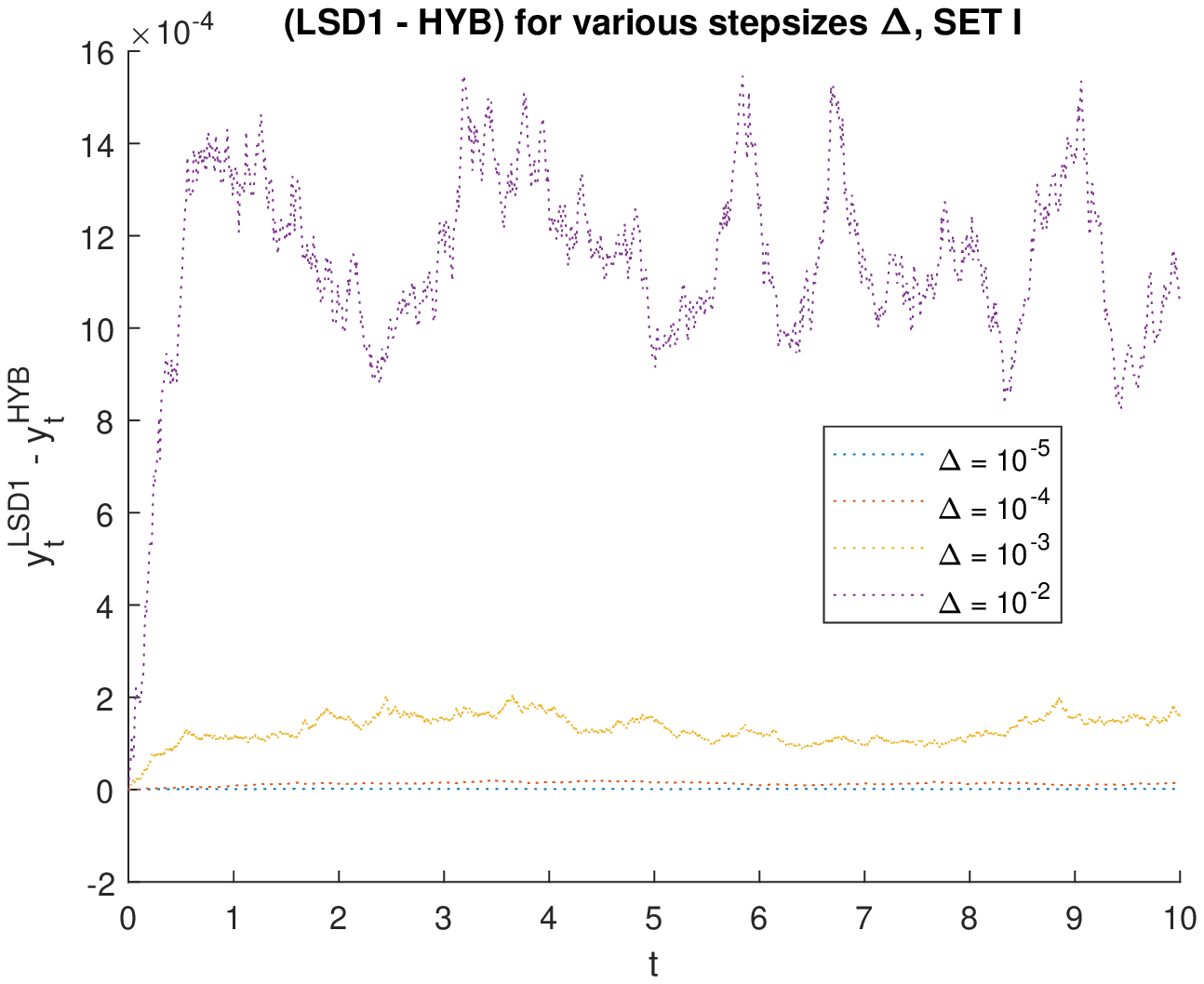}
		\caption{LSD1 - HYB}
	\end{subfigure}
	\begin{subfigure}{.3\textwidth}
	\includegraphics[width=1\textwidth]{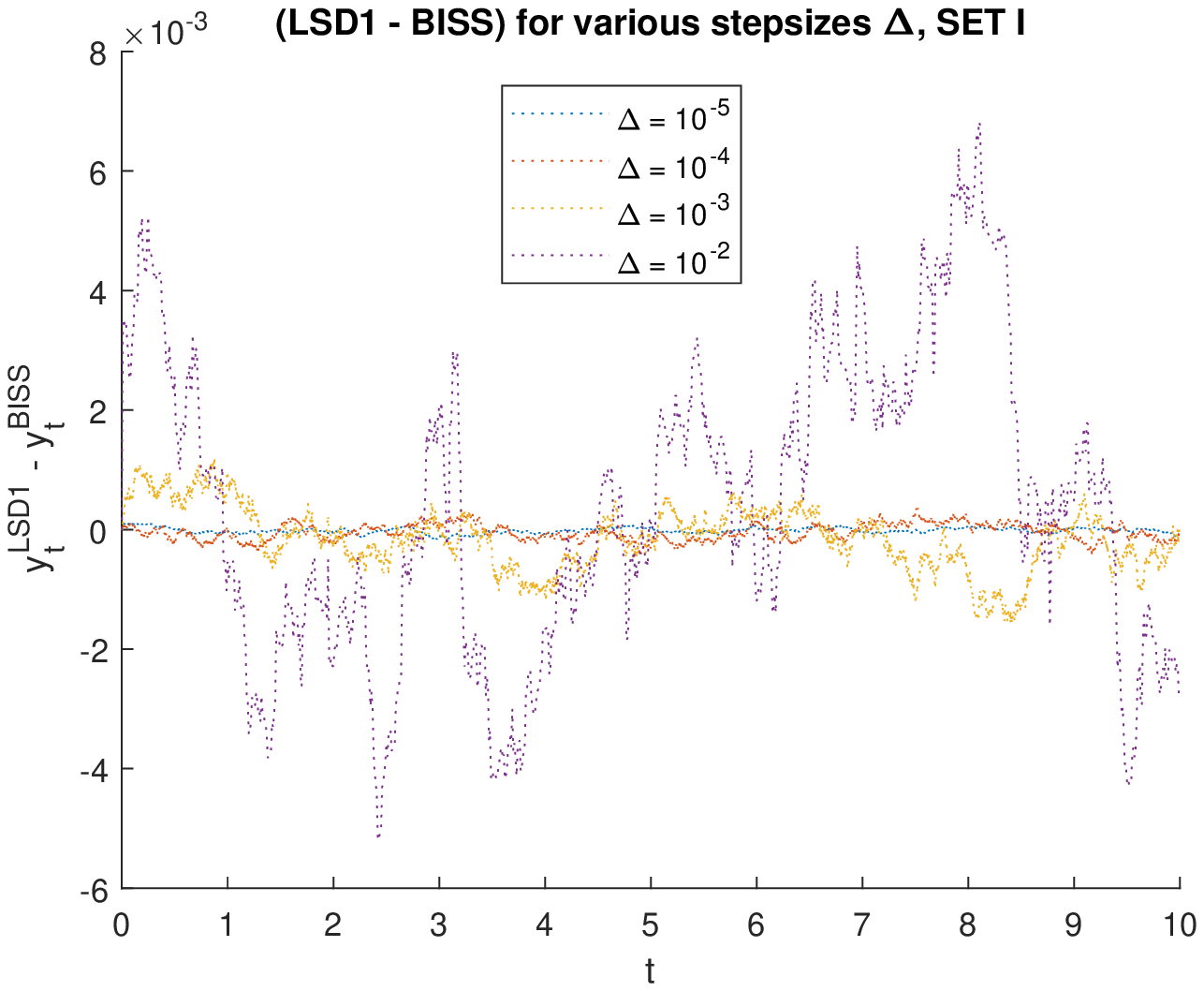}
	\caption{LSD1 - BISS}
\end{subfigure}
	\begin{subfigure}{.3\textwidth}
		\includegraphics[width=1\textwidth]{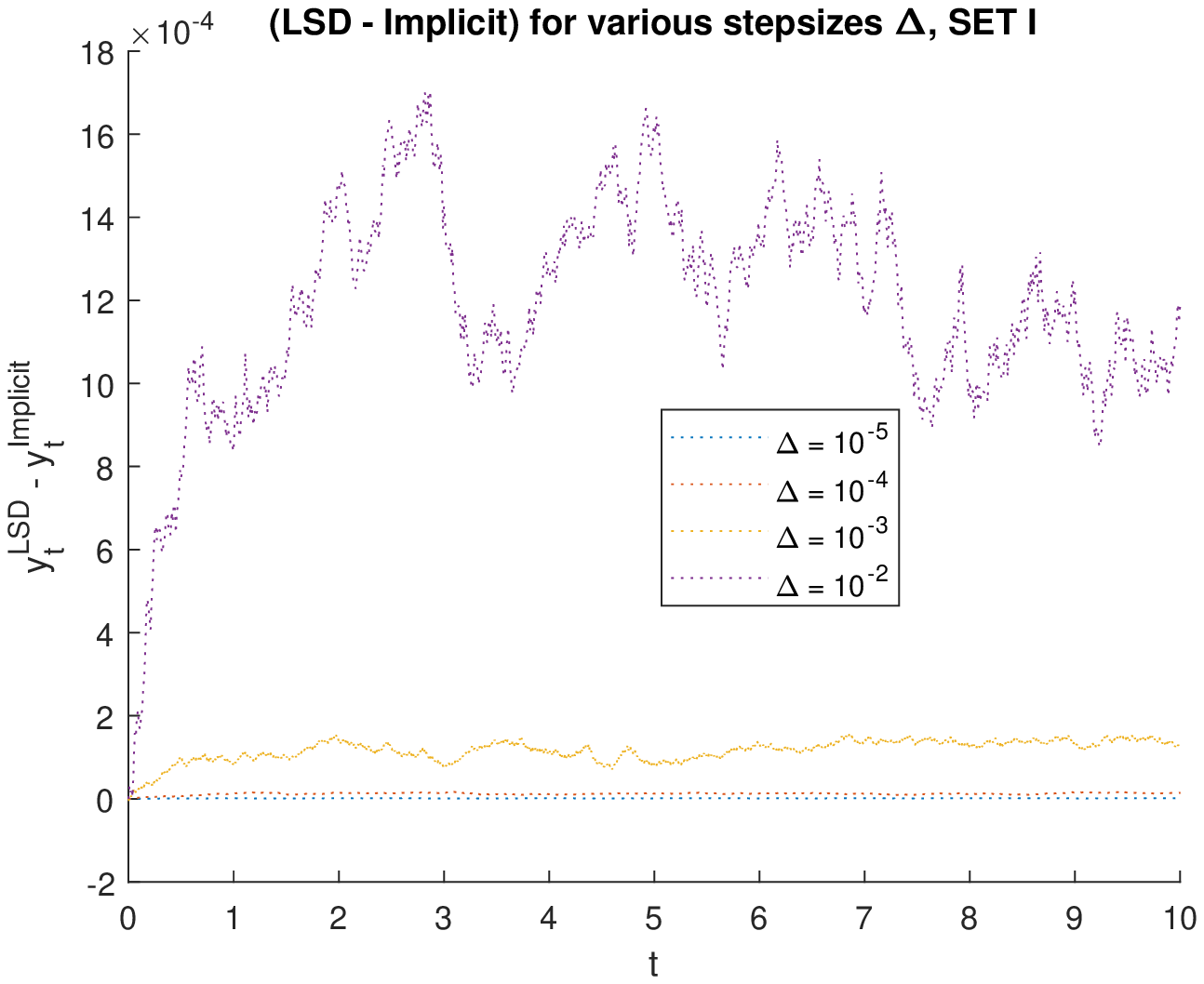}
		\caption{LSD1 - Implicit}
	\end{subfigure}
	\begin{subfigure}{.3\textwidth}
	\includegraphics[width=1\textwidth]{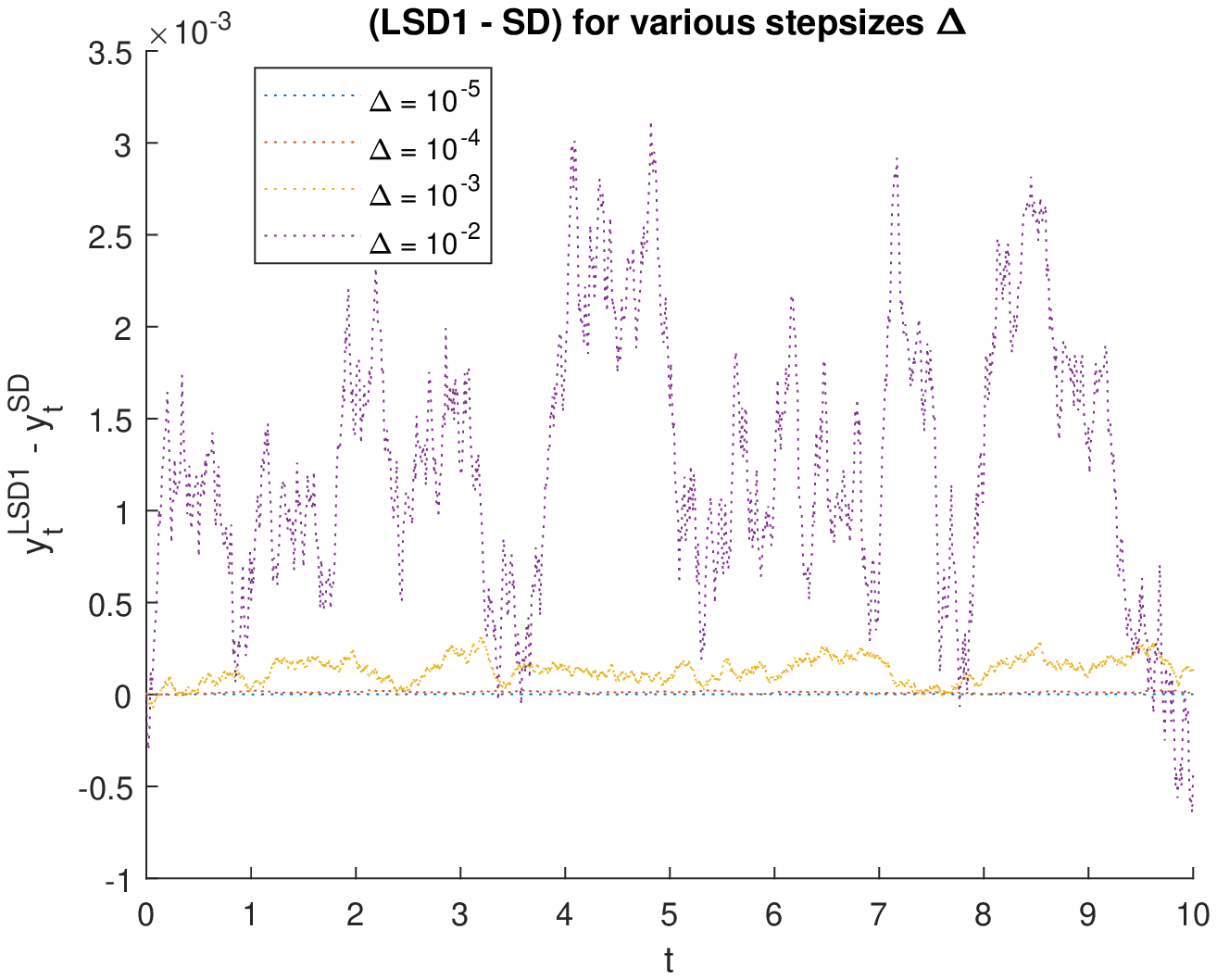}
	\caption{LSD1 - SD}
\end{subfigure}
\begin{subfigure}{.3\textwidth}
	\includegraphics[width=1\textwidth]{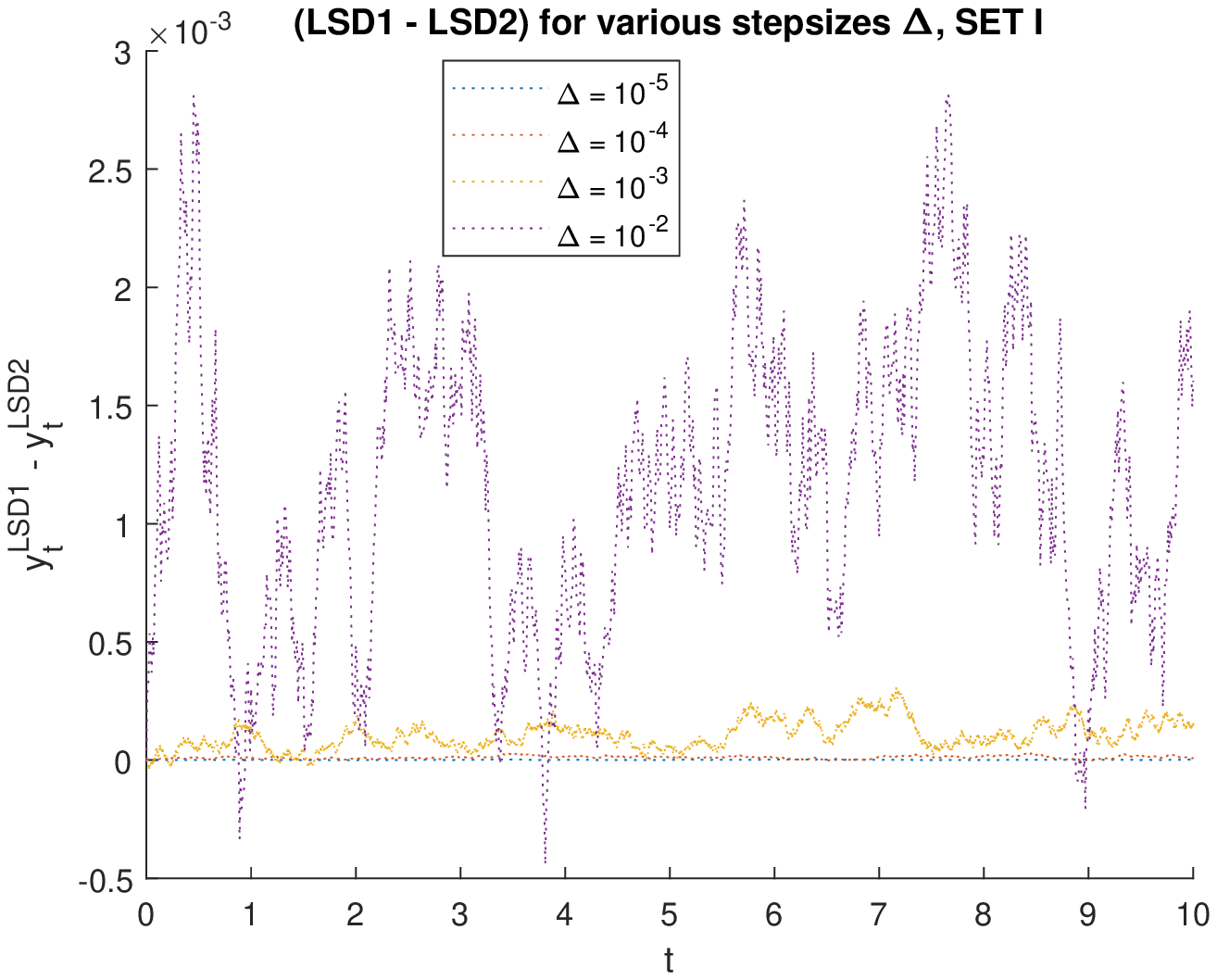}
	\caption{LSD1 - LSD2}
\end{subfigure}
	\begin{subfigure}{.3\textwidth}
	\includegraphics[width=1\textwidth]{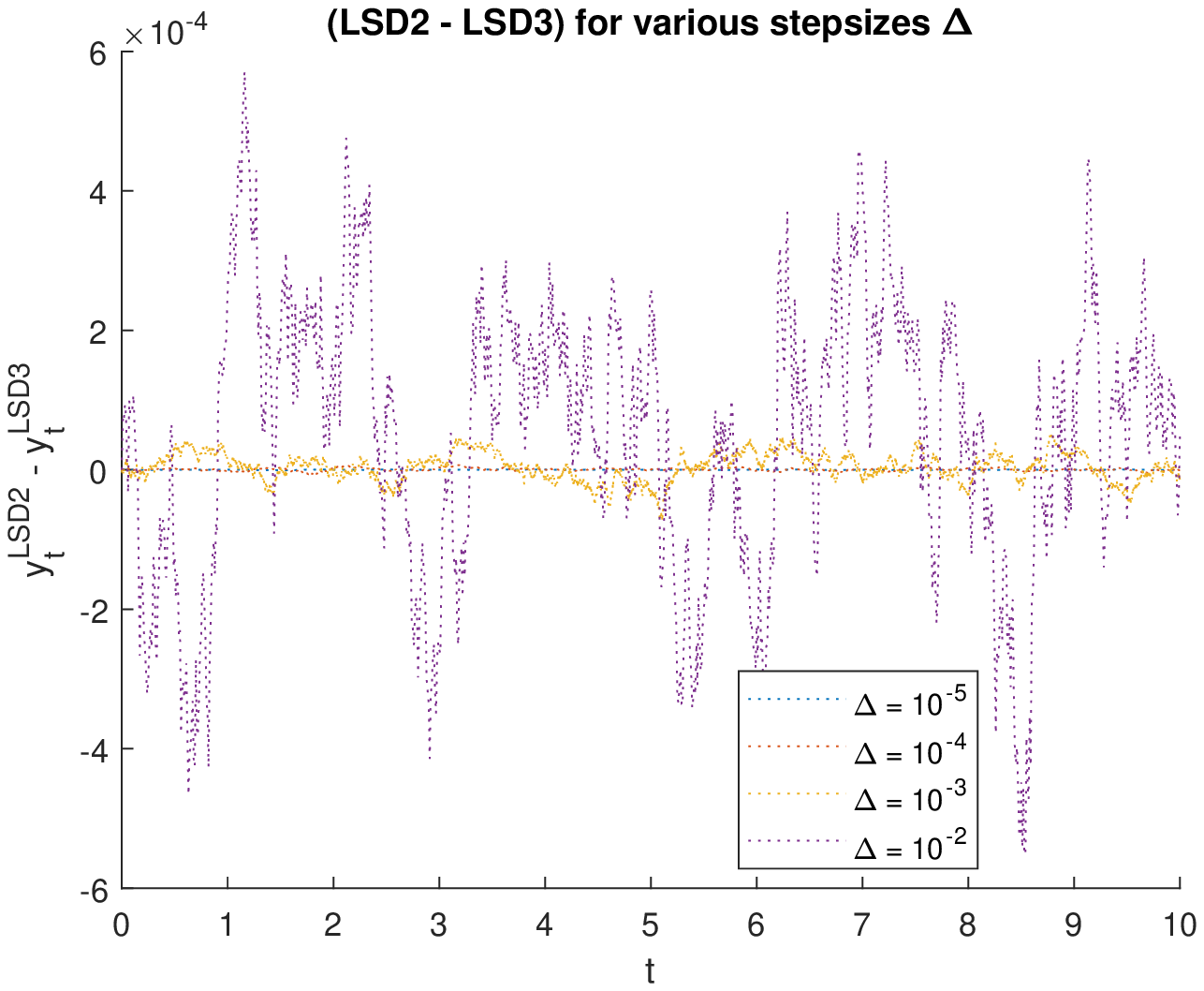}
	\caption{LSD2 - LSD3}
\end{subfigure}
\begin{subfigure}{.4\textwidth}
	\includegraphics[width=1\textwidth]{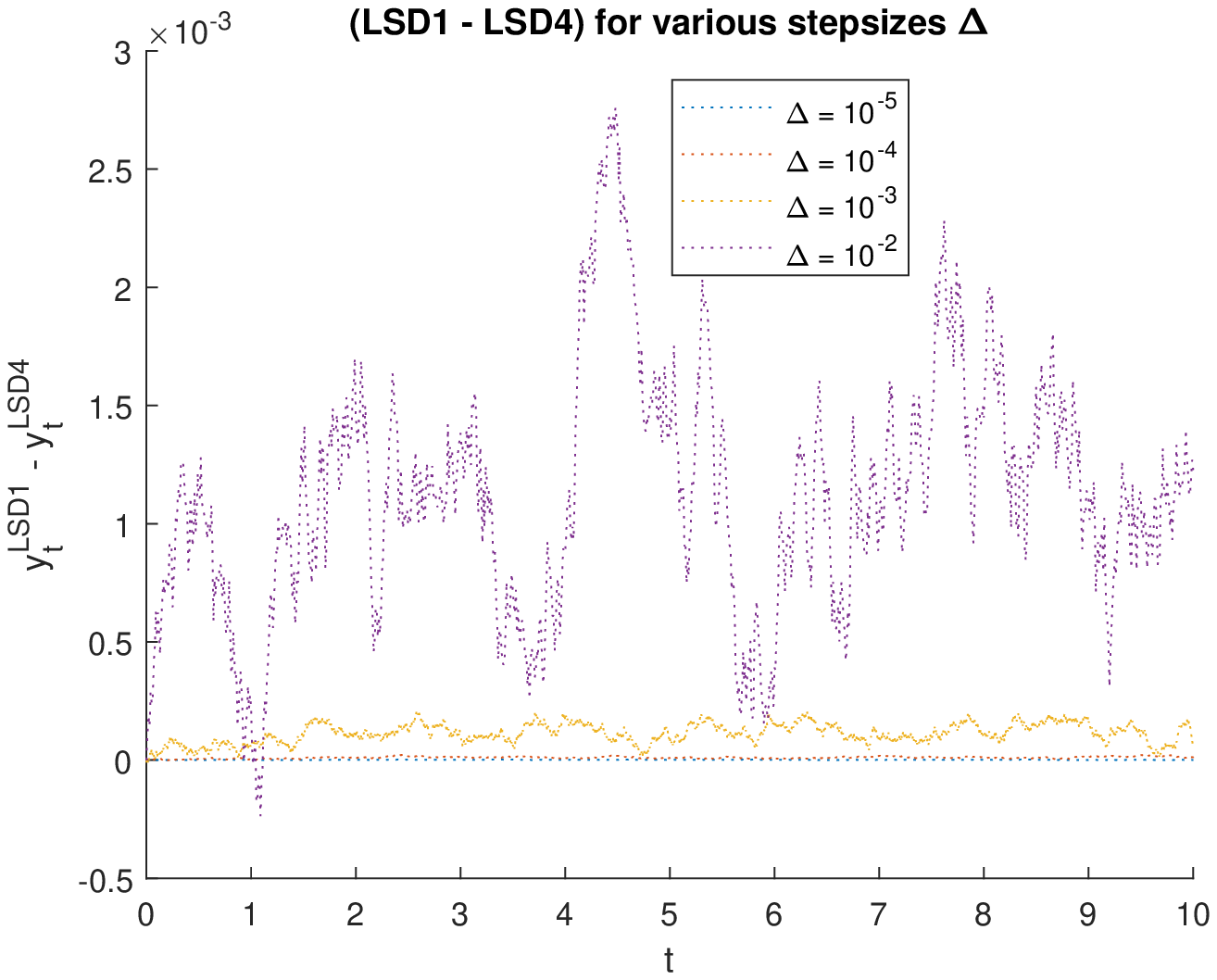}
	\caption{LSD1 - LSD4}
\end{subfigure}
\begin{subfigure}{.4\textwidth}
	\includegraphics[width=1\textwidth]{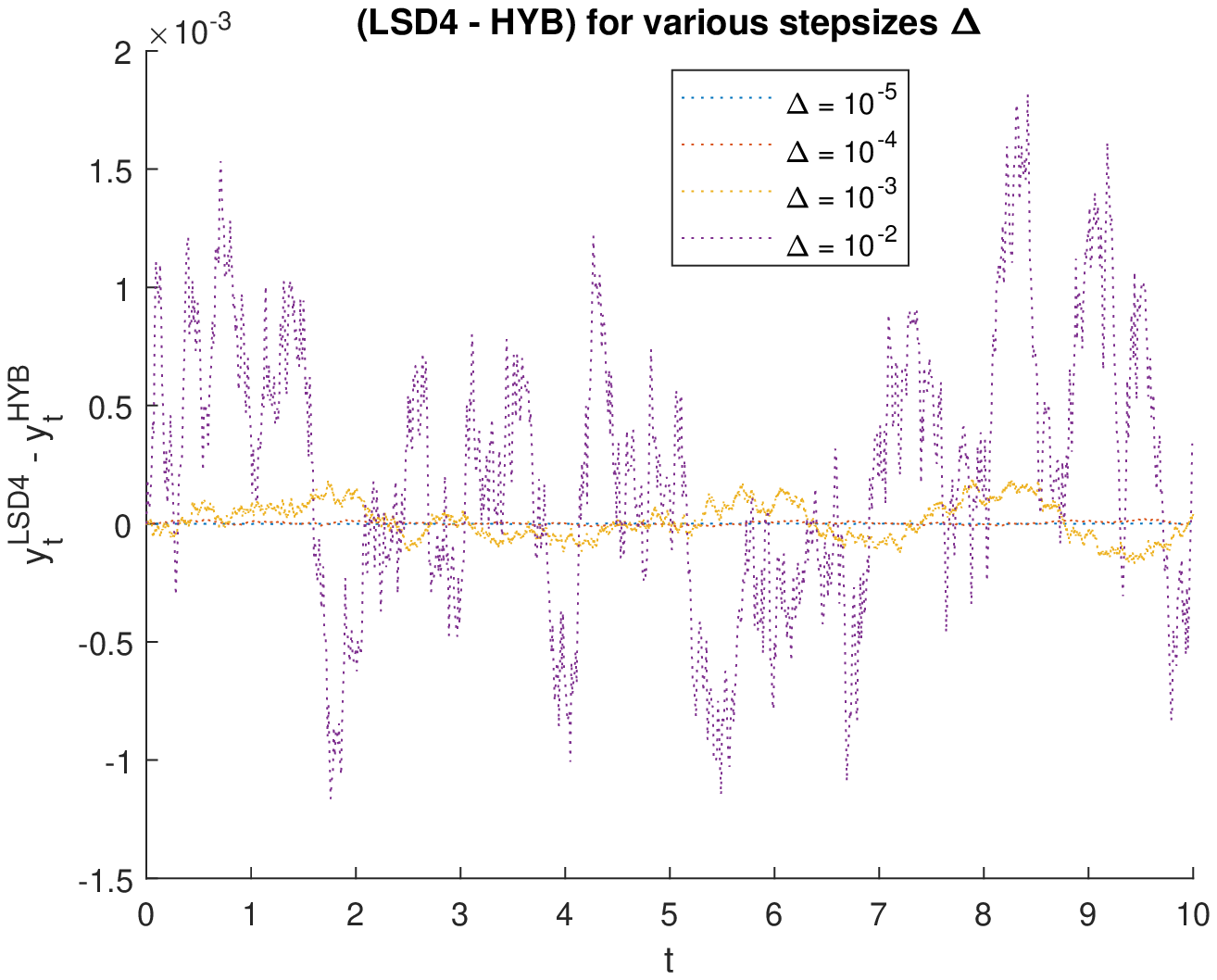}
	\caption{LSD4 - HYB}
\end{subfigure}
	\caption{Trajectories of the difference of numerical methods for the approximation of (\ref{LSD-eq:WFmodel}) with various $\D.$}\label{LSD-fig:LSDWFsminSDI}
\end{figure}

Figure \ref{LSD-fig:LSDWFSDI} shows paths for the proposed LSD and existing numerical methods for the Wright -Fisher model and Figure \ref{LSD-fig:LSDWFsminSDI}  shows a graphical estimation of the difference of the methods. 

\begin{figure}[h]
	\centering
	\begin{subfigure}{.47\textwidth}
		\includegraphics[width=1\textwidth]{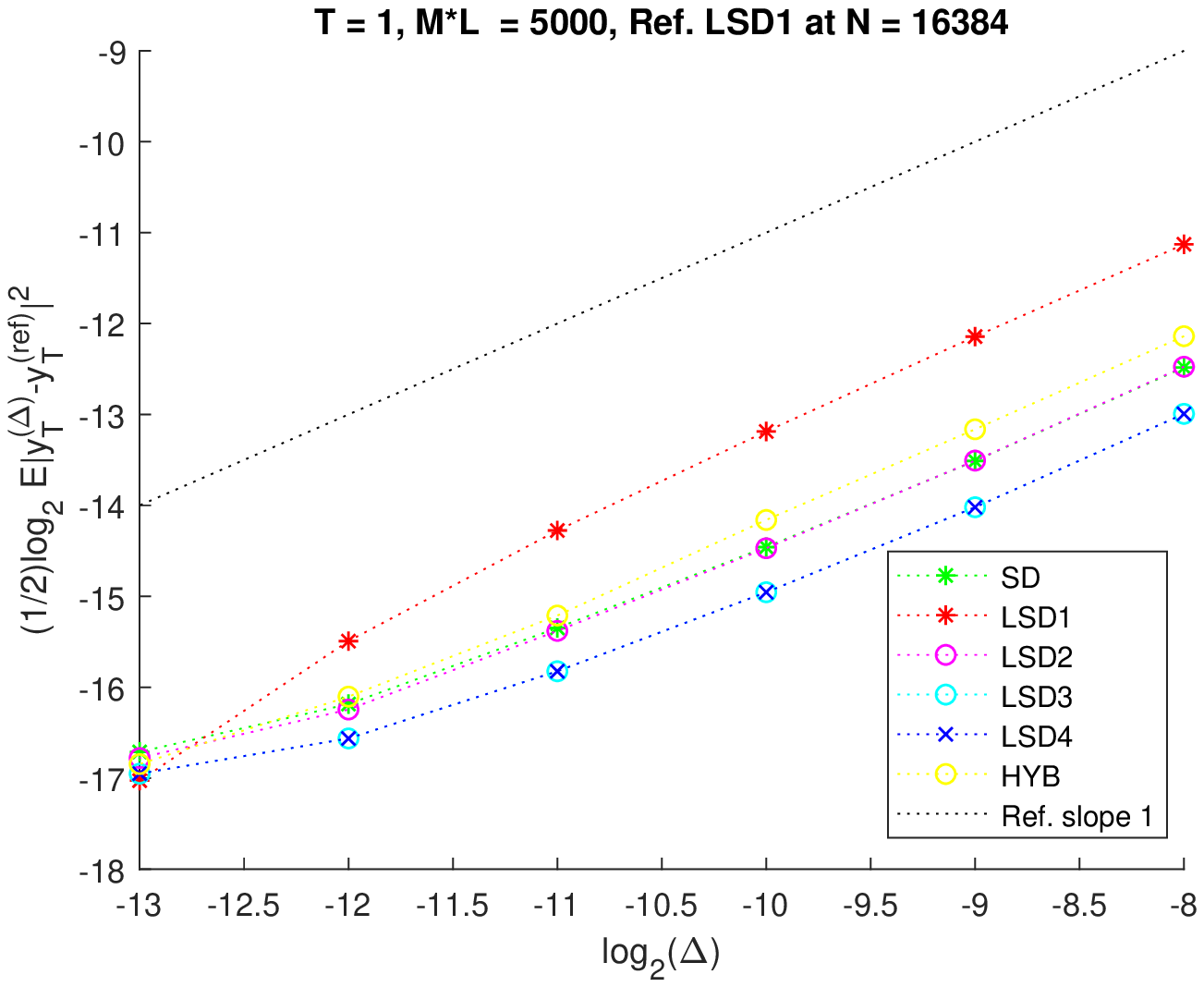}
		\caption{LSD1 reference solution.}
	\end{subfigure}
	\begin{subfigure}{.47\textwidth}
	\includegraphics[width=1\textwidth]{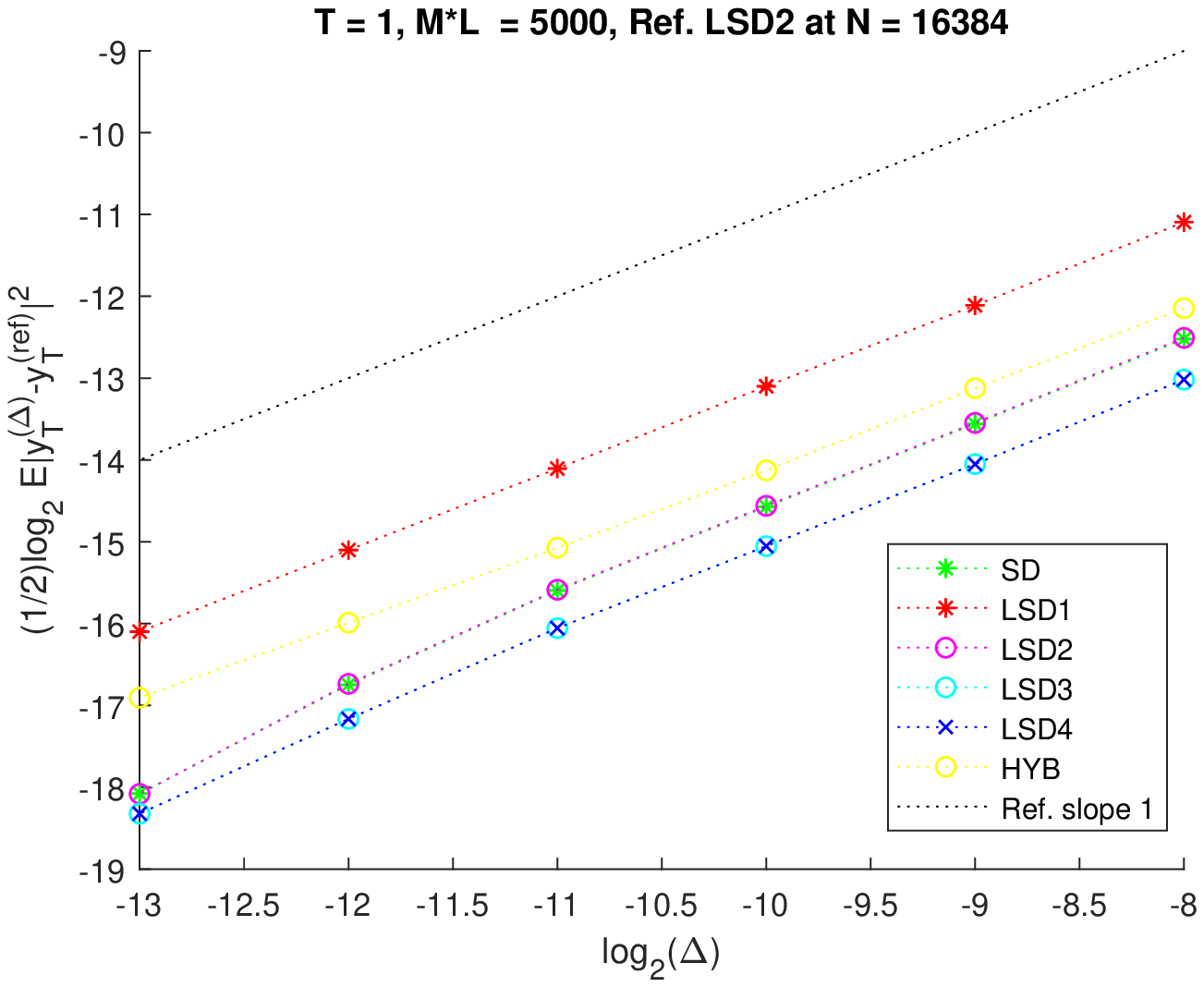}
	\caption{LSD2 reference solution.}
	\end{subfigure}
\begin{subfigure}{.47\textwidth}
		\includegraphics[width=1\textwidth]{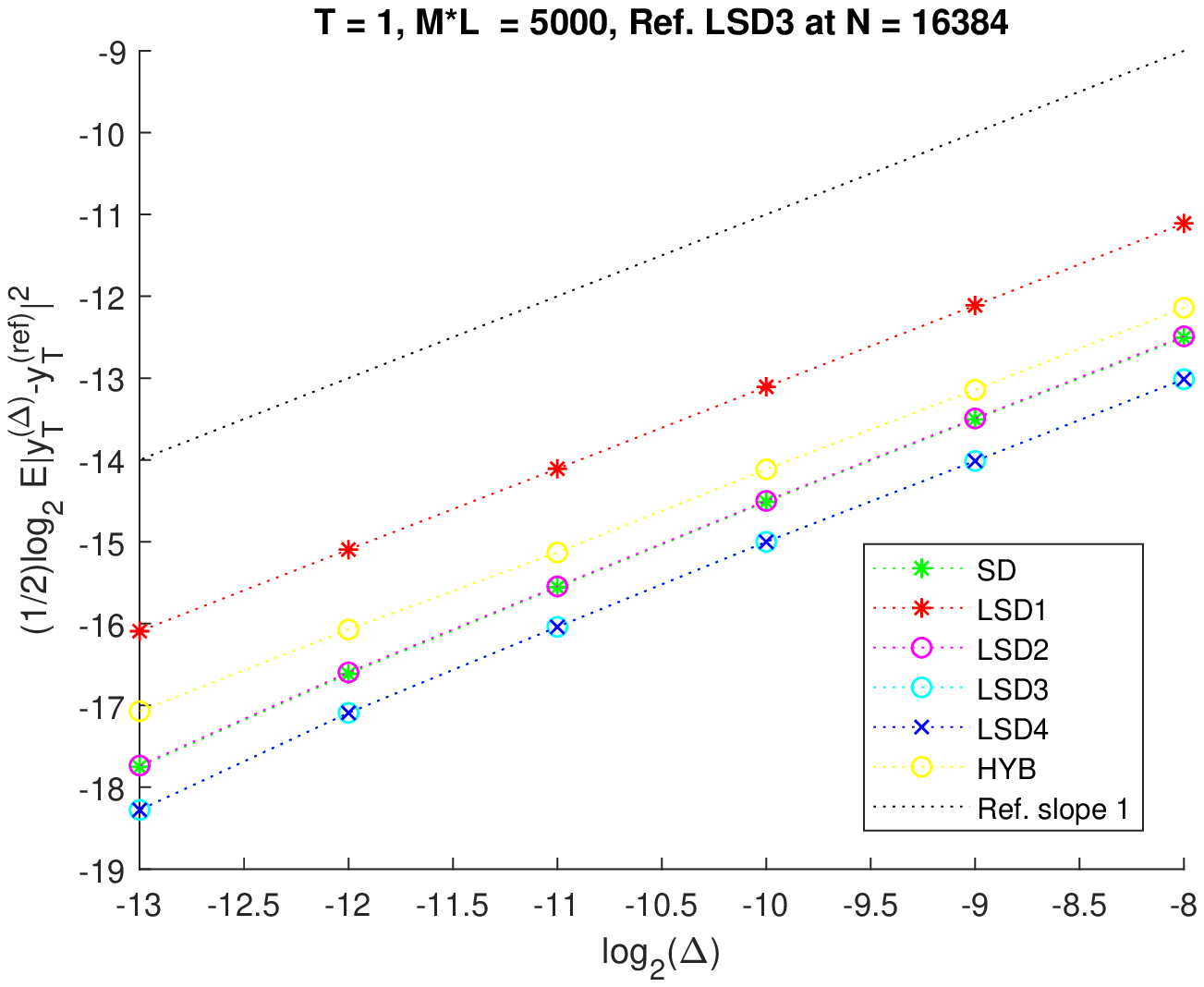}
		\caption{LSD3 reference solution.}
	\end{subfigure}
\begin{subfigure}{.47\textwidth}
	\includegraphics[width=1\textwidth]{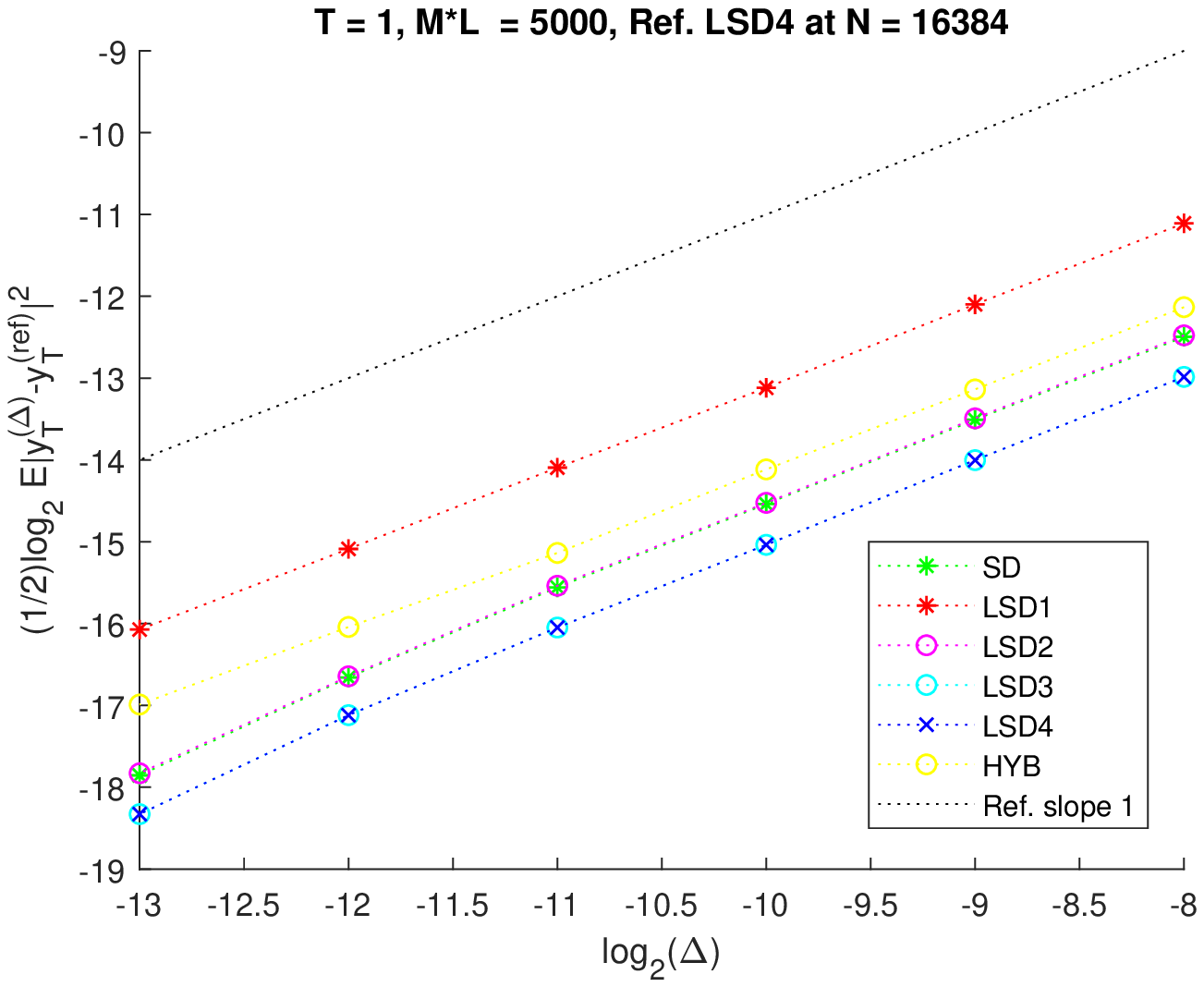}
	\caption{LSD4 reference solution.}
\end{subfigure}
\begin{subfigure}{.47\textwidth}
	\includegraphics[width=1\textwidth]{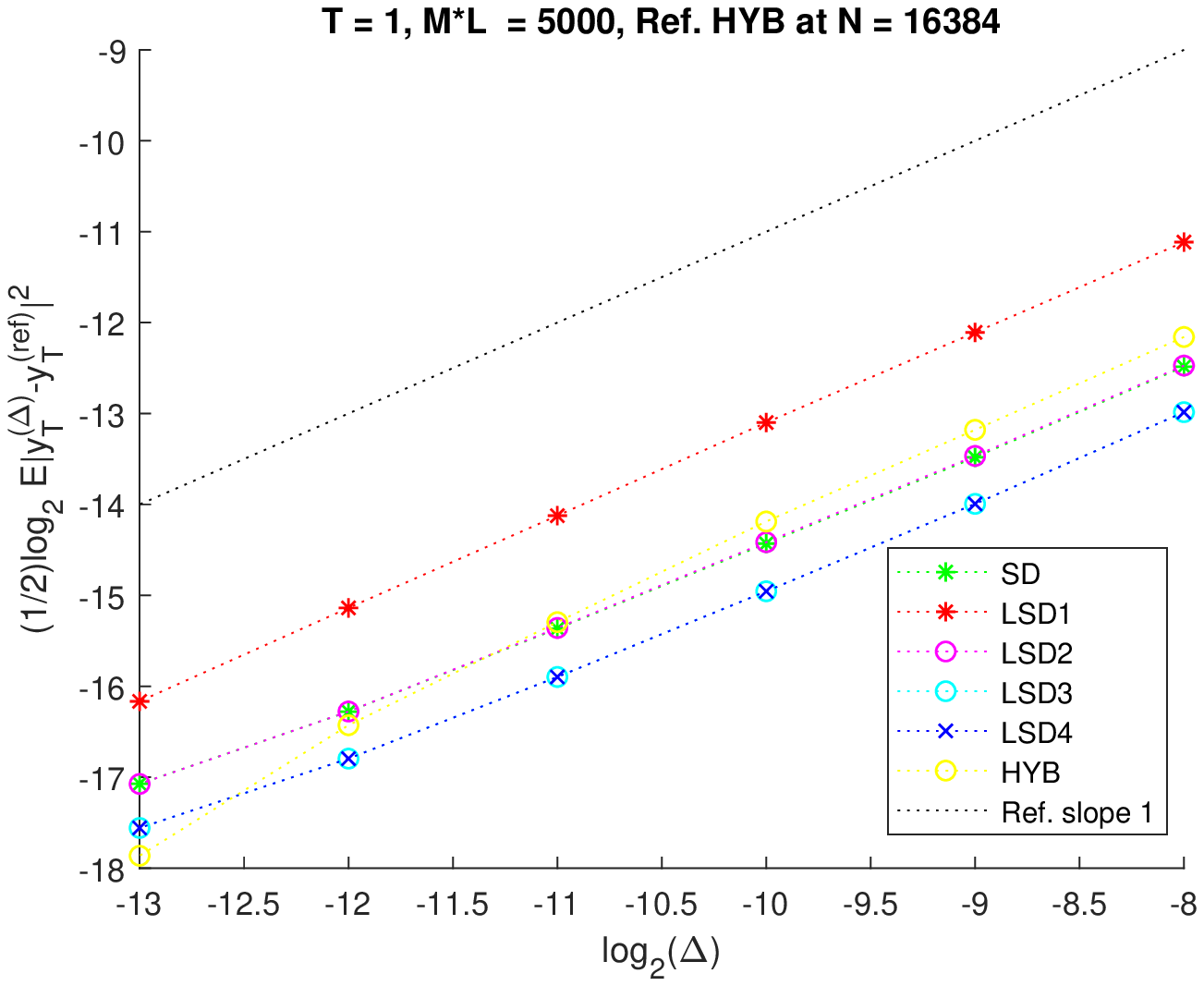}
	\caption{HYB reference solution.}
\end{subfigure}
\caption{Convergence of  LSD methods (\ref{LSD-eq:WFmodeloriginal}), (\ref{LSD-eq:WFmodeloriginal2}), (\ref{LSD-eq:WFmodeloriginal3}) and (\ref{LSD-eq:WFmodeloriginal4}) for the approximation of (\ref{LSD-eq:WFmodel}) with different reference solutions.}\label{LSD-fig:LSDWForder}
\end{figure}

We also examine numerically the order of strong convergence of the LSD method. The numerical results suggest that the LSD is mean-square convergent with order close to $1,$ see Figure \ref{LSD-fig:LSDWForder}. 

\section{Heston $3/2$-model}\label{LSD:ssec:Heston}
Let
\beqq  \label{LSD-eq:Heston}
x_t =x_0  + \int_0^t (k_1x_s - k_2(x_s)^2) ds + \int_0^t k_3(x_s)^{3/2} dW_s, \quad t\geq0,
\eeqq
where the coefficients $k_i, i=1,2,3$ are positive and $x_0>0.$  SDE (\ref{LSD-eq:Heston}) is known as the Heston $3/2$-model appearing in financial mathematics as a stochastic volatility process, see \cite{heston:1997}, and satisfies $x_t>0$ a.s. 
The Lamperti transformation of (\ref{LSD-eq:Heston}) is $z = \frac{2}{k_3}x^{-1/2}$ implying that, see Appendix \ref{LSD-ap:Lamperti_tranformation},  

\beqq  \label{LSD-eq:HestonLamperti}
z_t =z_0 + \int_0^t \left( (\frac{2k_2}{(k_3)^2} + 3)(z_s)^{-1} - \frac{k_1}{2}z_s\right)ds + \int_0^t  dW_s, \,\, t\geq0.
\eeqq

\subsection{Lamperti Semi-Discrete method $\wt{z}^1_{n}$ and $\wt{z}^2_{n}$ for Heston $3/2$-model}\label{LSD:subsec:Heston}

As in Section \ref{LSD:subsec:CIRz1z2} we consider the following two versions of the semi-discrete method for approximating  (\ref{LSD-eq:HestonLamperti}),

\beqq\label{LSD-eq:SDschemeHestonLT}
y_t = \D W_n + y_{t_n} - \frac{k_1}{2}y_{t_n}\D +  \int_{t_n}^t (\frac{2k_2}{(k_3)^2} + 3)(y_s)^{-1}ds, \,\, t\in (t_n, t_{n+1}],
\eeqq
with $y_{0}=z_0$ and 
\beqq\label{LSD-eq:SDschemeHestonLT2}
\hat{y}_t = \D W_n + \hat{y}_{t_n} + \int_{t_n}^t \left( (\frac{2k_2}{(k_3)^2} + 3)(\hat{y}_s)^{-1} - \frac{k_1}{2}\hat{y}_s\right)ds, \,\, t\in(t_n, t_{n+1}],
\eeqq
with $\hat{y}_{0}=z_0.$  (\ref{LSD-eq:SDschemeHestonLT}) and (\ref{LSD-eq:SDschemeHestonLT2}) are Bernoulli type equations with solutions satisfying, see Appendix \ref{LSD-ap:Bernoulli_sol},
\beqq\label{LSD-eq:SDschemeHestonLTsol}
(y_t)^2 = \left(\D W_n +  \left(1-\frac{k_1\D}{2}\right)y_{t_n}\right)^2 + (\frac{4k_2}{(k_3)^2} + 6)(t-t_n)
\eeqq
and
\beqq\label{LSD-eq:SDschemeHestonLTsol2}
(\hat{y}_t)^2 = (\D W_n + \hat{y}_{t_n})^2e^{-k_1(t-t_n)} + (\frac{4k_2}{(k_3)^2} + 6)\frac{1-e^{-k_1(t-t_n)}}{k_1},
\eeqq
respectively.
We propose the following versions of the semi-discrete method for the approximation of (\ref{LSD-eq:exampleSDELamperti}),  
\beqq\label{LSD-eq:SDschemeHestonLT_transf}
y_{t_{n+1}} = \sqrt{\left(\D W_n +  \left(1-\frac{k_1\D}{2}\right)y_{t_n}\right)^2 + (\frac{4k_2}{(k_3)^2} + 6)\D}
\eeqq
and 
\beqq\label{LSD-eq:SDschemeHestonLT_transf2}
\hat{y}_{t_{n+1}} = \sqrt{(\D W_n + \hat{y}_{t_n})^2e^{-k_1\D} + (\frac{4k_2}{(k_3)^2} + 6)\frac{1-e^{-k_1\D}}{k_1}},
\eeqq
which suggests the versions of the Lamperti semi-discrete method $(\wt{z}_n)_{n\in\bbN}$ for the approximation of (\ref{LSD-eq:exampleSDE})
\beqq\label{LSD-eq:SDschemeHestonLToriginal}
\wt{z}^1_{t_{n+1}} = \frac{4}{(k_3)^2}\left(\left(\D W_n +  \left(1-\frac{k_1\D}{2}\right)y_{t_n}\right)^2 + (\frac{4k_2}{(k_3)^2} + 6)\D\right)^{-1}
\eeqq

\beqq\label{LSD-eq:SDschemeHestonLToriginal2}
\wt{z}^2_{t_{n+1}} = \frac{4}{(k_3)^2}\left((\D W_n + \hat{y}_{t_n})^2e^{-k_1\D} + (\frac{4k_2}{(k_3)^2} + 6)\frac{1-e^{-k_1\D}}{k_1}\right)^{-1}
\eeqq

\subsection{Numerical experiment for Heston $3/2$-model}\label{LSD:subsec:Hestonnum}

For a minimal numerical experiment we present simulation paths for the numerical approximation of (\ref{LSD-eq:Heston}) with  $x_0=1$ and compare with the SD method proposed in \cite[Sec. 5]{halidias_stamatiou:2016}, which reads 
\beqq\label{LSD-eq:SDsol} 
\wt{y}_{t_{n+1}}=\wt{y}_{t_{n}}\exp\left\{\left(k_1 - k_2\wt{y}_{t_n} - \frac{(k_3)^2}{2}\wt{y}_{t_n}\right)\D  + k_3\sqrt{\wt{y}_{t_n}}\D W_n\right\}.
\eeqq
The semi-discrete method (\ref{LSD-eq:SDsol}) is strongly converging and positivity preserving, see \cite[Sec. 5]{halidias_stamatiou:2016}.    

Moreover, we compare with the implicit method proposed in \cite{neuenkirch_szpruch:2014}.
Set
 $$G(x) = \left(1 + \frac{k_1}{2}\D\right) x - \left(\frac{k_2}{2} + \frac{3(k_3)^2}{8}\right)\D x^{-1}$$ and compute
$$
y_{n+1} = G^{-1}\left(y_{n} - \frac{k_3}{2}\D W_n\right)
$$
and then transform back to get the following scheme
\beqq\label{LSD-eq:ImplicitschemeHeston}
y_{n+1}^{Impl}= (y_{n+1})^{-2}.
\eeqq

\begin{figure}[ht]
	\centering
	\begin{subfigure}{.47\textwidth}
		\includegraphics[width=1\textwidth]{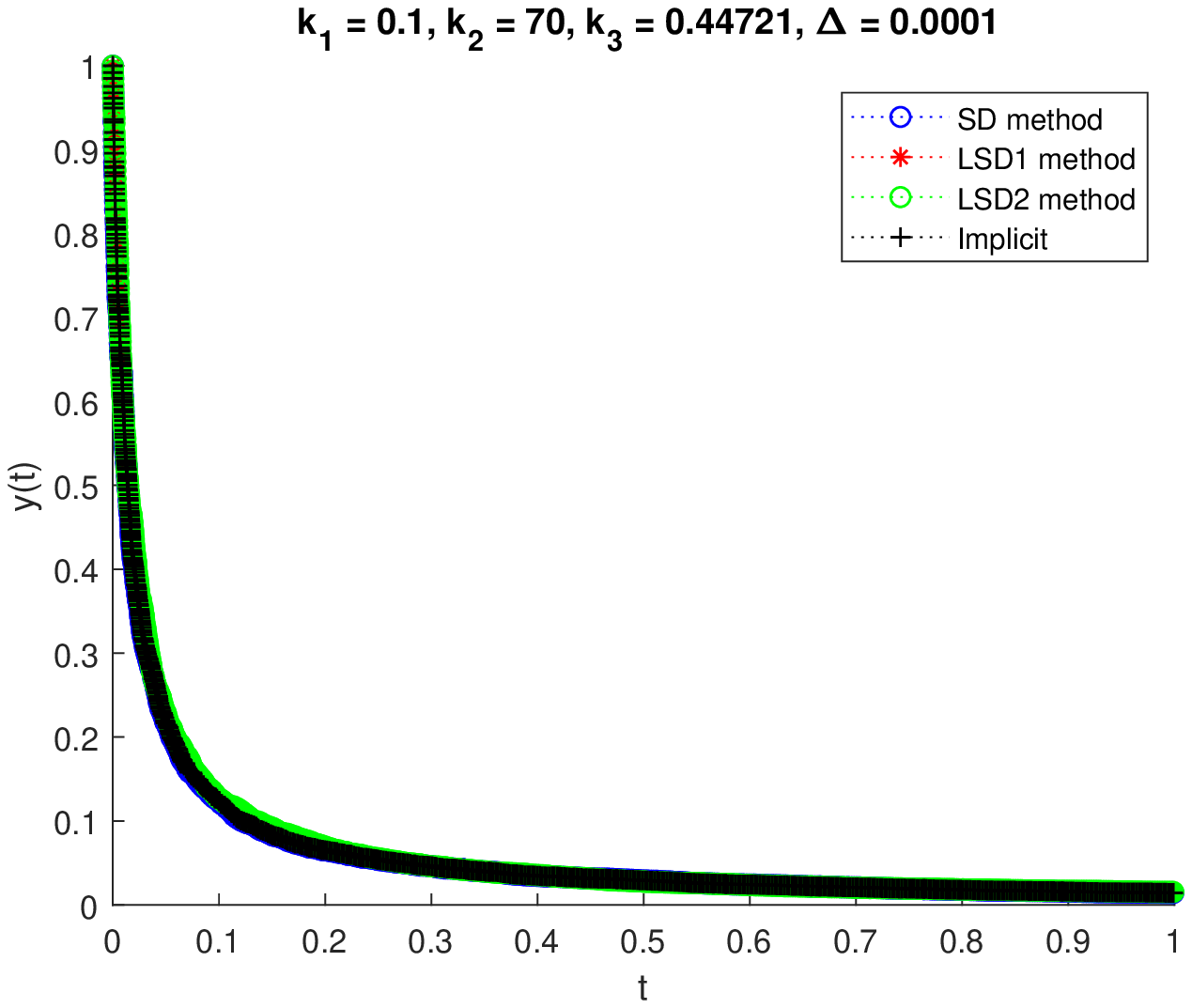}
		\caption{Trajectories  of  (\ref{LSD-eq:SDschemeHestonLToriginal})-(\ref{LSD-eq:ImplicitschemeHeston}).}
	\end{subfigure}
	\begin{subfigure}{.47\textwidth}
		\includegraphics[width=1\textwidth]{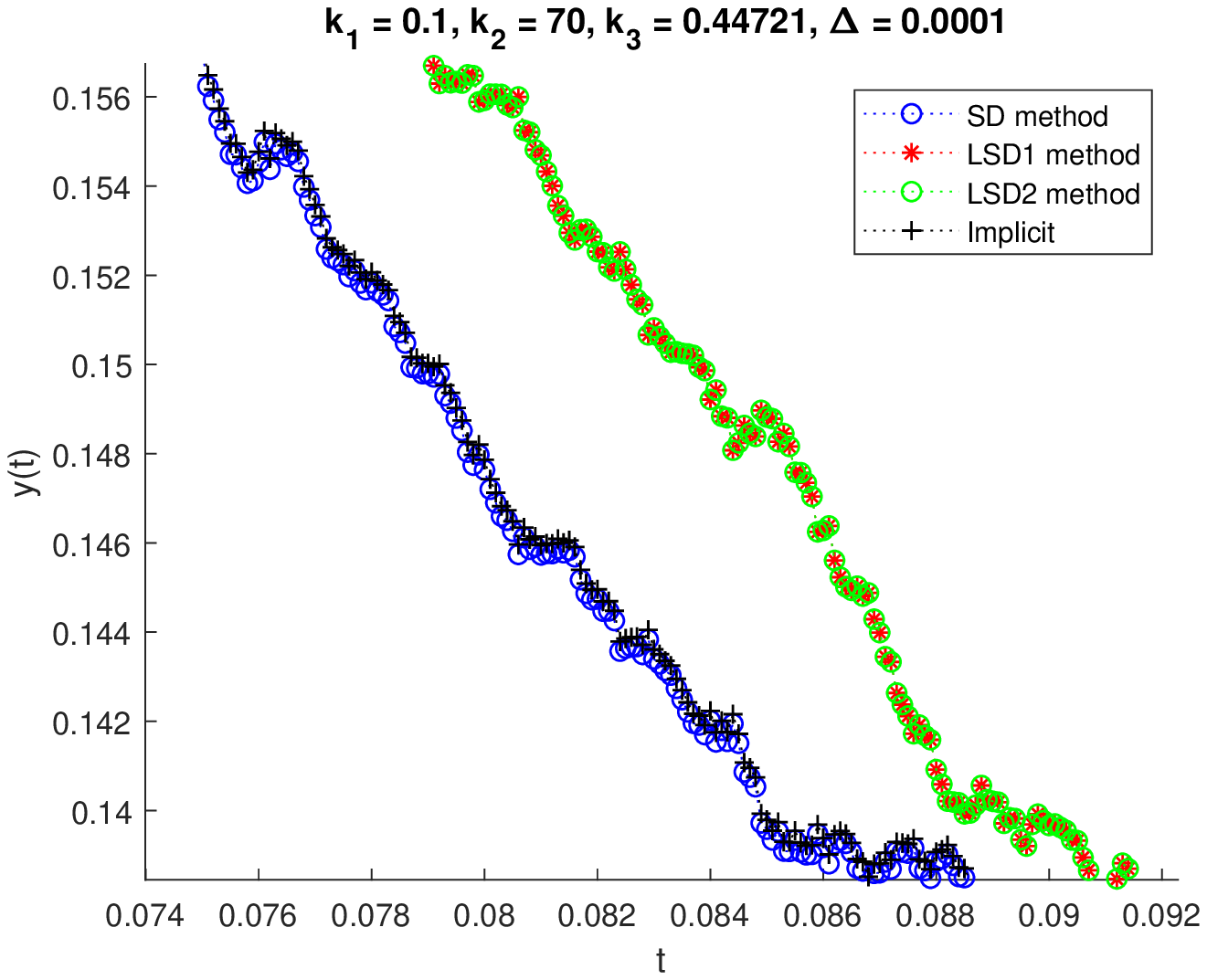}
		\caption{Zoom of Figure \ref{LSD-fig:LSDsSDHT}(A).}
	\end{subfigure}
	\caption{Trajectories  of  (\ref{LSD-eq:SDschemeHestonLToriginal}), (\ref{LSD-eq:SDschemeHestonLToriginal2}), (\ref{LSD-eq:SDsol}) and (\ref{LSD-eq:ImplicitschemeHeston}) for the approximation of (\ref{LSD-eq:Heston}) with $\D=10^{-4}$.}\label{LSD-fig:LSDsSDHT}
\end{figure}

\begin{figure}[ht]
	\centering
	\begin{subfigure}{.47\textwidth}
		\includegraphics[width=1\textwidth]{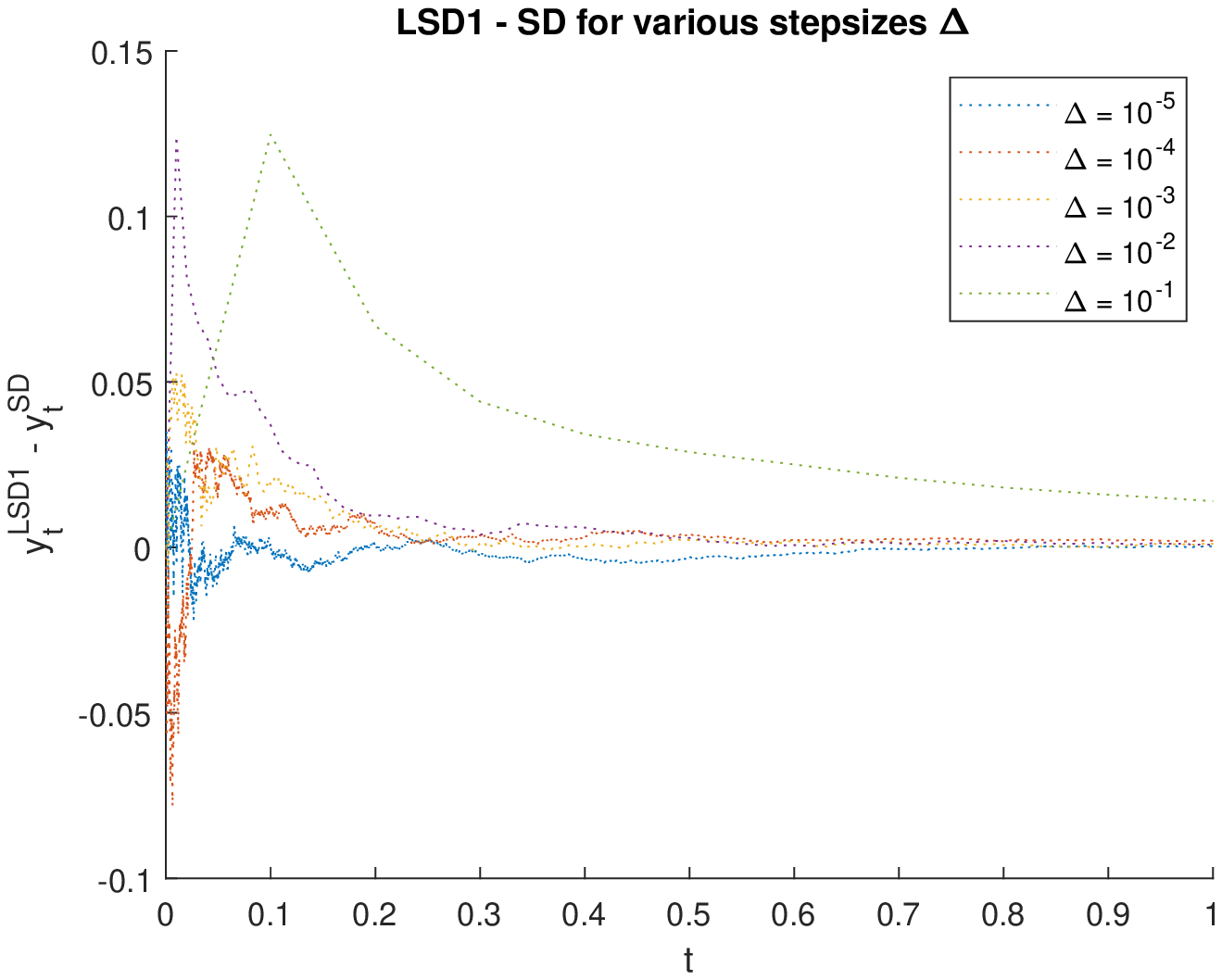}
		\caption{Difference (\ref{LSD-eq:SDschemeHestonLToriginal}) - (\ref{LSD-eq:SDsol}) for various step sizes.}
	\end{subfigure}
	\begin{subfigure}{.47\textwidth}
		\centering
		\includegraphics[width=1\textwidth]{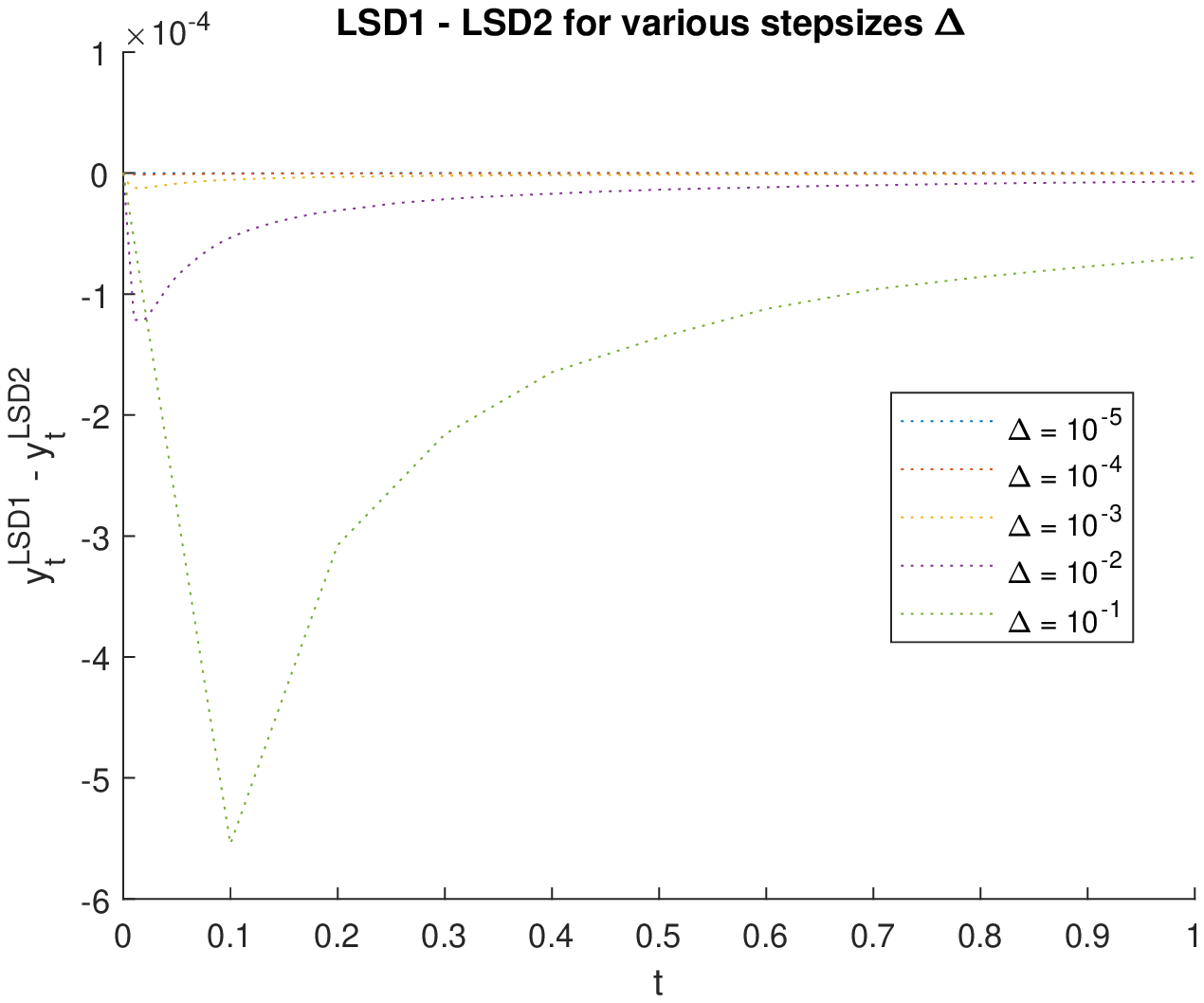}
		\caption{Difference (\ref{LSD-eq:SDschemeHestonLToriginal}) - (\ref{LSD-eq:SDschemeHestonLToriginal2}) for various step sizes.}
	\end{subfigure}
	\begin{subfigure}{.55\textwidth}
	\centering
	\includegraphics[width=1\textwidth]{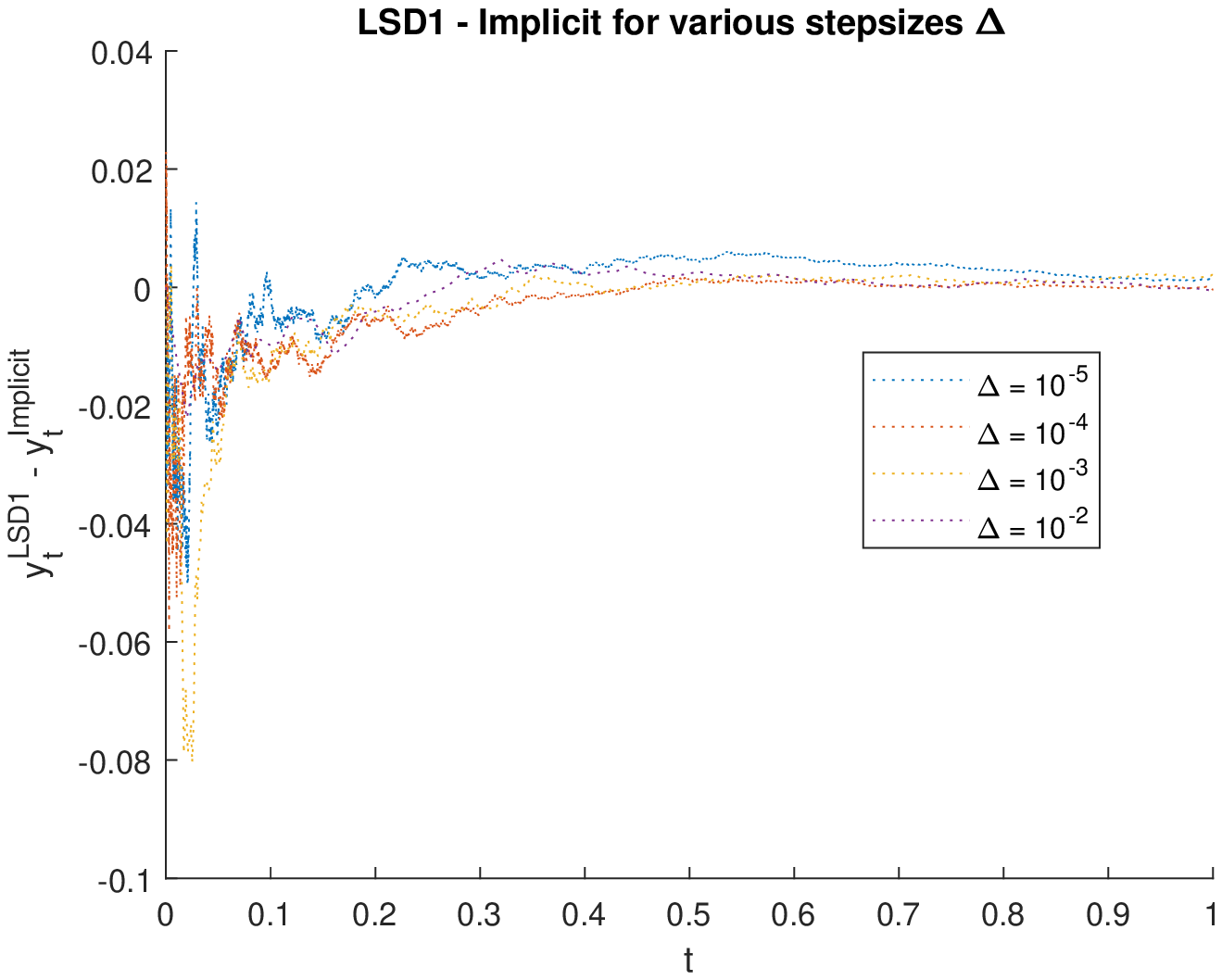}
	\caption{Difference (\ref{LSD-eq:SDschemeHestonLToriginal}) - (\ref{LSD-eq:ImplicitschemeHeston}) for various step sizes.}
\end{subfigure}
	\caption{Trajectories of the differences of the semi-discrete methods (\ref{LSD-eq:SDschemeHestonLToriginal}), (\ref{LSD-eq:SDschemeHestonLToriginal2}), (\ref{LSD-eq:SDsol}) and the implicit scheme (\ref{LSD-eq:ImplicitschemeHeston}) for the approximation of (\ref{LSD-eq:Heston}) with various step-sizes.}\label{LSD-fig:LSDsminSDHT}
\end{figure}

We use the set of parameters from \cite[Sec.5]{halidias_stamatiou:2016}; we take $k_1 = 0.1, k_2 =70, k_3 = \sqrt{0.2}, T = 1$ with various step-sizes. We compare the proposed two versions of LSD scheme (\ref{LSD-eq:SDschemeHestonLToriginal}) and (\ref{LSD-eq:SDschemeHestonLToriginal2}) with the SD scheme (\ref{LSD-eq:SDsol}) and the implicit method (\ref{LSD-eq:ImplicitschemeHeston}). Figure \ref{LSD-fig:LSDsSDHT} shows that the pair of LSD1 and LSD2 and the pair of SD with the implicit method are almost identical for a step-size $\D = 10^{-4}.$ Moreover, the two pairs are getting very close. We give a presentation of the difference of the various approximations in Figure \ref{LSD-fig:LSDsminSDHT}.

\begin{figure}[ht]
	\centering
	\begin{subfigure}{.47\textwidth}
		\centering
		\includegraphics[width=1\textwidth]{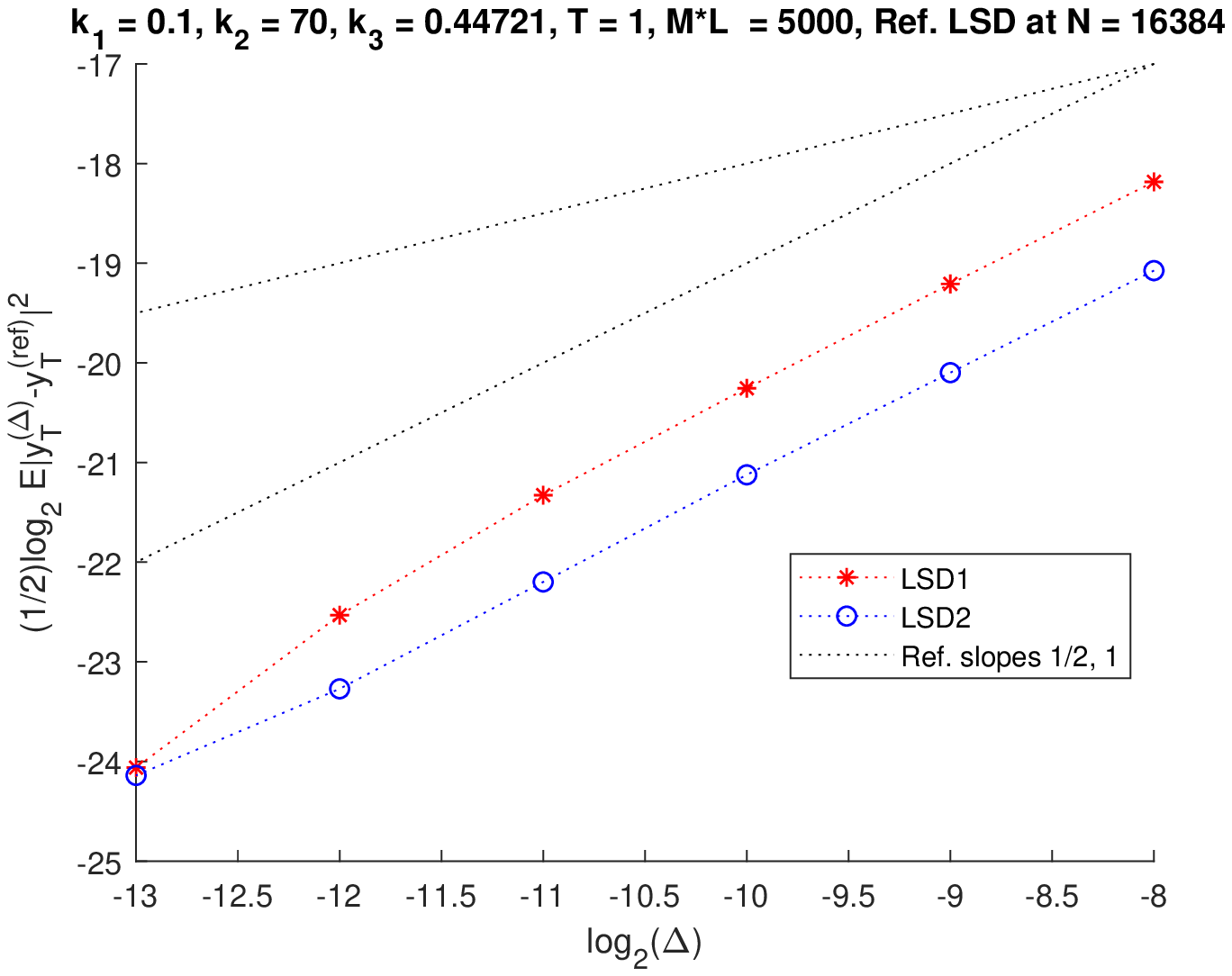}
		\caption{LSD1 as reference solution.}
	\end{subfigure}
	\begin{subfigure}{.47\textwidth}
		\centering
		\includegraphics[width=1\textwidth]{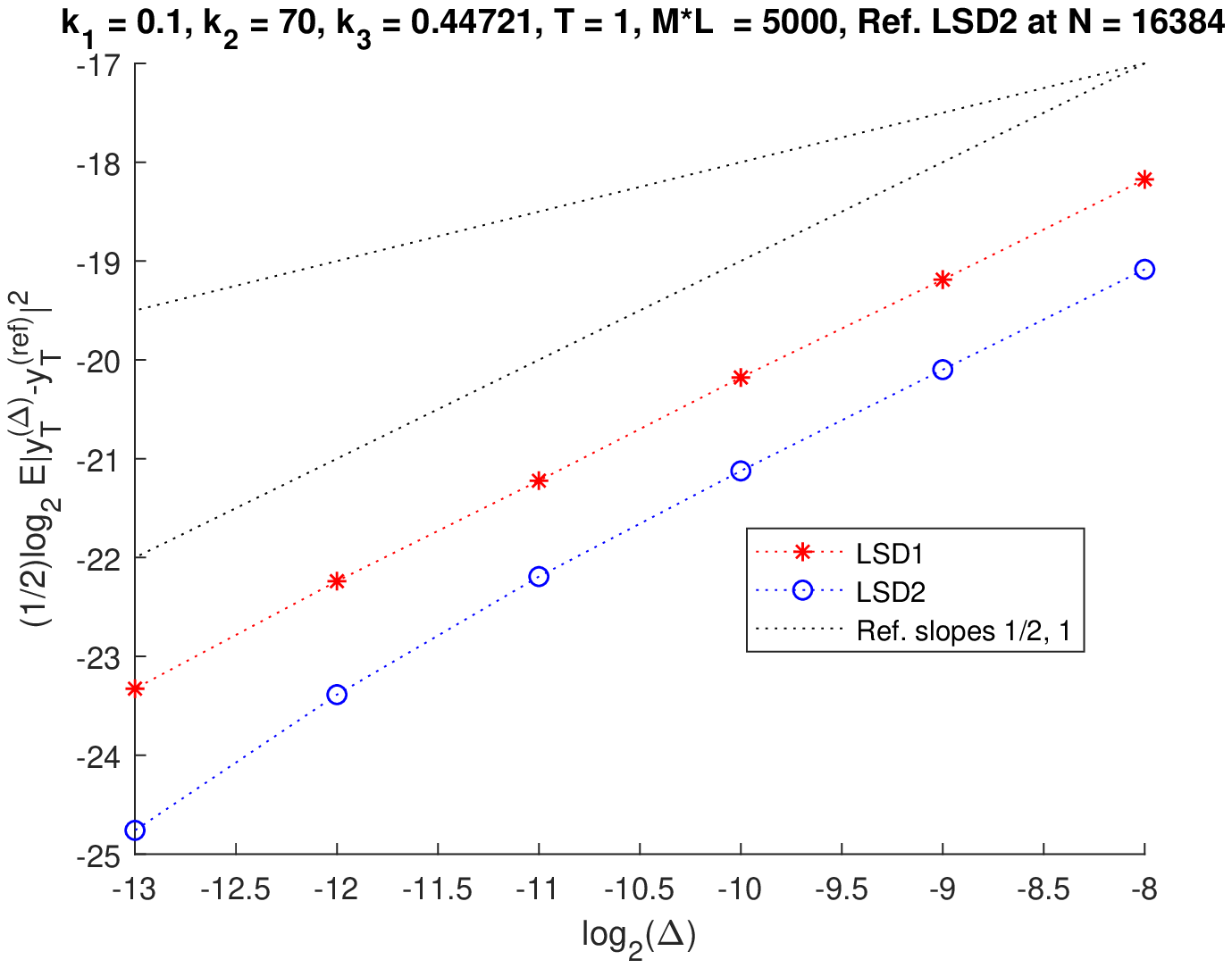}
		\caption{LSD2 as reference solution.}
	\end{subfigure}
	\caption{Convergence of  LSD1 and LSD2 method (\ref{LSD-eq:SDschemeHestonLToriginal}) and (\ref{LSD-eq:SDschemeHestonLToriginal}) for the approximation of (\ref{LSD-eq:Heston}) with different reference solution.}\label{LSD-fig:LSDorderHT}
\end{figure}

Finally, we examine numerically the order of strong convergence of the LSD method. The numerical results suggest that the LSD1 scheme as well as LSD2 converge in the mean-square sense with order close to $1,$ see Figure \ref{LSD-fig:LSDorderHT}. 

\section{A\"it-Sahalia model}\label{LSD:ssec:Ait}
Let
\beqq  \label{LSD-eq:Ait}
x_t =x_0  + \int_0^t (k_{-1}(x_s)^{-1} -k_0 + k_1x_s - k_2 (x_s)^r) ds + \int_0^t k_3(x_s)^{\rho} dW_s, \quad t\geq0.
\eeqq
where $k_{-1}, k_0, k_1, k_2$ and $k_3$ are positive constants with $r>1$ and $\rho>1.$ SDE (\ref{LSD:ssec:Ait}) is an A\"it-Sahalia model with superlinear coefficients and the property $x_t>0$ a.s. The Lamperti transformation of (\ref{LSD-eq:Ait}) is $z = x^{1-\rho}$ with dynamics, see Appendix \ref{LSD-ap:Lamperti_tranformation},  
\beam\nonumber  
z_t &=&z_0 + \int_0^t \Big(k_{-1}(1-\rho)(z_s)^{\frac{\rho+1}{\rho-1}} -k_{0}(1-\rho)(z_s)^{\frac{\rho}{\rho-1}} + k_{1}(1-\rho)z_s\\
\label{LSD-eq:AitLamperti}&& -k_{2}(1-\rho)(z_s)^{\frac{\rho - r}{\rho-1}} - \frac{\rho(1-\rho)(k_3)^2}{2}(z_s)^{-1} \Big)ds + \int_0^t k_3(1-\rho)dW_s.
\eeam 

Set $K_{i} = k_{i}(\rho-1), i=-1,\dots,3$ and $K_4 = \frac{\rho(\rho-1)(k_3)^2}{2}.$
We examine the new version $(y_t),$ of the semi-discrete method for approximating  (\ref{LSD-eq:AitLamperti}) see Section \ref{LSD:subsec:Aitz1}, where in each subinterval $(t_n, t_{n+1}]$ we solve an algebraic equation, producing a positive numerical scheme.

\subsection{Lamperti Semi-Discrete methods $\wt{z}^1_{n}, \wt{z}^2_{n}$ for A\"it-Sahalia}\label{LSD:subsec:Aitz1}

Rewrite (\ref{LSD-eq:AitLamperti}) as

\beam\nonumber  
z_t &=&z_0 + \int_0^t \left(-K_{-1}(z_s)^{\frac{\rho+1}{\rho-1}} +K_{0}(z_s)^{\frac{\rho}{\rho-1}} - K_{1}z_s  + K_2(z_s)^{\frac{\rho - r}{\rho-1}} + K_4(z_s)^{-1} \right)ds\\
\label{LSD-eq:AitLampertiK}&& \qquad \quad -  \int_0^t K_3dW_s.
\eeam 
Let $t\in(t_n, t_{n+1}]$ and 
\beam\nonumber
y_t &=&  -K_3(W_t - W_{t_n}) + y_{t_n} - K_{-1}(y_{t_n})^{\frac{2}{\rho-1}}y_t\D + K_{0}(y_{t_n})^{\frac{\rho}{\rho-1}}\D - K_{1}y_t\D\\
\label{LSD-eq:SDschemeAitLT}&&  + K_2(y_{t_n})^{\frac{2\rho - r-1}{\rho-1}}(y_t)^{-1}\D + K_4(y_t)^{-1}\D
\eeam
with $y_{0}=z_0$ and
\beam\nonumber
\hat{y}_t &=&  -K_3(W_t - W_{t_n}) + \hat{y}_{t_n} - K_{-1}(\hat{y}_{t_n})^{\frac{\rho+1}{\rho-1}}\D + K_{0}(\hat{y}_{t_n})^{\frac{\rho}{\rho-1}}\D - K_{1}\hat{y}_t\D\\
\label{LSD-eq:SDschemeAitLT2}&&  + K_2(\hat{y}_{t_n})^{\frac{2\rho - r-1}{\rho-1}}(\hat{y}_t)^{-1}\D + K_4(\hat{y}_t)^{-1}\D
\eeam
with $\hat{y}_{0}=z_0.$  The solutions of (\ref{LSD-eq:SDschemeAitLT}) and  (\ref{LSD-eq:SDschemeAitLT2}) are such that 
\beam\nonumber
&&(1+K_{-1}(y_{t_n})^{\frac{2}{\rho-1}}\D + K_1\D)(y_t)^{2} - \left(-K_3(W_t - W_{t_n}) + y_{t_n} + K_{0}(y_{t_n})^{\frac{\rho}{\rho-1}}\D\right)y_t\\
\label{LSD-eq:SDschemeAitLTsol}&&\qquad\qquad\quad -\left(K_2(y_{t_n})^{\frac{2\rho - r-1}{\rho-1}} + K_4 \right)\D = 0
\eeam
and
\beam\nonumber
&&(1+  K_1\D)(\hat{y}_t)^{2} - \left(-K_3(W_t - W_{t_n}) + \hat{y}_{t_n} + K_{0}(\hat{y}_{t_n})^{\frac{\rho}{\rho-1}}\D - K_{-1}(\hat{y}_{t_n})^{\frac{\rho+1}{\rho-1}}\D\right)\hat{y}_t\\
\label{LSD-eq:SDschemeAitLTsol2}&&\qquad\qquad\quad -\left(K_2(\hat{y}_{t_n})^{\frac{2\rho - r-1}{\rho-1}} + K_4 \right)\D = 0
\eeam
respectively.

We propose the following versions of the semi-discrete method for the approximation of (\ref{LSD-eq:AitLampertiK}),  
\beqq\label{LSD-eq:SDschemeAitLT_transf}
y_{t_{n+1}} = \frac{\phi_\D(y_{t_n},\D W_n) + \sqrt{\phi^2_\D(y_{t_n},\D W_n) + 4C_1(y_{t_n})C_2(y_{t_n})}}{2C_1(y_{t_n})},
\eeqq
with $\phi_\D(x,y) = -K_3y + x  + K_0x^{\frac{\rho}{\rho-1}}\D, C_1(x) = (1+K_{-1}x^{\frac{2}{\rho-1}}\D + K_1\D)$ and $C_2(x) = \left(K_2x^{\frac{2\rho - r-1}{\rho-1}} + K_4 \right)\D$ 
and 
\beqq\label{LSD-eq:SDschemeAitLT_transf2}
\hat{y}_{t_{n+1}} = \frac{\hat{\phi}_\D(\hat{y}_{t_n},\D W_n) + \sqrt{\hat{\phi}^2_\D(\hat{y}_{t_n},\D W_n) + 4(1 + K_1\D)C_2(\hat{y}_{t_n})}}{2(1 + K_1\D)},
\eeqq
with $\hat{\phi}_\D(x,y) = \phi_\D(x,y)  - K_{-1}x^{\frac{\rho+1}{\rho-1}}\D,$
which suggests the versions of the Lamperti semi-discrete method $(\wt{z}^1_n)_{n\in\bbN}$ and $(\wt{z}^2_n)_{n\in\bbN}$ for the approximation of (\ref{LSD-eq:Ait}) with $\wt{z}^1_{n} = (y_{n})^{1/(1-\rho)}, \wt{z}^2_{n} = (\hat{y}_{n})^{1/(1-\rho)}$
or
\beqq \label{LSD-eq:SDschemeAitLToriginal}
\wt{z}^1_{t_{n+1}} = \left| \frac{\phi_\D(y_{t_n},\D W_n) + \sqrt{\phi^2_\D(y_{t_n},\D W_n) + 4C_1(y_{t_n})C_2(y_{t_n})}}{2C_1(y_{t_n})} \right|^{\frac{1}{1-\rho}}.
\eeqq
\beqq \label{LSD-eq:SDschemeAitLToriginal2}
\wt{z}^2_{t_{n+1}} = \left| \frac{\hat{\phi}_\D(\hat{y}_{t_n},\D W_n) + \sqrt{\hat{\phi}^2_\D(\hat{y}_{t_n},\D W_n) + 4(1+K_1\D)C_2(\hat{y}_{t_n})}}{2(1+K_1\D)} \right|^{\frac{1}{1-\rho}}.
\eeqq
\subsection{Numerical experiment for A\"it-Sahalia}\label{LSD:subsec:Aitnum}

For a minimal numerical experiment we present simulation paths for the numerical approximation of (\ref{LSD-eq:Ait}) with  $x_0 = 4$ and compare with the implicit method proposed in \cite{neuenkirch_szpruch:2014}.
Set
$$G(x) = x + \left(1 + K_{-1}x^{\frac{\rho+1}{\rho-1}}  - K_{0}x^{\frac{\rho}{\rho-1}} + K_{1}x  - K_2x^{\frac{\rho - r}{\rho-1}} + K_4x^{-1}\right)\D 
$$
and compute
$$
y_{n+1} = G^{-1}\left(y_{n} -K_3\D W_n\right)
$$
and then transform back to get the following scheme
\beqq\label{LSD-eq:ImplicitschemeAIT}
y_{n+1}^{Impl}= (y_{n+1})^{1/(1-\rho)}.
\eeqq

\begin{figure}[ht]
	\centering
	\begin{subfigure}{.47\textwidth}
		\includegraphics[width=1\textwidth]{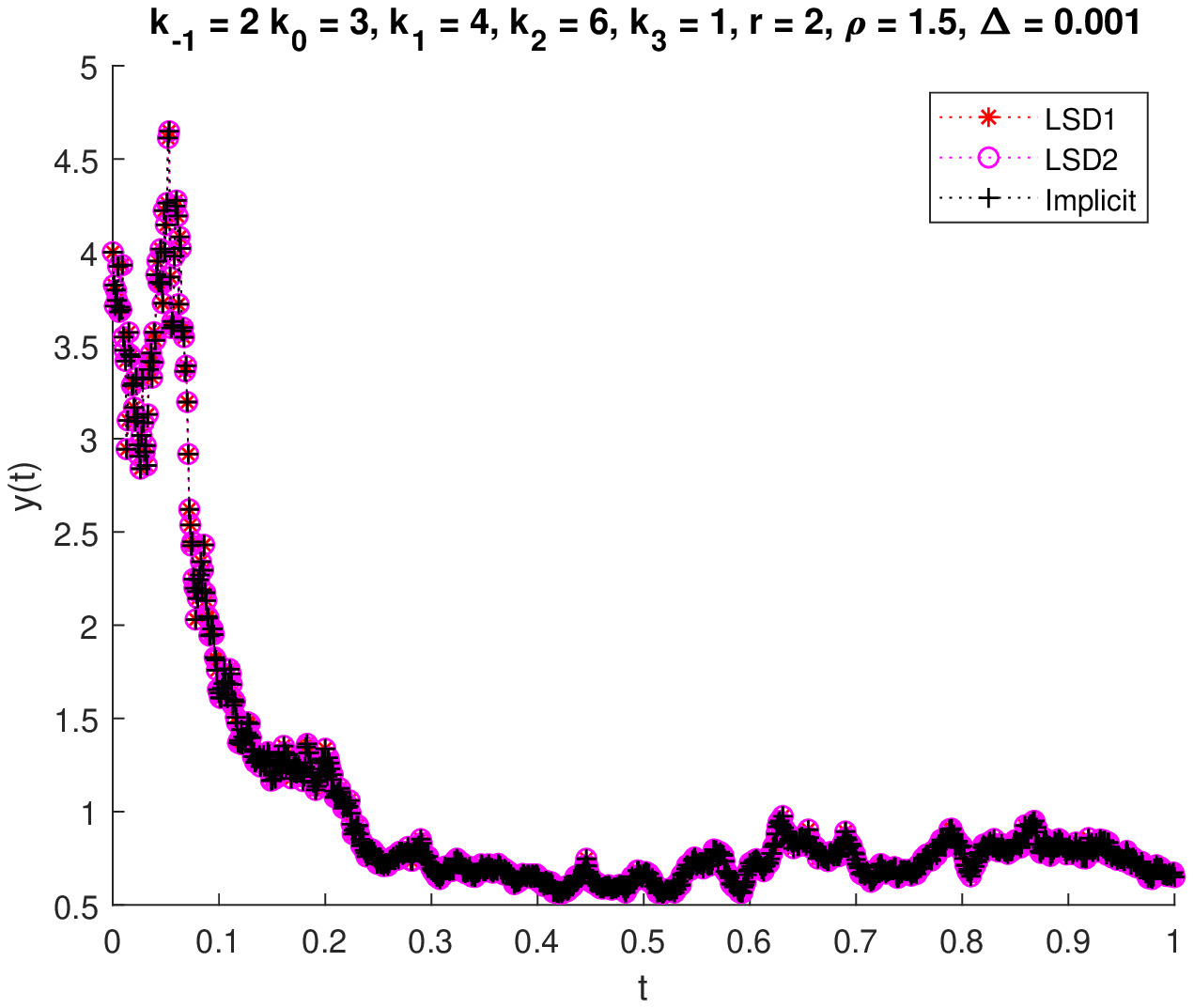}
		\caption{A path with $\D = 10^{-3}$.}
	\end{subfigure}
	\begin{subfigure}{.47\textwidth}
		\includegraphics[width=1\textwidth]{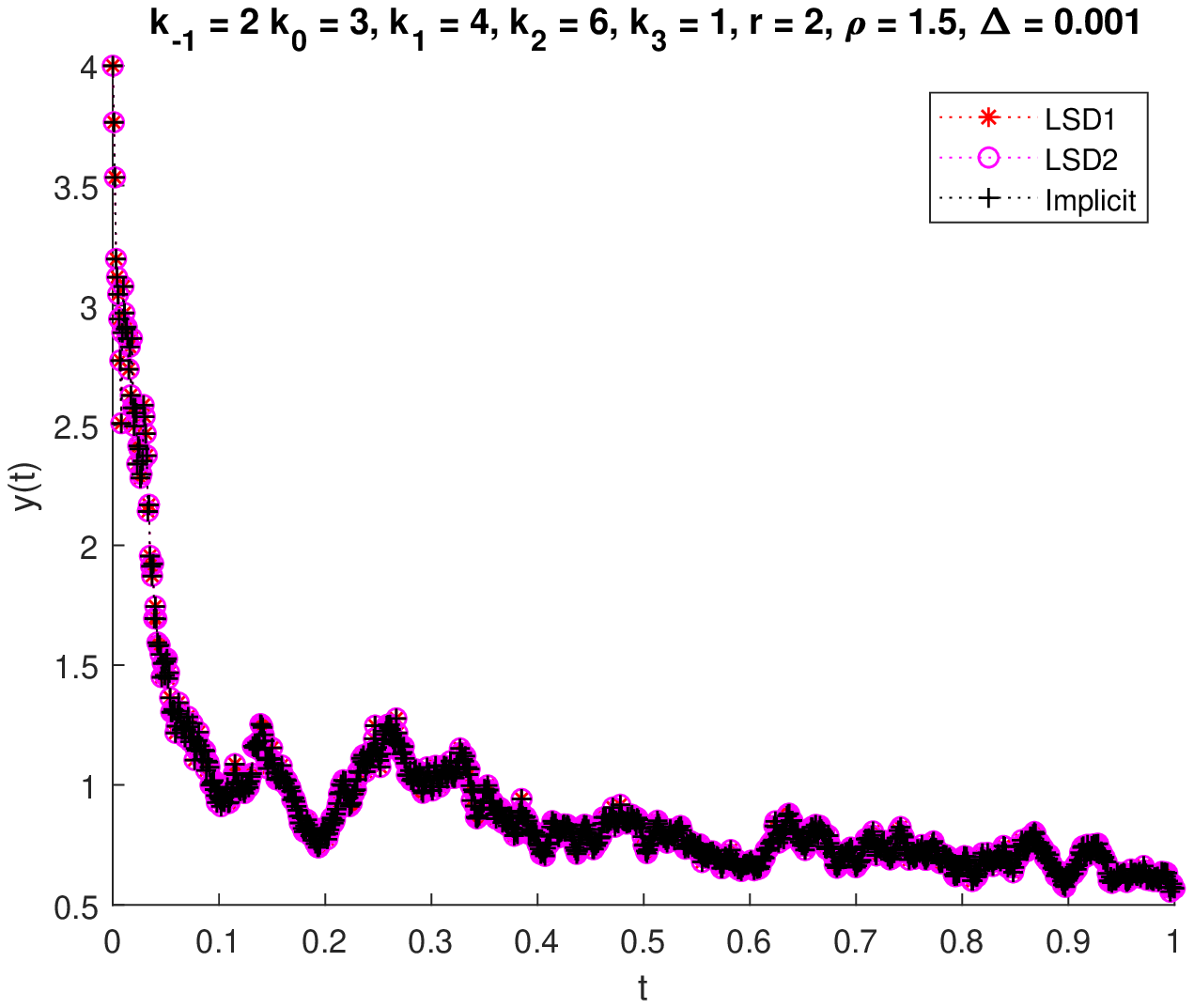}
		\caption{A path with $\D = 10^{-3}$.}
	\end{subfigure}
	\begin{subfigure}{.47\textwidth}
	\includegraphics[width=1\textwidth]{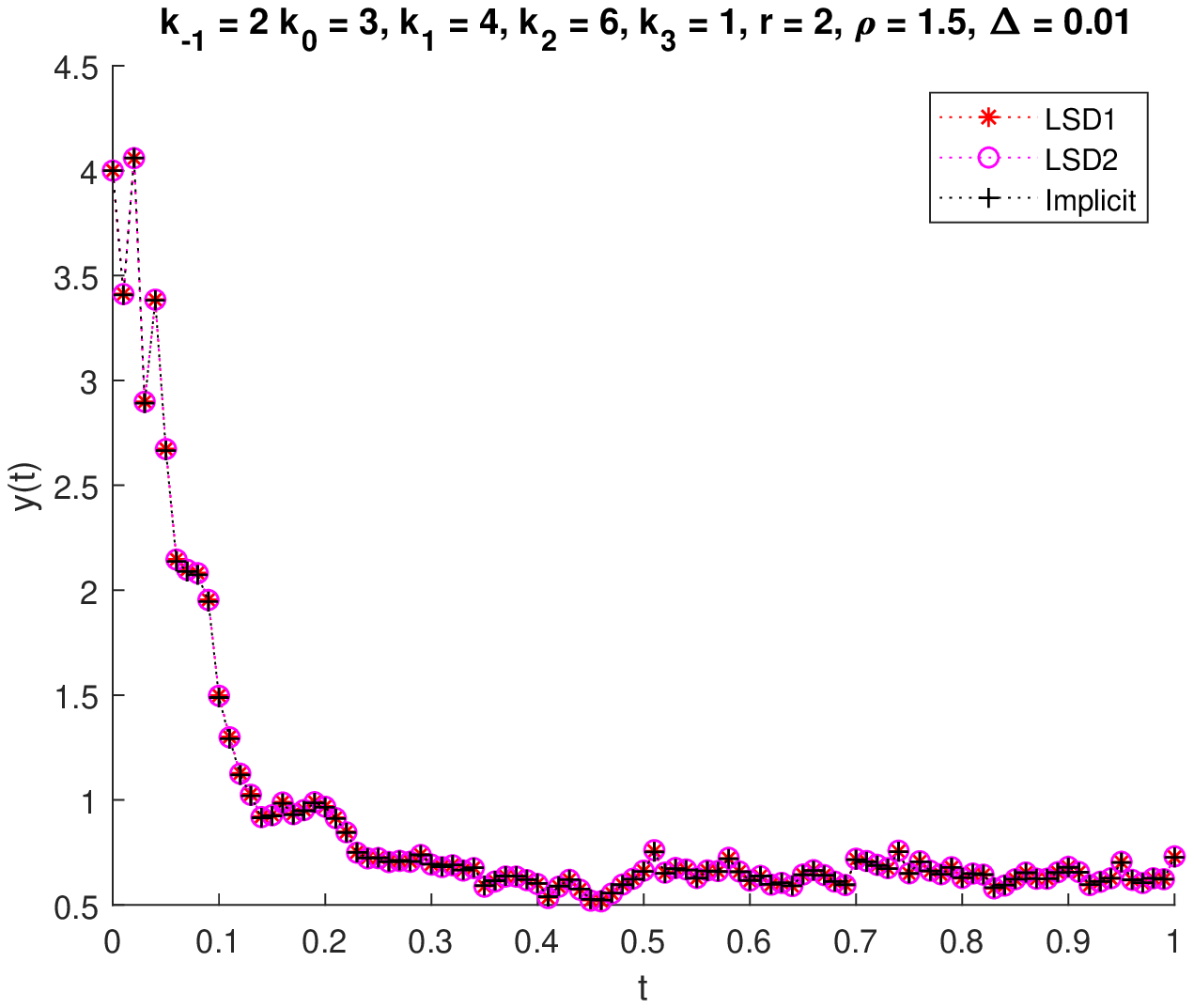}
	\caption{A path with $\D = 10^{-2}$.}
\end{subfigure}
\begin{subfigure}{.47\textwidth}
	\includegraphics[width=1\textwidth]{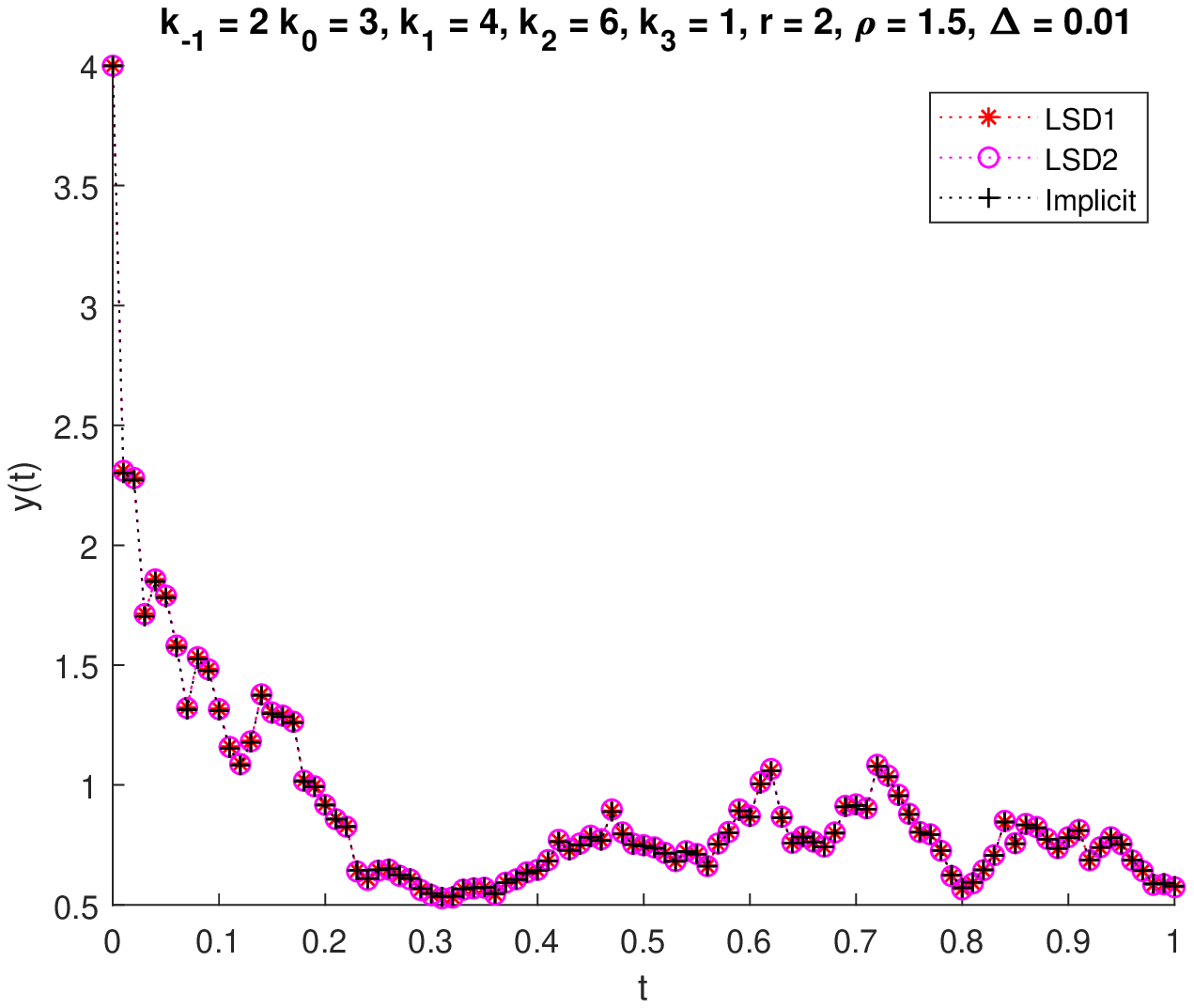}
	\caption{A path with $\D = 10^{-2}$.}
\end{subfigure}
	\caption{Trajectories  of  (\ref{LSD-eq:SDschemeAitLToriginal}), (\ref{LSD-eq:SDschemeAitLToriginal2}) and (\ref{LSD-eq:ImplicitschemeAIT}) for the approximation of (\ref{LSD-eq:Ait}) with various $\D$.}\label{LSD-fig:LSDAIT}
\end{figure}

We use a set of parameters so that (\ref{LSD-eq:ImplicitschemeAIT}) works; we take the coefficients $k_{-1} = 2, k_0 =3, k_1 = 4, k_2 =6, k_3 = 1$ the exponents $r = 2$ and $\rho = 3/2$ and $T =1.$ We compare the proposed versions of LSD schemes (\ref{LSD-eq:SDschemeAitLToriginal}) and  (\ref{LSD-eq:SDschemeAitLToriginal2}) with the implicit method (\ref{LSD-eq:ImplicitschemeAIT}). Figure \ref{LSD-fig:LSDAIT} shows that the LSD1 and LSD2 are very close to the implicit method. We give a presentation of the difference of the two methods in Figure \ref{LSD-fig:LSDsminSDHT}.

\begin{figure}[ht]
	\centering
		\centering
			\begin{subfigure}{.47\textwidth}
			\includegraphics[width=1\textwidth]{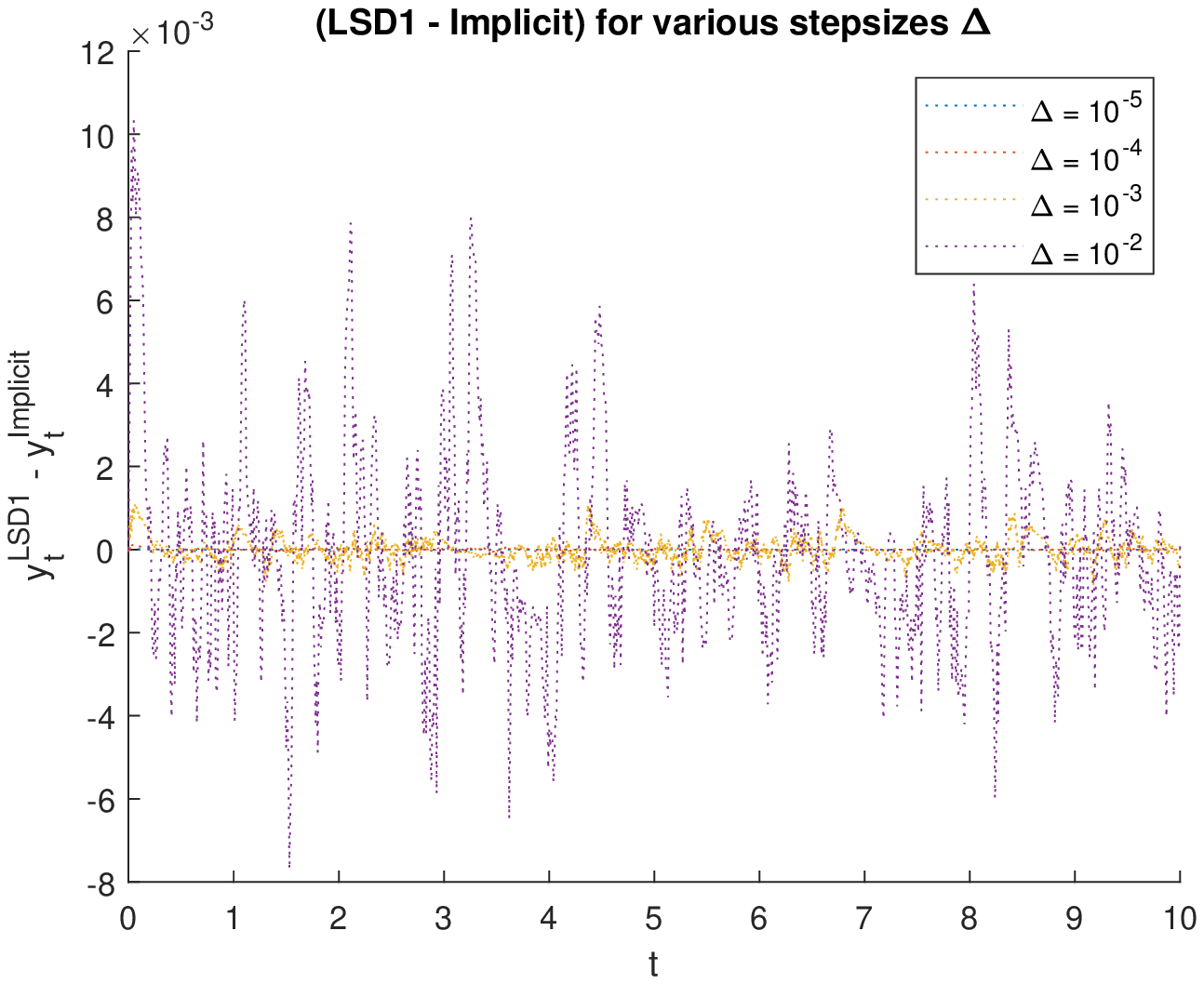}
			\caption{(\ref{LSD-eq:SDschemeAitLToriginal}) - (\ref{LSD-eq:ImplicitschemeAIT})}
		\end{subfigure}
		\begin{subfigure}{.47\textwidth}
			\includegraphics[width=1\textwidth]{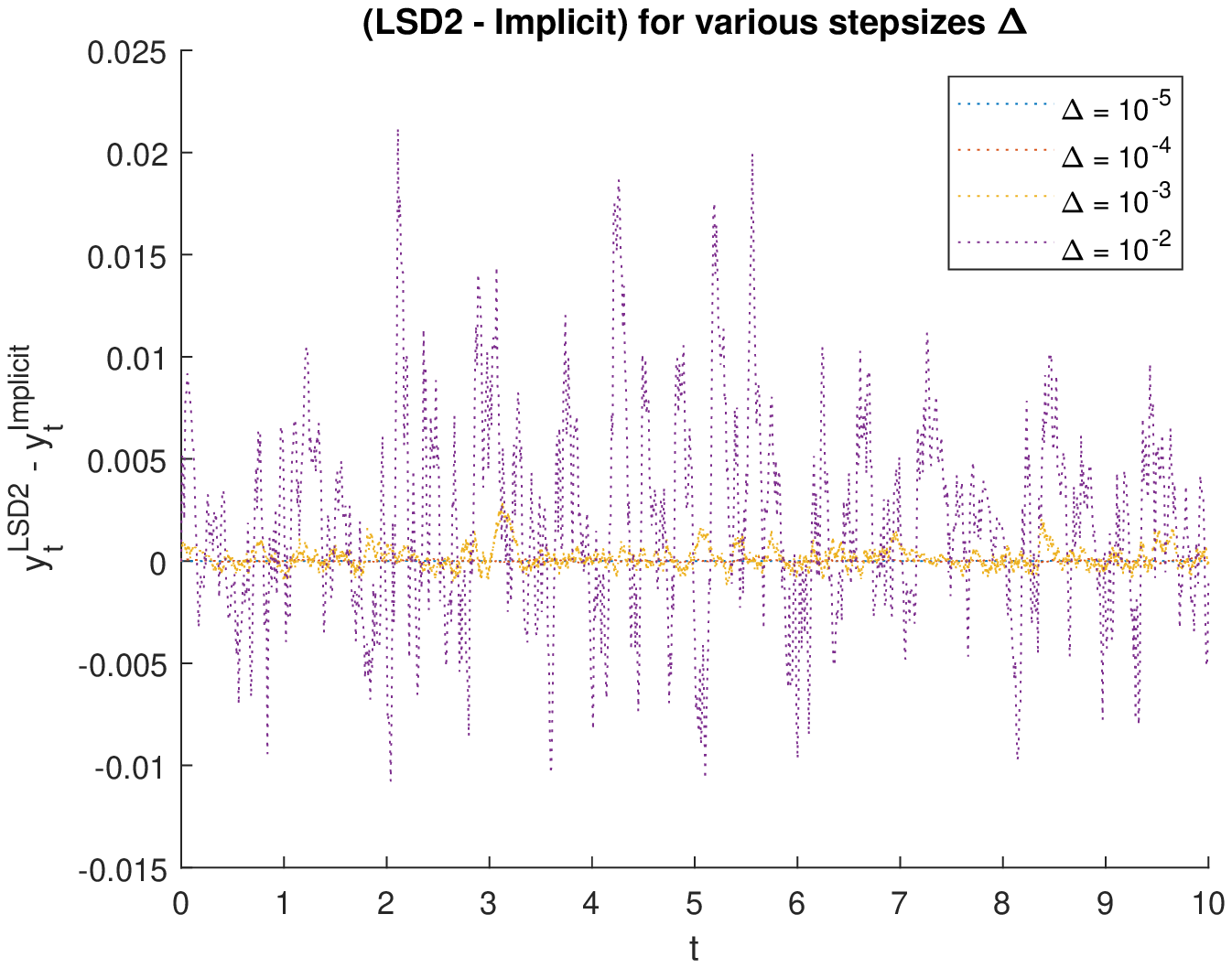}
			\caption{(\ref{LSD-eq:SDschemeAitLToriginal2}) - (\ref{LSD-eq:ImplicitschemeAIT})}
		\end{subfigure}	
		\begin{subfigure}{.55\textwidth}
		\includegraphics[width=1\textwidth]{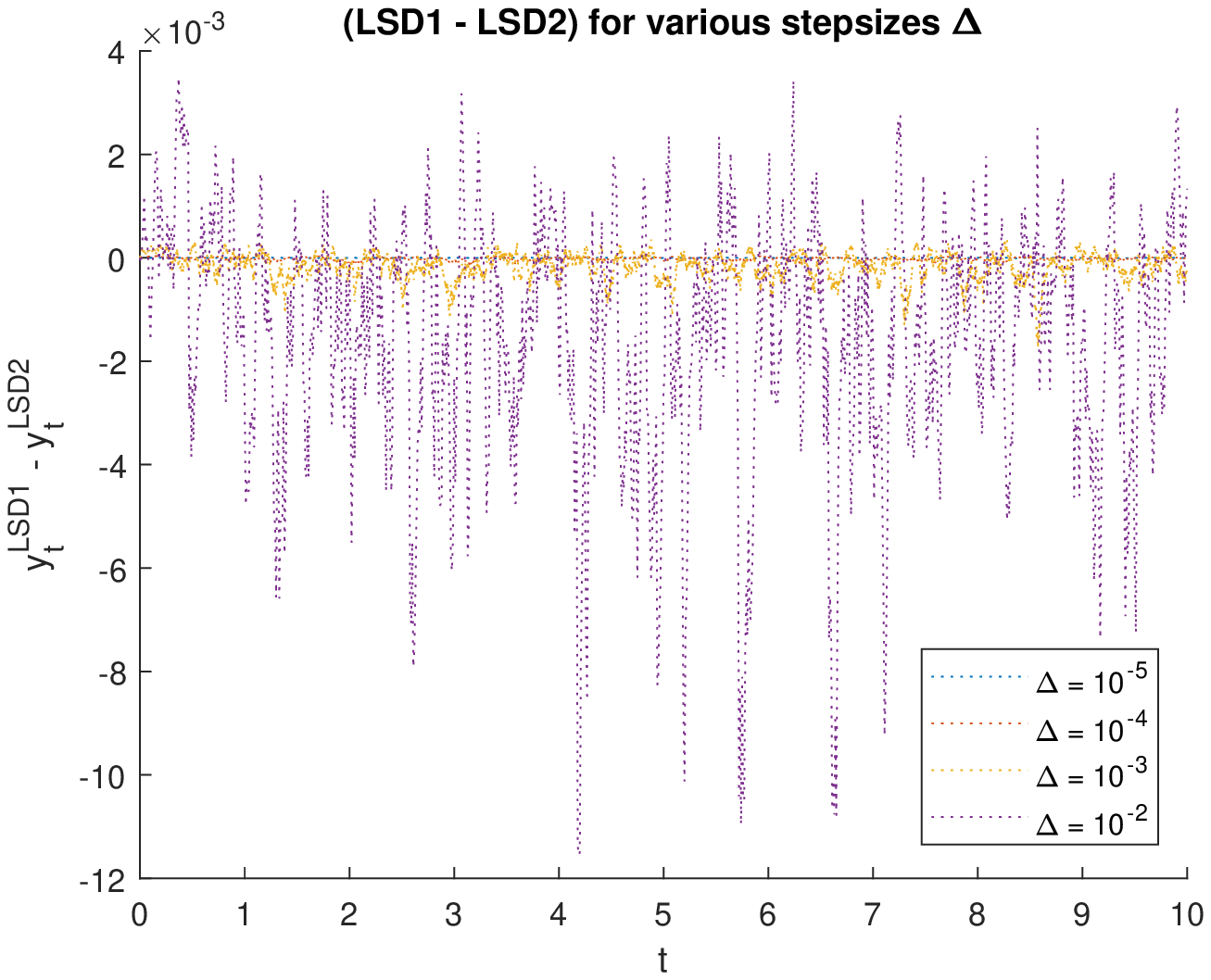}
		\caption{(\ref{LSD-eq:SDschemeAitLToriginal}) - (\ref{LSD-eq:ImplicitschemeAIT})}
	\end{subfigure}	
		\caption{Differences between (\ref{LSD-eq:SDschemeAitLToriginal}), (\ref{LSD-eq:SDschemeAitLToriginal2}) and (\ref{LSD-eq:ImplicitschemeAIT}) for the approximation of (\ref{LSD-eq:Ait}) with various step-sizes.}\label{LSD-fig:LSDsminImplAIT}
\end{figure}

\begin{figure}[ht]
	\centering

			\begin{subfigure}{.47\textwidth}
			\includegraphics[width=1\textwidth]{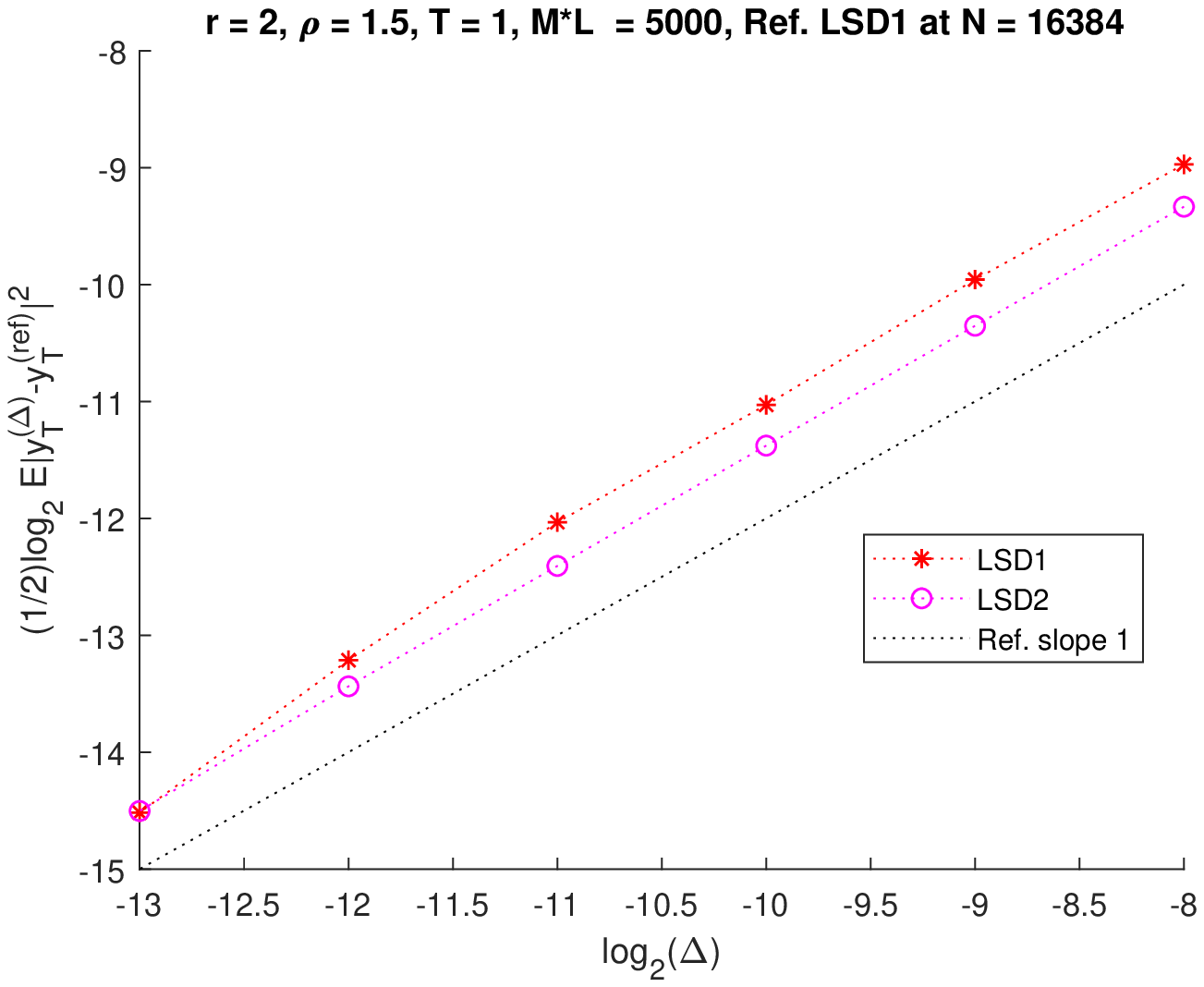}
			\caption{LSD1 as reference solution}
		\end{subfigure}
		\begin{subfigure}{.47\textwidth}
			\includegraphics[width=1\textwidth]{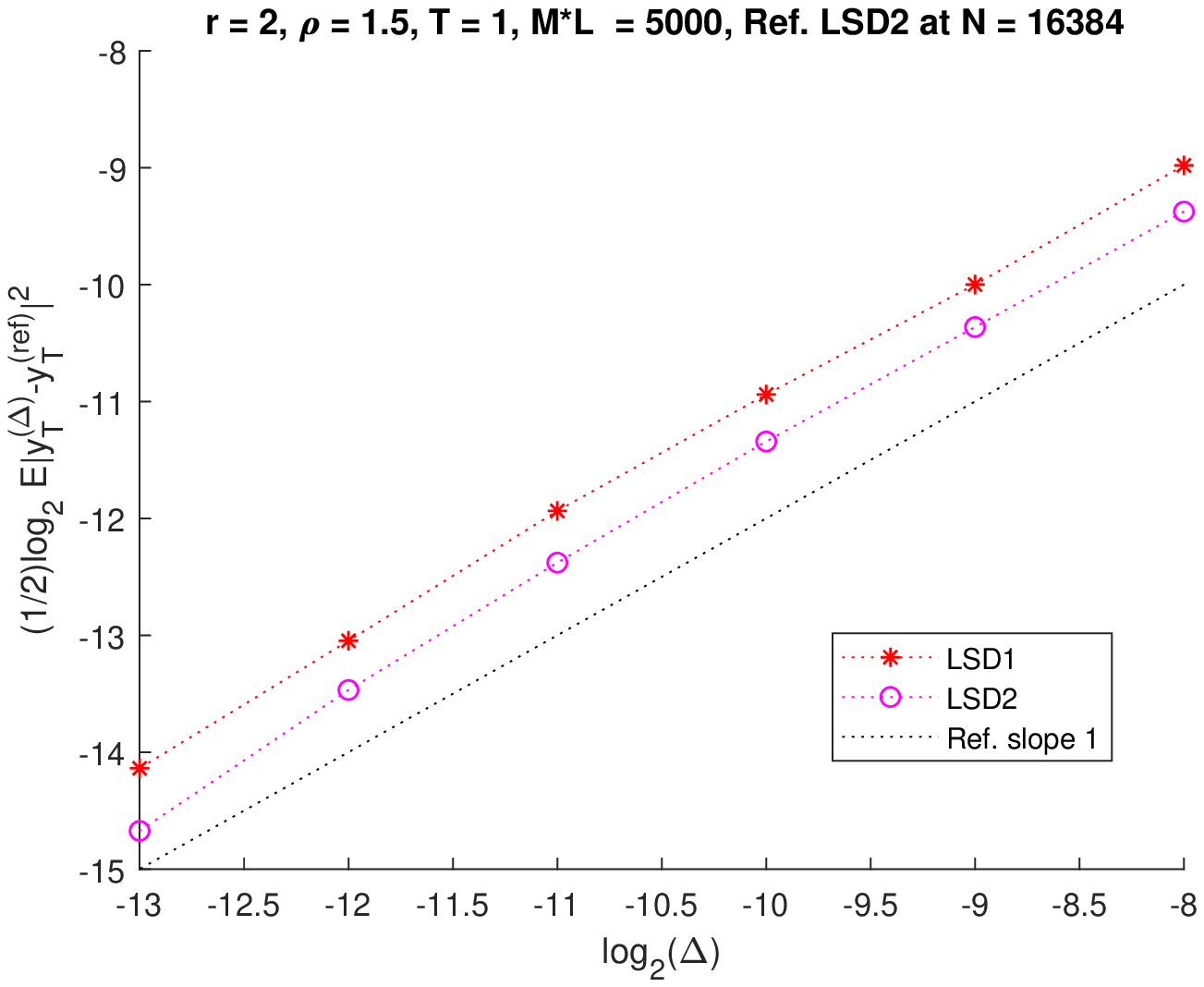}
			\caption{LSD2 as reference solution}
		\end{subfigure}	
		\caption{Convergence of  LSD methods (\ref{LSD-eq:SDschemeAitLToriginal} and (\ref{LSD-eq:SDschemeAitLToriginal2}) for the approximation of (\ref{LSD-eq:Ait}) with different reference solutions.}\label{LSD-fig:LSDorderAIT}
\end{figure}

Finally, we examine numerically the order of strong convergence of the LSD method. The numerical results suggest that the LSD1 and LSD2 schemes converge in the mean-square sense with order close to $1,$ see Figure \ref{LSD-fig:LSDorderAIT}. 

\bibliographystyle{unsrt}\baselineskip12pt
\bibliography{LSD}

\appendix

\section{Lamperti Tranformation of (\ref{LSD-eq:exampleSDE}), (\ref{LSD-eq:Heston}), (\ref{LSD-eq:Ait})}	\label{LSD-ap:Lamperti_tranformation}

Applying the It\^o formula to the transformation $z(x) = \frac{2}{k_3}\sqrt{x}$ of (\ref{LSD-eq:exampleSDE}) we obtain 

\beao
dz_t & = & \left( \frac{1}{k_3}(x_t)^{-1/2}(k_1 - k_2x_t)  + \frac{1}{k_3}\frac{1}{2}(-\frac{1}{2})(x_t)^{-3/2}(k_3)^2(x_t)\right)dt +  \frac{1}{k_3}(x_t)^{-1/2}k_3(x_t)^{1/2}dW_t\\
& = &  \left( \frac{k_1}{k_3}(x_t)^{-1/2} - (\frac{k_2}{k_3} + \frac{1}{4}k_3)\sqrt{x_t}\right)dt + dW_t\\
& = &  \left( \frac{2k_1}{(k_3)^2}(z_t)^{-1} - (\frac{k_2}{2} + \frac{(k_3)^2}{8})z_t\right)dt + dW_t,
\eeao
or for $t\geq t_0$
\beao
z_t & = & z_{t_0} + \int_{t_0}^t \left( \frac{2k_1}{(k_3)^2}(z_s)^{-1} - (\frac{k_2}{2} + \frac{(k_3)^2}{8})z_s\right)ds + \int_{t_0}^t dW_s\\
 & = & z_{t_0} + \int_{t_0}^t \left( \frac{2k_1}{(k_3)^2}(z_s)^{-1} - (\frac{k_2}{2} + \frac{(k_3)^2}{8})z_s\right)ds  + W_t - W_{t_0}.
\eeao

Analogously, the transformation $z(x) = \frac{2}{k_3}(x)^{-1/2}$ of (\ref{LSD-eq:Heston}) has the following dynamics, 

\beao
dz_t & = & \left( \frac{-1}{k_3}(x_t)^{-3/2}(k_1x_t - k_2(x_t)^2)  + \frac{3}{4k_3}(x_t)^{-5/2}(k_3)^2(x_t)^3\right)dt +  \frac{-1}{k_3}(x_t)^{-3/2}k_3(x_t)^{3/2}dW_t\\
& = &  \left(-\frac{k_1}{k_3}(x_t)^{-1/2} + (\frac{k_2}{k_3} + \frac{3}{2}k_3)\sqrt{x_t}\right)dt + dW_t\\
& = &  \left( (\frac{2k_2}{(k_3)^2} + 3)(z_t)^{-1} - \frac{k_1}{2}z_t\right)dt + dW_t,
\eeao
or for $t\geq t_0,$
$$
z_t  = z_{t_0} + \int_{t_0}^t \left( (\frac{2k_2}{(k_3)^2} + 3)(z_s)^{-1} - \frac{k_1}{2}z_s\right)ds  + W_t - W_{t_0}.
$$

Finally, the transformation  $z(x) = x^{1-\rho}$ of (\ref{LSD-eq:Ait}) is such that

\beao  
dz_t &=& \left( (1-\rho)(x_t)^{-\rho} (k_{-1}(x_t)^{-1} -k_0 + k_1x_t - k_2 (x_t)^r) - \frac{\rho(1-\rho)(k_3)^2}{2}(x_t)^{\rho+1} \right)dt\\
&& + k_3(1-\rho)(x_t)^{-\rho}(x_t)^{\rho}dW_t\\
&=& (1-\rho)\left( k_{-1}(x_t)^{-\rho-1} - k_0(x_t)^{-\rho} + k_1(x_t)^{-\rho+1} - k_2 (x_t)^{-\rho+r} - \frac{\rho(1-\rho)(k_3)^2}{2}(x_t)^{\rho+1} \right)dt\\
&& + k_3(1-\rho)dW_t\\
&=& (1-\rho)\left( k_{-1}(z_t)^{\frac{\rho+1}{\rho-1}} - k_0(z_t)^{\frac{\rho}{\rho-1}} + k_1z_t - k_2 (z_t)^{-\frac{r-\rho}{\rho-1}} - \frac{\rho(1-\rho)(k_3)^2}{2}(z_t)^{-1} \right)dt\\
&& + k_3(1-\rho)dW_t,
\eeao 
or for $t\geq t_0,$

\beao 
z_t &=& z_{t_0} + \int_{t_0}^t \Big(k_{-1}(1-\rho)(z_s)^{\frac{\rho+1}{\rho-1}} -k_{0}(1-\rho)(z_s)^{\frac{\rho}{\rho-1}} + k_{1}(1-\rho)z_s\\
&& -k_{2}(1-\rho)(z_s)^{\frac{\rho - r}{\rho-1}} - \frac{\rho(1-\rho)(k_3)^2}{2}(z_s)^{-1} \Big)ds + k_3(1-\rho)(W_t - W_{t_0}).
\eeao

\section{Solution of Bernoulli equations (\ref{LSD-eq:SD schemeExampleLT}), (\ref{LSD-eq:SD schemeExampleLT2}), (\ref{LSD-eq:SD schemeCEVLT})}\label{LSD-ap:Bernoulli_sol}

Consider the following differential equation
\beqq\label{LSD-eq:SD Bernoulli_de}
y_t = A_n  + \int_{t_n}^t \left(B_n(y_s)^{-l} + C_ny_s\right) ds,
\eeqq
with $l>0.$ The dynamics for the transformation $r = y^{1+l}$ are
$$dr_t = \left((1+l)B_n + (1+l)C_nr_t\right) dt,
$$
that is a linear equation with solution
\beao
r_t &=& \frac{\int_{t_n}^t e^{-(1+l)C_n(s-t_n)}(1+l)B_n ds + (A_n)^{1+l}}{e^{-(1+l)C_n(t-t_n)}}\\ 
&=& (1+l)B_ne^{(1+l)C_nt}\int_{t_n}^t e^{-(1+l)C_ns}ds + (A_n)^{1+l}e^{(1+l)C_n(t-t_n)}\\
&=& -\frac{(1+l)B_n}{(1+l)C_n}e^{(1+l)C_nt}(e^{-(1+l)C_nt} - e^{-(1+l)C_nt_n}) + (A_n)^{(1+l)}e^{(1+l)C_n(t-t_n)}\\
&=& -\frac{B_n}{C_n}(1 - e^{(1+l)C_n(t-t_n)}) + (A_n)^{1+l}e^{(1+l)C_n(t-t_n)},
\eeao
where for the case $C_n = 0 $ we read
$$
r_t = (1+l){B_n}(t-t_n) + (A_n)^{1+l}.
$$ 

\section{Solution of the equation (\ref{LSD-eq:WFmodelLT})}\label{LSD-ap:WFmodelLT_sol}

We rewrite equation 
(\ref{LSD-eq:WFmodelLT}) as
$$
\frac{1}{\cot (y/2)}dy = A_ndt
$$
and integrate between $[t_n, t]$ to get
\beao
\int_{t_n}^t \tan(y/2)d(y/2) &=& \frac{A_n}{2}(t-t_n)\\
-\ln|\cos (y_t/2)|& = & \frac{A_n}{2}(t-t_n)  + C\\
|\cos (y_t/2)|& = & Ce^{-\frac{A_n}{2}(t -t_n)}.
\eeao

\end{document}